\documentclass[12pt]{article}
\usepackage{amssymb}
\usepackage{latexsym}
\usepackage{amsfonts}
\usepackage{amsmath}
\newtheorem{theorem}{Theorem}[section]
\newtheorem{lemma}[theorem]{Lemma}

\newtheorem{proposition}[theorem]{Proposition}

\newtheorem{corollary}[theorem]{Corollary}
\newtheorem{remark}[theorem]{Remark}

\numberwithin{equation}{section}

\usepackage[dvips]{color}
\usepackage{mathtools}
\mathtoolsset{showonlyrefs}
\usepackage{mathrsfs}
\usepackage{comment}
\usepackage{hyperref}
\usepackage{xcolor}

\def\P{\mathbb{P} }
\def\R{\mathbb{R} }

\def\bE{{\bf {E} }}
\def\E{\mathbb{E} }

\def\bP{{\bf P} }

\def\cL{{\cal L} }

\def\cP{\mathcal{P}}

\begin{document}
	\begin{doublespace}
	\title{Moments of additive martingales  of branching L\'evy processes and applications}
		\date{}
	\author{{\bf Yan-Xia Ren}\thanks{Y.-X. Ren:  LMAM School of Mathematical Sciences \& Center for
		Statistical Science, Peking
		University,  Beijing, 100871, P.R. China.
		Email: {\texttt yxren@math.pku.edu.cn}}
	\quad
	{\bf Renming Song}\thanks{R. Song: Department of Mathematics,
		University of Illinois,
		Urbana, IL 61801, U.S.A.
		Email: {\texttt rsong@illinois.edu}}
	\quad
	{\bf Rui Zhang*}\thanks{R. Zhang*:  corresponding author. School of Mathematical Sciences \& Academy for Multidisciplinary Studies, Capital Normal
		University,  Beijing, 100048, P.R. China.
		Email: {\texttt zhangrui27@cnu.edu.cn}}}		
		
		\maketitle
		
		\abstract{Let $W_t(\theta)$ be
		the Biggins martingale of
		a supercritical branching L\'evy process with non-local branching mechanism,  and denote by  $W_\infty(\theta)$ its limit.
		In this paper, we first study moment properties of $W_t(\theta)$ and $W_\infty(\theta)$, and the tail behavior of $W_\infty(\theta)$. We then apply these results to establish
  central limit theorems for $W_t(\theta)-W_\infty(\theta)$.
		}
		\section{Introduction}

		Suppose that $\beta>0$ is a constant. Suppose that on some probability space $(\Omega,\mathcal{F}, \bP)$, we are given three independent quantities: an exponential random variable $\tau$ of parameter $\beta$, a L\'evy process $\xi=(\xi_t, t\ge0)$ on $\R$
		starting from the origin
		and a point process  $\mathcal{P}=\sum_{i=1}^N \delta_{S_i}$ on $\R$.
		We emphasize that $N$ is allowed to be $\infty$.
		We will use $\bE$ to denote the expectation with respect to $\bP$.
		For $x\in \R$, we will use $\bP_x$ to denote the law of $\xi+x=(\xi_t+x, t\ge0)$
		and use $\bE_x$ to denote the corresponding expectation.

		In this paper, we consider the following branching system on $\R$:
		(i) initially there is an individual $u$ located at some point $x\in \R$ and it moves according to the L\'evy process $x+\xi$;
		(ii) after an exponential time, independent of the spatial motion, of parameter $\beta>0$, this individual dies and gives birth to a random number of individuals, and the offspring are displaced, relative to the death location of $u$,  according to a point process $\mathcal{P}^u=\sum_{i=1}^{N^u} \delta_{S_i^u}$ which has the same law as $\mathcal{P}$;
		(iii) any individual, once born, evolves independently from its birth place with the
			same branching rule
		and the same spatial motion as its parent.

	Let $\mathcal{M}_a(\R)$ be the space of integer-valued atomic measures on $\R$,
		and $\mathcal{B}_b(\R)(\mathcal{B}_b^+(\R))$  be the set of bounded (non-negative) Borel functions on $\R$.
		Let $X_t(B)$ be the number of particles alive at time $t$ and located in $B$.
		Then $X=\{X_t,t\geq 0\}$ is an $\mathcal{M}_a(\R)$-valued Markov process,
which is called a
$(\xi, \beta, \mathcal{P})$-branching L\'evy process.
We will use $\{\mathcal{F}_t: t\ge0\}$ to denote the natural filtration of $\{X_t,t\geq 0\}$.
For any $\nu\in\mathcal{M}_a(\R)$, we denote the law of $X$ with initial configuration $\nu$ by $\P_\nu$. We write $\P:=\P_{\delta_0}$ and will use $\E$ to denote the expectation with respect to $\P$.

If  $\bE(N)<\infty$, we  say  that the branching L\'evy process has finite birth intensity.
Recently a more general branching L\'evy process was introduced in \cite{BM}. This more general branching L\'evy process may not have  a first branching time, and can be constructed as increasing limit of a sequence of so-called nested branching L\'evy processes with finite birth intensities.
For the branching L\'evy process of this paper,
we have a first branching time, but $N$ may be infinite.
Therefore, our model is general than  the  classical branching L\'evy process, as introduced by Kyprianou \cite{Kyp99}, where $N$ is finite.

It is easy to see from the spatial homogeneity of $\xi$ that $(X,\P_{\delta_x})$ has the same law as $(x+X,\P)$. Let $\cL_t$ be the set of all the individuals alive at time $t$. For any $u\in \mathcal{L}_s$, we use $z_u(s)$ to denote the position of $u$ at time $s$ and $\{z_u(s)+X_t^u,t\ge 0\}$ to denote the subtree generated by $u$. Then $\{X_t^u,t\ge0\}$ has the same law as $X$ under $\P$. It is well known that the total number of individuals at time $t$ is a continuous time Galton-Watson process. In this paper we always assume that the process is supercritical, that is,  $\bE N>1$, which guarantees that the population survives with positive probability. Let $\mathcal{E}$ be the event of extinction.

		For any $\theta\in\R$, we write $e_\theta(x):=e^{\theta x}$. Define
\begin{equation}\label{def-chi}
\chi(\theta):=\bE\langle e_\theta,\cP\rangle=\bE\Big(\sum_{i=1}^N e^{\theta S_i}\Big) \mbox{ and } \varphi(\theta):=\log \bE(e^{\theta \xi_1}).
\end{equation}
		Let $n(dy)$ be the L\'evy measure of $\xi$. Note that
		$$
		\varphi(\theta)<\infty\Longleftrightarrow \int_{-\infty}^\infty {\bf 1}_{|y|>1}e^{\theta y}n(dy)<\infty.$$
		For any $\theta$ with $\varphi(\theta)<\infty$, $\varphi(\theta)$ can be written as:
		$$
		\varphi(\theta)=a\theta+\frac{1}{2} b^2\theta^2+\int_{\R}\left( e^{\theta y}-1-\theta y{\bf 1}_{|y|\le 1}\right)\, n(dy),
		$$
		where $a\in\R, b\ge0.$
		
		Define
\begin{equation}\label{def-kappa}
\kappa(\theta):=\beta(\chi(\theta)-1)+\varphi(\theta)\in(-\infty,\infty].
\end{equation}
		Let $\Theta:=\{\theta>0: \kappa(\theta)<\infty\}$.
		In this paper we always assume that $\Theta$ is non-empty.
		It is easy to see from their definitions that $\chi(\theta)$, $\varphi(\theta)$ and $\kappa(\theta)$ are all convex in $\Theta$.
		Let $\Theta_0$ be the interior of $\Theta$. It is clear that $\Theta_0$ is an open interval and we denote it by $(\theta_-,\theta_+)$.

The following elementary result will be used to introduce the Biggins martingale below. We will give
the proof of this result in the Appendix.
\begin{lemma}\label{f-moment}
Let $f(x)=|x|^re^{\theta x}$ or $f(x)=(x\vee 0)^r e^{\theta x}$ with $r\ge0$ and $\theta>0$. Then $\E(\langle f,X_t\rangle)<\infty$ for some $t>0$ if and only if $\kappa(\theta)<\infty$,
\begin{align}\label{condition-f}
\int_{|y|>1} f(y)n(dy)<\infty \mbox{ and } \bE\left(\sum_{i=1}^N f(S_i)\right)<\infty.
\end{align}
In particular, if $\E(\langle f,X_t\rangle)<\infty$ for some $t>0$, then $\E(\langle f,X_t\rangle)<\infty$ for all $t>0$.			
\end{lemma}

		For any $\theta\in \Theta$, define
		$$
		W_t(\theta):=e^{-\kappa(\theta)t}\langle e_\theta,X_t\rangle.
		$$

		\begin{lemma}
			If $\theta\in \Theta$, then
			$\{W_t(\theta), t\ge0\}$ is a non-negative martingale of mean 1.
		\end{lemma}
		{\bf Proof:}
		Applying  Lemma \ref{f-moment} with $r=0$,  we get
		 that $\E \langle e_\theta,X_t\rangle<\infty$ for any $\theta\in\Theta$ and $t\ge0$.
		Let $\tau_1$ be the first splitting time. By the Markov property and the branching property of $X$, we have that
		\begin{align*}
			\E \langle e_\theta,X_t\rangle&=\E [\langle e_\theta,X_t\rangle; \tau_1>t]+\E[ \langle e_\theta,X_t\rangle; \tau_1\le t]\nonumber\\
			&=e^{-\beta t}\bE(e^{\theta \xi_t})+\int_0^t \beta e^{-\beta s} \bE\left(\sum_{i=1}^N e^{\theta(\xi_s+S_i )}\right)\E \langle e_\theta,X_{t-s}\rangle\, ds\nonumber\\
			&=e^{(\varphi(\theta)-\beta)t}+\chi(\theta)\int_0^t \beta e^{-\beta s} e^{\varphi(\theta) s}\E \langle e_\theta,X_{t-s}\rangle\, ds\\
			&=e^{(\varphi(\theta)-\beta)t}+\chi(\theta)e^{(\varphi(\theta)-\beta)t}\int_0^t \beta  e^{(\beta-\varphi(\theta) )s}\E \langle e_\theta,X_{s}
						\rangle\, ds.
		\end{align*}
		Let $g(t)=\E \langle e_\theta,X_t\rangle$. It is easy to get that $$g'(t)=(\beta(\chi(\theta)-1)+\varphi(\theta))g(t)=\kappa(\theta) g(t),$$
		which implies that
		$$\E \langle e_\theta,X_t\rangle=g(t)=e^{\kappa(\theta)t}.$$
		By the Markov property and the branching property of $X$, we have that, for any $t>s>0$,
		\begin{align}
			\E (\langle e_\theta,X_t\rangle|\mathcal{F}_s)=\sum_{u\in\mathcal{L}_s} e^{\theta z_u(s)}	\E \langle e_\theta,X_{t-s}\rangle=e^{\kappa(\theta)(t-s)}\langle e_\theta,X_s\rangle,
		\end{align}
		which implies that $\{W_t(\theta), t\ge0\}$ is a non-negative martingale .
		\hfill$\Box$

		\bigskip

		The martingale $\{W_t(\theta): t\ge0\}$ is called the Biggins martingale or additive martingale of our branching L\'evy process.
		Since $\{W_t(\theta): t\ge0\}$ is a non-negative martingale, $W_\infty(\theta):=\lim_{t\to\infty}W_t(\theta)$ exists almost surely.

		For the general model in \cite{BM}, necessary and sufficient conditions for the non-degeneracy of the counterpart of $W_\infty(\theta)$ were given in \cite{BM18,IB19}.  We now specialize their conditions to our model.
		We write
		\begin{align}\label{def-kappa'}
			\kappa'(\theta):=\beta \bE\Big(\sum_{i=1}^N S_ie^{\theta S_i}\Big)+a+b^2\theta+\int_\R y(e^{\theta y}-1_{|y|\le 1})n(dy),
		\end{align}
		whenever
		\begin{align}\label{kappa'}
			\bE\Big(\sum_{i=1}^N |S_i|e^{\theta S_i}\Big)<\infty \mbox{ and } \int_{|y|>1} |y|e^{\theta y}n(dy)<\infty.
		\end{align}
		It was proved in \cite{BM18} that, for any for $\theta\in\Theta$ satisfying \eqref{kappa'}, $\{W_t(\theta),t\ge 0\}$ is uniformly integrable if and only if
		\begin{align}\label{LLogL}
			\theta\kappa'(\theta)<\kappa(\theta) \mbox{ and }\bE\left[ \langle e_\theta,\cP\rangle\log_+\langle e_\theta,\cP\rangle\right]<\infty,
		\end{align}
		where $\log_+x:=\max\{\log x,0\}$.
		When $\theta\in\Theta$ does not satisfy \eqref{kappa'},
		Iksanov and Mallein \cite[Theorem 4.1]{IB19} provided a more general  necessary and sufficient condition for the non-degeneracy of $W_\infty(\theta)$.
		
		The purpose of this paper is to investigate
		moment properties of $W_t(\theta)$ and $W_\infty(\theta)$, and the tail behaviors of  $W_\infty(\theta)$,
		and then use these results to study the fluctuations of $W_t(\theta)-W_\infty(\theta)$.
		\bigskip

Now we give a summary of the main results of this paper. Throughout this paper,  we always assume that $\theta\in\Theta$.

\begin{description}
\item{(a)} {\bf Moment properties of $W_t(\theta)$}. In Theorem  \ref{wtp}, we give necessary and
sufficient conditions, in terms of $(\beta,\xi,\mathcal{P})$, for the finiteness of $W_t(\theta)^p$, $p>1$. In Theorems \ref{prop:moment7} and \ref{proposition-log}, we give sufficient conditions for the finiteness of $\E(W_t(\theta)H(W_t(\theta)))$, where $H$ is a slowly varying function at $\infty$.

\item{(b)} {\bf Moment properties of $W_\infty(\theta)$}. In Theorem \ref{main-Lp}, we give necessary and
sufficient conditions for the finiteness of $W_\infty(\theta)^p$, $p>1$. In Theorem  \ref{Theorem:WH(W)} we provide sufficient
conditions for the finiteness of  $\E(W_\infty(\theta)^pH(W_\infty(\theta)))$,  where $H$ is a slowly varying function at $\infty$.

\item{(c)} {\bf Tail behaviors of $W_\infty(\theta)$}. In the case when there exists $p_*>1$ such that $\kappa(\theta p_*)=p_*\kappa(\theta)$, we give sufficient conditions for $\P(W_\infty(\theta)>x)  \sim cx^{-p^*}$ as $x\to \infty$ in
Proposition \ref{boundry case}.  In Theorem \ref{Theorem:tailW} we give necessary and sufficient conditions for
$$
\P(W_\infty(\theta)>x)  \sim x^{-p}L(x),\quad x\to\infty,
$$
where $L$ is a lowly varying function at $\infty$, in the case when $p\in (1, 2)$ satisfies $p\theta\in \Theta_0$ and $\kappa(p\theta)<p\kappa(\theta)$.

\item{(d)}
{\bf Extremal process of $X$}. Propositions \ref{prop:Maximum}, \ref{maximum-weak} and \ref{prop:extremal} are results on the extreme  process of $X$.

\item{(e)} {\bf Central limit-type theorems for $W_t(\theta)-W_\infty(\theta)$}. Theorem \ref{them-clt2}
and its corollaries say that, if $\kappa(2\theta)<2\kappa(\theta)$ and $\bE(\langle e_{\theta},\cP\rangle^2)<\infty$, then after appropriate normalization, conditioned on $\mathcal{F}_t, $ $W_t(\theta)-W_\infty(\theta)$ converges weakly to a normal random variable.
Theorem \ref{main-stableCLT}
and its corollaries say that, when the tail distribution of $W_\infty(\theta)$ belongs to the domain of attraction of a $p$-stable distribution, $p\in (1, 2)$, under some other conditions, after appropriate  normalization, conditioned on $\mathcal{F}_t$, $W_t(\theta)-W_\infty(\theta)$ converges weakly to a $p$-stable random variable.
\end{description}

Set
\begin{align}\label{kappa''}
	\kappa''(\theta):=\beta \bE\Big(\sum_{i=1}^N (S_i)^2e^{\theta S_i}\Big)+b^2+\int_{\R} y^2 e^{\theta y} n(dy)\in[0,\infty].
\end{align}
Note that if $\theta\in\Theta_0$,
then \eqref{kappa'} holds,
$\kappa'(\theta)$ and $\kappa''(\theta)$ are finite, and they are
the first  and second derivatives of $\kappa(s)$ respectively.
 Let $\Theta':=\Theta\cap\{\theta>0: \eqref{kappa'} \mbox{ holds} \}$.
Clearly, the interior of $\Theta'$ is also $\Theta_0=(\theta_-,\theta_+).$

If there exists $\theta>0$ such that $\kappa''(\theta)=0$, then $\mathcal{P}=N\delta_0$ and $\xi_t=at$.
Therefore, for all
$\theta\in\R$, $\kappa(\theta)=\beta(\bE(N)-1)+a\theta$. Thus, $\Theta$ is nonempty if and only if $\bE(N)<\infty$.
In this case, $X_t=\|X_t\|\delta_{at},$
thus $W_t(\theta)=e^{-\beta(\bE(N)-1)t}\|X_t\|$ for all $\theta\in\R$.
Note that $\theta\kappa'(\theta)<\kappa(\theta)$, $\kappa(p\theta)<p\kappa(\theta)$  for all $\theta>0$ and $p>1$,
and that \eqref{LLogL} is equivalent to
$\bE(N\log_+N)<\infty$.

If $\kappa''(\theta)>0$ for all $\theta>0$. Then $\theta\kappa'(\theta)-\kappa(\theta)$ is strictly increasing on $\Theta'$.
If, in addition,  there exists $\theta_*\in\Theta'$ such that $\theta_*\kappa'(\theta_*)=\kappa(\theta_*)$,  then
$\theta\kappa'(\theta)<\kappa(\theta)$ for $\theta\in (\theta_-, \theta_*)$ and $\theta\kappa'(\theta)>\kappa(\theta)$ for $\theta\in (\theta_*, \theta_+)$.
We emphasize that there is at most one such $\theta_*$, which is
also called the boundary point in some literature.
Note also that, even when such $\theta_*$ exists, it may not be in $\Theta_0$.
Consider the function $u(\theta)=\kappa(\theta)/\theta$. Note that for $\theta\in\Theta_0$, $u'(\theta)=\frac{\theta\kappa'(\theta)-\kappa(\theta)}{\theta^2}$.
Thus $u(\theta)$ is decreasing on
$(\theta_-, \theta_*)$ and increasing on $(\theta_*, \theta_+)$.
So
\begin{equation}\label{kappa-theta*}
	\frac{\kappa(\theta_*)}{\theta_*}=\inf_{\theta>0}\frac{\kappa(\theta)}{\theta}.
\end{equation}
\begin{remark}\label{rem:kappa}
	(1) For $\theta\in\Theta$ and $p>1$, if  $\kappa(p\theta)<p\kappa(\theta)$, then  $\kappa(r\theta)<r\kappa(\theta)$ for all $1<r\le p$.  If, in addition, $\theta\in\Theta'$ then  $\theta\kappa'(\theta)<\kappa(\theta)$.
	
	(2)	Condition \eqref{LLogL} fails in the boundary case $\theta=\theta_*$, thus
	$W_\infty(\theta_*)=0$, $\P$-a.s.  In Subsection \ref{subsection: derivative martingale}, we will introduce another important martingale in the boundary case.
	
	(3)
	In the case of branching Brownian motion with supercritical local branching mechanism, that is, $\cP=N\delta_0$ with $m:=\bE(N)\in(1,\infty)$, we have
	\begin{equation}\label{local-kappa}
		\kappa(\theta)=\beta\Big(m-1\Big)+\frac{1}{2} \theta^2.
	\end{equation}
	In this case, $\theta_*=\sqrt{2\beta(m-1)}$.  One can easily check that, for $p>1$, the inequality
	$\kappa(p\theta)<p\kappa(\theta)$ holds if and only if  $0<\theta<\theta_*$ and $1<p<2\beta(m-1)\theta^{-2}$. For $0<\theta<\theta_*$,
	$p_*:=2\beta(m-1)\theta^{-2}$ satisfies $\kappa(p_*\theta)=p_*\kappa(\theta)$.
\end{remark}

Now we give some related literature.  The  counterpart of Theorem \ref{main-Lp} for branching random walks was proved in \cite{Liu2000}.  In \cite{BM18}, Bertoin and Mallein gave a sufficient condition for the $L^p$,  $p\in(1,2]$, convergence of the Biggins martingale of the general
branching L\'evy process introduced in \cite{BM}.
Subsequently, \cite[Theorem 4.1]{IB19} gave a sufficient and necessary condition for the $L^p$,  $p\in(1,2]$, convergence of the Biggins martingale of the general branching L\'evy process.
The counterpart of Theorem \ref{Theorem:WH(W)} for Galton-Watson processes was proved in \cite{BD74}. The counterparts of  Proposition \ref{boundry case} and Theorem \ref{Theorem:tailW} for branching random walks were established in  \cite{Liu2000, IKM}. Normal central limit theorems were
proved for Galton-Watson processes  in \cite{Athreya}
and for branching random walks in \cite{IK16, IKM20}.
Stable central limit theorems were proved for Galton-Watson processes in \cite{Heyde71},
			for  branching random walks  it was established in \cite{IKM, IKM20}, and
			for super-Brownian motions can be found in \cite{Yang}.

This paper is organized as follows. In Section 2, we give some preliminaries. Moment properties of $W_t(\theta)$ are studied in Section 3. In Section 4, we first study moment properties of $W_\infty(\theta)$ and then study the tail behaviors of $W_\infty(\theta)$. In Subsections 5.1 and 5.2, we use results of Section 3 to study the derivative martingale  and the extremal process of $X$.
In  Subsection 5.3, as applications of the previously mentioned results, we prove some central limit-type theorems for $W_t(\infty)-W_\infty(\theta)$. In the Appendix, we give the proofs of some auxiliary results.

\section{Preliminaries}

			Let  $\mathbb{L}$  be the class of all the  positive,  locally bounded
			functions on $[0,\infty)$ that are slowly varying at $\infty$.
			For  $L\in\mathbb{L}$ ,  we define
			\begin{equation}\label{def-L*-tildeL}
				L^*(x):=\int_{1}^{x\vee 1}\frac{L(t)}{t}\,dt \mbox{ and } \tilde{L}(x):=\int_{x\vee 1}^{\infty}\frac{L(t)}{t}\,dt, \quad x\ge 0,
			\end{equation}
			whenever $\int_{1}^{\infty}\frac{L(t)}{t}\,dt<\infty.$
			It is well known that $L^*$ and $\tilde{L}$ (when $\tilde{L}$ is finite) are slowly varying and
			\begin{align}\label{compare}
				\frac{L(x)}{L^*(x)}\to0,\quad \frac{L(x)}{\tilde{L}(x)}\to0, \quad  x\to\infty.
			\end{align}
			By \eqref{compare}, there exists $c>0$ such that for any $x>0,$
			\begin{align}\label{compare2}
				L(x)\le c(1+L^*(x)), \mbox{ and } L(x)\le c(1+\tilde{L}(x)).
			\end{align}
			We put
			$$\mathbb{L}^*:=\{L^*| L\in \mathbb{L}\} \mbox{ and } \widetilde{\mathbb{L}}:=\left\{\tilde{L}| L\in \mathbb{L},\int_{1}^{\infty}\frac{L(t)}{t}\,dt<\infty\right\}.
			$$
			We  now recall some useful properties of slowly vary functions.
			For any $\epsilon>0$, we put
			\begin{align}
				\mathbb{L}_\epsilon:=\{L| &L \in\mathbb{L}, x^\epsilon L(x) \mbox{ is  increasing and } x^{-\epsilon}L(x)  \mbox{ is  decreasing on } [0,\infty)\}.
			\end{align}
			It is clear that  if $L\in \mathbb{L}_\epsilon$, then
			\begin{align}\label{L-epsilon}
				\frac{L(x)}{L (y)}\le \left(\frac{x}{y}\right)^\epsilon\vee \left(\frac{y}{x}\right)^\epsilon, x>0,y>0.
			\end{align}
			For any $L\in\mathbb{L}$ and $\epsilon>0$, it follows from \cite[Theorem 1.5.5]{Bingham} that  there exists $x_\epsilon>0$ such that $x^\epsilon L(x)$ is  increasing and $x^{-\epsilon}L(x)$ is  decreasing on $[x_\epsilon,\infty)$. Define $L_\epsilon(x)=L(x\vee x_\epsilon)$. Then $ L_\epsilon\in \mathbb{L}_\epsilon$. Without loss of generality, in the remainder of this paper, we always assume that $L\in \mathbb{L}_\epsilon$ for some $\epsilon>0$.

			For any integer $n\ge0$, define
			\begin{align}\label{Tn} \mathcal{T}_n(x):=(-1)^{n+1}\Big(e^{-x}-\sum_{k=0}^{n}\frac{(-1)^kx^k}{k!}\Big),\quad x\in\R.
				\end{align}
			Using Taylor's formula, it is easy to check that for any $n\ge0$ and $x\ge0$,
			\begin{equation}\label{est:T}
				0\le \mathcal{T}_n(x)\le \frac{2x^n}{n!}\wedge\frac{x^{n+1}}{(n+1)!}.
			\end{equation}
			For any  nonnegative random variable $\mathcal{Y}$ with
			$E(\mathcal{Y}^n)<\infty$, we have  by \eqref{Tn} that for all $\lambda\ge 0$,
			$$E(\mathcal{T}_n(\lambda \mathcal{Y}))=(-1)^{n+1}\Big(
			E(e^{-\lambda \mathcal{Y}})
			-\sum_{k=0}^{n}\frac{(-1)^kE(\mathcal{Y}^k)}{k!}\lambda^k\Big)
			<\infty.
			$$

The following two results are \cite[Theorems A and Theorem B]{BD74}. They will play key
roles in the proofs of
Theorems \ref{prop:moment7}, \ref{proposition-log},  \ref{Theorem:WH(W)}  and \ref{Theorem:tailW} below.
			
			\begin{lemma}\label{lemma:tb}
				Let $\mathcal{Y}$ be  a nonnegative random variable and $L\in\mathbb{L}$. For any  integer $n\ge 1$ and $p\in[n,n+1]$, the following statements are equivalent:
				\begin{itemize}
					\item[(1)] $E(\mathcal{Y}^n)<\infty$ and
					\begin{align*}
						\int_0^1 E(\mathcal{T}_n(\lambda \mathcal{Y}))L(1/\lambda)\lambda^{-(1+p)}\,dt<\infty.
					\end{align*}
					\item[(2)]  $E(\mathcal{Y}^pH(\mathcal{Y}))<\infty$,
					where
					\begin{align}\label{def:H}
						H(x)=\left\{\begin{array}{ll}
							L^*(x), & \mbox{ when } p=n; \\
							L(x),  & \mbox{ when } p\in(n,n+1);\\
							\tilde{L}(x),  & \mbox{ when } p=n+1.
						\end{array}\right.
					\end{align}
				\end{itemize}
			\end{lemma}
			
			\begin{lemma}\label{lem:taub}
				Let $\mathcal{Y}$ be a non-negative random variable with $E(\mathcal{Y}^n)<\infty$ for some integer $n\ge0$. For any $p\in(n,n+1)$ and $L\in \mathbb{L}$, the following  statements are equivalent:
				\begin{itemize}
					\item [(1)] $E(\mathcal{T}_n(\lambda \mathcal{Y}))\sim \lambda^{p}L(1/\lambda),\quad \lambda\to0;$
					\item [(2)] $P(\mathcal{Y}>x)\sim \frac{(p-1)(p-2)\cdots(p-n)}{\Gamma(1+n-p)}x^{-p}L(x),\quad x\to\infty.$
				\end{itemize}
		\end{lemma}

	 \section{Moments of $W_t(\theta)$}\label{moments of wt}

\subsection{$p$-th moments of  $W_t(\theta)$ with  $p>1$}

		Let  $n\ge 1$ be an integer.  In our $(\xi, \beta, \mathcal{P})$-branching L\'evy process, if we kill all offspring (along with all their possible descendants) with birth-places at least $n$ units
		to the left of the death-places of their parents, we get a
		$(\xi,\beta,\mathcal{P}^{(n)})$-branching L\'evy process $X^{(n)}$ with
		$$\mathcal{P}^{(n)}=\sum_{i=1}^N{\bf 1}_{\{S_i>-n\}}\delta_{S_i}.$$
		For any $\theta\in\Theta$, we have
		$$\kappa^{(n)}(\theta)=\beta(\bE(\langle e_\theta,\mathcal{P}^{(n)}\rangle-1)+\varphi(\theta)\to\kappa(\theta), \quad n\to\infty.$$
		If $\theta\in\Theta$, then  we have
		$$\bE\left( \sum_{i=1}^N{\bf 1}_{\{S_i>-n\}}\right) \le \bE\left(\sum_{i=1}^N e^{\theta(S_i-n)}\right)<\infty.$$
Thus the $(\xi,\beta,\mathcal{P}^{(n)})$-branching L\'evy process has finite birth intensity.

		\begin{theorem}\label{wtp}
	Let $\theta\in\Theta$ and $p>1$. $\E(W_t(\theta)^p)<\infty$ for some $t>0$  if and only if $\varphi(p\theta)<\infty$ and
	$\E(\langle e_\theta,\mathcal{P}\rangle^p)<\infty.$
	In particular, if $\E(W_t(\theta)^p)<\infty$ for some $t>0$, then $\E(W_t(\theta)^p)<\infty$ for all $t>0$.		
		\end{theorem}
\noindent{\bf Proof :}
(1) We first prove the necessity.
Assume that $\E(W_t(\theta)^p)<\infty$ for some $t>0$.
Since $\{W_s(\theta)^p: s\ge 0\}$ is a submartingale, by the optional stopping  theorem, we have
		$$
		\E(W_{t\wedge \tau_1}(\theta)^p)\le \E(W_t(\theta)^p)<\infty,
		$$ which implies that
		\begin{align*}
			\infty>\E(W_{\tau_1}(\theta)^p;\tau_1<t)&=\bE(e^{-p\kappa(\theta)\tau}e^{p\theta\xi_\tau}\langle e_\theta,\mathcal{P}\rangle^p,\tau<t)\\
			&=\bE(\langle e_\theta,\mathcal{P}\rangle^p)\int_0^t \beta e^{-\beta s}e^{-p\kappa(\theta)s}e^{\varphi(p\theta)s}\,ds.
		\end{align*}
		Thus
		$$
		\varphi(p\theta)<\infty, \mbox{ and } \bE(\langle e_\theta,\mathcal{P}\rangle^p)<\infty.
		$$
		
(2) We now prove the sufficiency in the case $\bE(N)<\infty$.
		Let
		$$
		V_{m}(t):=\E(\langle e_\theta, X_t\rangle^p;t<\tau_m).
		$$
		It is clear that $\E(\langle e_\theta, X_t\rangle^p)=\sup_{m}V_{m}(t). $
		We now derive a recursive formula for $V_m(t)$ and then give an upper bound of $V_m(t)$ by induction.
		Observe  that for all $t>0$,
		$$
		V_1(t)=e^{-\beta t}e^{\varphi(p\theta)t}<\infty.
		$$
		By the strong Markov property of $X_t$, we have
		\begin{align}
			&V_{m+1}(t) =\E (\langle e_\theta,X_t\rangle^{p};t<\tau_1)+\E (\langle e_\theta,X_t\rangle^{p};\tau_1\le t<\tau_{m+1}) \nonumber\\
			& =e^{-\beta t}e^{\varphi(p\theta)t}+\bE \int_0^t \beta e^{-\beta s}\E_{\xi_s+\cP} (\langle e_\theta,X_{t-s}\rangle^{p};t-s<\tau_{m})\,ds,\label{formula-V}
		\end{align}
		where in the last equality, we use that fact that $(X_{\tau_1},\P)\overset{d}{=}(\xi_\tau+\mathcal{P},\bP)$.
		For any $\nu=\sum_{i=1}^K\delta_{x_i}\in\mathcal{M}_a$, we have that,  under $\P_{\nu}$,
		$$
		X_t=\sum_{i=1}^K(x_i+X_t^{i}),
		$$
		where $X_t^i$, $i=1,\cdots,K,$ are iid with the same law as $X_t$ under $\P$. Let $\tau_m^i, i=1,\cdots, K$  be the $m$-th branching time of the subtree $X^i$.
			Then, we have that
		\begin{align*}
			&\E_{\nu} (\langle e_\theta,X_{t}\rangle^{p}; t<\tau_{m})=  \E_{\nu}\left[ \left(\sum_{i=1}^Ke^{\theta x_i}\langle e_\theta,X_t^i\rangle \right)^{p};t<\tau_m\right]\\
			\le &\E_{\nu} \left[\left(\sum_{i=1}^Ke^{\theta x_i}\langle e_\theta,X_t^i\rangle\textbf{1}_{t< \tau_{m}^i}\right)^{p}\right]
			\le\langle e_\theta,\nu\rangle^{p-1}\E_{\nu} \left(\sum_{i=1}^Ke^{\theta x_i}\langle e_\theta,X_t^i\rangle^p\textbf{1}_{t< \tau_{m}^i}\right)\\
			=&\langle e_\theta,\nu\rangle^{p}V_m(t),
		\end{align*}
		where in the second inequality, we used H\"older's  inequality.
		Plugging this  into \eqref{formula-V}, we get that
		\begin{align*}
			V_{m+1}(t) &\le e^{-\beta t}e^{\varphi(p\theta)t}+\bE \int_0^t \beta e^{-\beta s}e^{p\theta\xi_s}\langle e_\theta,\mathcal{P}\rangle^p V_m(t-s)\,ds\\
			&=e^{-\beta t}e^{\varphi(p\theta)t}+\bE(\langle e_\theta,\mathcal{P}\rangle^p) \int_0^t \beta e^{-\beta s}e^{\varphi(p\theta)s} V_m(t-s)\,ds.
		\end{align*}
		Using  induction on $m$,  we  can obtain  that for all $t>0$,
		$$
		V_m(t)\le e^{c_{p,\theta} t},
		$$
		where $c_{p,\theta}=\beta(\bE(\langle e_\theta,\mathcal{P}\rangle^p)-1)+\varphi(p\theta).$
		Therefore
		\begin{align}\label{222}
			\E[ \langle e_\theta,X_t\rangle^{p}]=\sup_{m}V_{m}(t)\le e^{c_{p,\theta} t}<\infty.
		\end{align}

(3) We finally prove the sufficiency in the case $\bE(N)=\infty$.
Using \eqref{222} for $X^{(n)}$, we get that
		\begin{align*}
			\E(\langle e_\theta,X_t^{(n)}\rangle^p)&\le \exp\{(\beta(\bE(\langle e_\theta,\mathcal{P}^{(n)}\rangle^p)-1)+\varphi(p\theta)) t\}\\
			&\le \exp\{(\beta(\bE(\langle e_\theta,\mathcal{P}\rangle^p)-1)+\varphi(p\theta)) t\}.
		\end{align*}
		Thus
		$$\E(\langle e_\theta,X_t\rangle^p)=\lim_{n\to\infty}\E(\langle e_\theta,X_t^{(n)}\rangle^p)\le \exp\{(\beta(\bE(\langle e_\theta,\mathcal{P}\rangle^p)-1)+\varphi(p\theta)) t\}.$$
		
		\hfill$\Box$

\subsection{Finiteness of  $\E(W_t L^*(W_t))$ }

	In this subsection, we study the finiteness of
		$\E(W_t(\theta) L^*(W_t(\theta)))$ with $L\in \mathbb{L}$.

\begin{theorem}\label{prop:moment7}
	Let  $\theta\in\Theta_0$ and  $L\in\mathbb{L}$. Then
	$$
	\E(W_t(\theta) L^*(W_t(\theta)))<\infty \quad \mbox{ for some } t>0
	$$
	if and only if
	\begin{align}\label{condition: moment7}
		\bE( \langle e_\theta,\mathcal{P}\rangle L^*(\langle e_\theta,\mathcal{P}\rangle))<\infty.
	\end{align}
	In particular, if  $
	\E(W_t(\theta) L^*(W_t(\theta)))<\infty
	$ for some $t>0$, then $
	\E(W_t(\theta) L^*(W_t(\theta)))<\infty
	$ for all $t>0$.
\end{theorem}
	
To prove Theorem \ref{prop:moment7}, we need some preparations.
	For $\theta\in\Theta$, put
\begin{equation}\label{def-Xi'}
	Y_\theta:=e^{-\kappa(\theta)\tau}e^{\theta \xi_\tau}.
\end{equation}
Note that
\begin{align}\label{W-t-tau}
	W_{t\wedge \tau_1}(\theta)\overset{d}{=}{\bf 1}_{\{\tau< t\}}Y_\theta\langle e_\theta,\mathcal{P}\rangle+{\bf 1}_{\{\tau\ge t\}}e^{-\kappa(\theta)t}e^{\theta\xi_t}.
\end{align}
If  $\theta\in\Theta_0$, then there exists  $r>1$ such that $\kappa(r\theta)<\infty$, which implies that $\bE(e^{r\theta\xi_t})<\infty$ and
\begin{align}\label{exp-Y-t}
\bE(Y_\theta^r,\tau<t)=\int_0^t \beta e^{-\beta u}e^{-\kappa(\theta)ru}e^{\varphi(r\theta)u}\,du<\infty.
\end{align}
Thus by \eqref{W-t-tau}, the independence of $\xi, \tau$ and $\cP$, we get that  for $\theta\in\Theta_0$,
\begin{align}\label{W-t-tau-2}
	\bE(\langle e_\theta,\mathcal{P}\rangle L^*(\langle e_\theta,\mathcal{P}\rangle))<\infty\Leftrightarrow \bE(W_{t\wedge \tau_1}(\theta)L^*(W_{t\wedge \tau_1}(\theta)))<\infty.
\end{align}

	Let $\phi(t,\lambda):=\E(e^{-\lambda W_t(\theta)})$ and
		\begin{align}\label{def-psi-t}
			\phi_1(t,\lambda):=\phi(t,\lambda)-1+\lambda.
		\end{align}
		By Lemma \ref{lemma:tb}, $\E(W_t L^*(W_t))<\infty$ if and only if
		\begin{align}
			\int_0^1 \phi_1(t,\lambda)L(1/\lambda)\lambda^{-2}\,d\lambda<\infty.
		\end{align}
		
		We first give some properties of the Laplace transform of $W_t(\theta)$.
Recall that $\tau_1$ is the first branching time.  Let $\varnothing$ be
the initial particle at time 0.
By the strong Markov property and the branching property of $X$, we have
		\begin{align}\label{decomp:Wt}
			W_t(\theta) & ={\bf 1}_{\tau_1<t}e^{-\kappa(\theta)\tau_1}e^{\theta z_\varnothing(\tau_1) }\sum_{j=1}^{N^\varnothing}e^{\theta S_i^\varnothing}W_{t-\tau_1}^j(\theta)+{\bf 1}_{\tau_1 \ge t}W_t(\theta),
		\end{align}
		where $W_t^j(\theta)$, $j=1,\cdots,N^\varnothing$, are iid with the same law as $W_t(\theta)$ under $\P$, and are independent of $X_{\tau_1}$.
	Then we have that
		\begin{align}\label{phi-t}
			\phi(t,\lambda)&=e^{-\beta t}\bE\Big(e^{-\lambda e^{-\kappa(\theta)t}e^{\theta\xi_t}}\Big)+\bE\left(\prod_{i=1}^N\phi\left(t-\tau,\lambda Y_\theta e^{\theta S_i}\right);\tau< t\right)\nonumber\\
			&=e^{-\beta t}\bE\Big(e^{-\lambda e^{-\kappa(\theta)t}e^{\theta\xi_t}}\Big)+\bE\left(e^{-I(t,\lambda)};\tau<  t\right),
		\end{align}
		where
		$$I(t,\lambda):=-\sum_{i=1}^N\log \phi\left(t-\tau,\lambda Y_\theta e^{\theta S_i}\right){\bf 1}_{\tau< t}.$$
Let
\begin{equation}\label{def-Xi'-2}
\Xi_\theta=Y_\theta \langle e_\theta,\mathcal{P}\rangle.
\end{equation}
Note hat $\Xi_\theta\overset{d}{=}W_{\tau_1}(\theta)$,  where $\tau_1$ is the first splitting  time.

		\begin{lemma}\label{lemma5}
If $\theta\in \Theta$, then for any $t>0$ and $\lambda\ge 0$,
			\begin{align}\label{est:phi1-t}
				0&\le \phi_1(t,\lambda)-\bE\left(\sum_{i=1}^N\phi_1(t-\tau,\lambda Y_\theta e^{\theta S_i};\tau<t)\right)
				\le \E \Big(e^{-\lambda W_{t\wedge \tau_1}(\theta)}-1+\lambda W_{t\wedge \tau_1}(\theta)\Big).
			\end{align}
		\end{lemma}
		{\bf Proof:}
		By  \eqref{def-psi-t} and \eqref{phi-t}, we have
		\begin{align}\label{11.0}
			\phi(t,\lambda)=&e^{-\beta t} \bE \Big(e^{-\lambda e^{-\kappa(\theta)t}e^{\theta\xi_t}}-1+\lambda e^{-\kappa(\theta)t}e^{\theta\xi_t}\Big)+e^{-\beta t}(1-\lambda \bE (e^{-\kappa(\theta)t}e^{\theta\xi_t}))\nonumber\\
			&+\bE\Big(e^{-I(t,\lambda)}-1+\sum_{i=1}^N(1-\phi(t-\tau,\lambda Y_\theta e^{\theta S_i}));\tau<t\Big)\nonumber\\
			&+\bP(\tau<t)+\bE\Big(\sum_{i=1}^N(\phi_1(t-\tau,\lambda Y_\theta e^{\theta S_i}));\tau<t\Big)\nonumber\\
&-\lambda \bE\Big(\sum_{i=1}^N Y_\theta e^{\theta S_i};\tau<t\Big).
		\end{align}
				Elementary calculations give that $\bP(\tau<t)=1-e^{-\beta t}$,
		\begin{align*}
			e^{-\beta t}(1-\lambda \bE (e^{-\kappa(\theta)t}e^{\theta\xi_t}))=e^{-\beta t}-\lambda e^{-\beta t-\kappa(\theta)t+\varphi(\theta)t}=e^{-\beta t}-\lambda e^{-\beta\chi(\theta)t}
		\end{align*}
and
		\begin{align}\label{11.20}
			\bE\left(\sum_{i=1}^N Y_\theta e^{\theta S_i};\tau<t\right)=\chi(\theta)\int_0^t \beta e^{-\beta s} e^{-\kappa(\theta)s}e^{\varphi(\theta) s}\,ds=1-e^{-\beta\chi(\theta)t}.
		\end{align}
		Thus by \eqref{11.0} we have
		\begin{align}\label{11.2}
			\phi_1(t,\lambda)=&\phi(t,\lambda)-1+\lambda=\bE\left(\sum_{i=1}^N\phi_1(t-\tau,\lambda Y_\theta e^{\theta S_i};\tau<t\right)\nonumber\\
			&+e^{-\beta t} \bE \Big(e^{-\lambda e^{-\kappa(\theta)t}e^{\theta\xi_t}}-1+\lambda e^{-\kappa(\theta)t}e^{\theta\xi_t}\Big)\nonumber\\
			&+\bE\Big(e^{-I(t,\lambda)}-1+\sum_{i=1}^N(1-\phi(t-\tau,\lambda Y_\theta e^{\theta S_i}));\tau<t\Big).
		\end{align}
		For any $x_i\in (0,  1], i\ge1$, one can easily check that
		$$1-\prod_{i=1}^n x_i\le \sum_{i=1}^n(1-x_i)\le \sum_{i=1}^n (-\log x_i).$$
		Thus on the event $\{\tau<t\}$, we have
		\begin{align}\label{11.5}
			0\le e^{-I(t,\lambda)}-1+\sum_{i=1}^N(1-\phi(t-\tau,\lambda Y_\theta e^{\theta S_i}))\le e^{-I(t,\lambda)}-1+I(t,\lambda).
		\end{align}
		By Jensen's inequality, $-\log\phi(t,\lambda)\le \lambda.$ Thus on the event $\{\tau<t\}$, we have
		$$
		I(t,\lambda)\le \sum_{i=1}^N \lambda Y_\theta e^{\theta S_i}=\lambda \Xi_\theta.
		$$
		Since $e^{-x}-1+x$ is increasing on $(0,\infty)$,
we have
		\begin{align}\label{11.6}
			e^{-I(t,\lambda)}-1+I(t,\lambda)\le e^{-\lambda \Xi_\theta}-1+\lambda\Xi_\theta.
		\end{align}
		Combining \eqref{11.2}, \eqref{11.5} and \eqref{11.6}, we get that
		\begin{align*}
			0&\le \phi_1(t,\lambda)-\bE(\sum_{i=1}^N\phi_1(t-\tau,\lambda Y_\theta e^{\theta S_i};\tau<t)\\
			 &\le \bE \Big(e^{-\lambda e^{-\kappa(\theta)t}e^{\theta\xi_t}}-1+\lambda e^{-\kappa(\theta)t}e^{\theta\xi_t};\tau\ge t\Big)+\bE(e^{-\lambda \Xi_\theta}-1+\lambda\Xi_\theta;\tau<t),
		\end{align*}
where in the last inequality we used the fact that $\bP(\tau\ge t)=e^{-\beta t}$.
		Since $W_{t\wedge \tau_1}(\theta)\overset{d}{=}{\bf 1}_{\{\tau<t\}}\Xi_\theta+{\bf 1}_{\{\tau\ge t\}}e^{-\kappa(\theta)t}e^{\theta\xi_t}$,
		the  desired result  follows immediately.
		\hfill$\Box$
		
		\bigskip

		For  fixed $t>0$, we define a  measure $\mu_t$ on $\R^+$ by
		\begin{align}\label{def:mut}
			\mu_t(B):=(1-e^{-\beta\chi(\theta)t})^{-1}\bE \Big(\sum_{i=1}^N{\bf 1}_B(Y_\theta e^{\theta S_i})Y_\theta e^{\theta S_i};\tau<t\Big),
		\end{align}
	for any Borel subset of $\R^+$. 	
		By \eqref{11.20}, we have $\mu_t(\R^+)=1$. Thus $\mu_t$ is a probability measure on $\R^+$.
		Let $Z_t$ be a random variable with distribution $\mu_t$.
		For any $r\in\R$, if $\kappa((r+1)\theta)<\infty$, then
		\begin{align}\label{moment-mu-t}
			 E(Z_t^r)&=(1-e^{-\beta\chi(\theta)t})^{-1}\chi((1+r)\theta)\bE(Y_\theta^{1+r};\tau<t)<\infty.
		\end{align}
		
		For any nonnegative Borel function $F(t,\lambda)$ on $\R^+\times \R^+$, define
		\begin{align}\label{def:U}
			U(F)(t,\lambda):=\bE\left(\sum_{i=1}^N F(t-\tau,\lambda Y_\theta e^{\theta S_i});\tau<t\right).
		\end{align}
		We also define $U^{(1)}(F):=U(F)$ and, for $k\ge 2$,
		$U^{(k)}(F):=U(U^{(k-1)}(F)).$
		Note that if $F$ is increasing in $t$, then
		\begin{align}\label{11.9}
			U(F)(t,\lambda)\le \bE\left(\sum_{i=1}^N F(t,\lambda Y_\theta e^{\theta S_i});\tau<t\right)=\lambda(1-e^{-\beta\chi(\theta)t})E(\tilde{F}(t,\lambda Z_t)),
		\end{align}
		where $\tilde{F}(t,\lambda):=F(t,\lambda)/\lambda$.

		\begin{remark}\label{remark:delta-t}
			If $\theta\in\Theta_0$, then there exists $r\in (0, 1)$ such that $\kappa((1+r)\theta)<\infty$ and $\kappa((1-r)\theta)<\infty$.
			By \eqref{moment-mu-t},  we have $E(Z_t^r)<\infty$ and $E(Z_t^{-r})<\infty$ for all $t>0$. For any $t>0$, by the dominated convergence theorem, we have
			$$\lim_{\delta\to0}(1-e^{-\beta\chi(\theta)t})E (Z_t^{\delta}\vee Z_t^{-\delta})=1-e^{-\beta\chi(\theta)t}<1.$$
			Thus we can choose $\delta_t\in(0,1)$ small enough such that
			$$(1-e^{-\beta\chi(\theta)t})E (Z_t^{\delta_t}\vee Z_t^{-\delta_t})<1.$$

			For $L_1,L_2\in\mathbb{L}$, if $L_1(x)\sim L_2(x)$ as $x\to\infty$, then one can show that $L_1^*(x)\sim L_2^*(x)$, as $x\to\infty.$  Thus without loss of generality,
in the proofs of Lemma \ref{lem:psito0-6} and Theorem \ref{prop:moment7},
we can assume that, for any fixed $t>0$, $L\in\mathbb{L}_{\delta_t}$ and $L$ is differentiable.
			Then
			$xL(x)$ is increasing. Put $h(x)=xL^*(x)$. Then for $x>1$, $h''(x)=L(x)/x+L'(x)=x^{-1}(xL(x))'>0.$ Thus $h$ is convex.
			\end{remark}

		\begin{lemma}\label{lem:psito0-6}
			Let  $\theta\in\Theta_0$ and  $L\in\mathbb{L}$.  If  \eqref{condition: moment7} holds,
			then for any $t>0$,
			$$\lim_{\lambda\to 0}\lambda^{-1}L^*(1/\lambda)\phi_1(t,\lambda)=0.$$
		\end{lemma}
		{\bf Proof:}
For $s>0$ and $\lambda>0$, define \begin{equation}\label{def-Delta1}
\Delta_1(s,\lambda):=\phi_1(s,\lambda)-U(\phi_1)(s,\lambda),
\end{equation}
and for $k\ge 2$,
		\begin{align*}
			\Delta_k(s,\lambda)&:=U(\Delta_{k-1})(s,\lambda)=U^{(k-1)}(\phi_1)(s,\lambda)-U^{(k)}(\phi_1)(s,\lambda).
		\end{align*}
		We claim that $\lim_{k\to\infty}U^{(k)}(\phi_1)(s,\lambda)=0$ for any $s>0$ and $\lambda\ge0$.
		
		In fact, since $e^{-\lambda x} $ is convex,  $\{e^{-\lambda W_s(\theta)}, s\ge0\}$ is a submartingale. This implies that $\phi_1(s,\lambda)$ is increasing in $s$. Thus  $U(\phi_1)(s,\lambda)$ is increasing in $s$, which implies that $U^{(k)}(\phi_1)(s,\lambda)$ is increasing in $s$ for any $k\ge 1$.
		For any $s>0$, let  $Z_{s,k}, k\ge 1,$ be iid random variables  with the law $\mu_s$.
		Using \eqref{11.9} and induction on $k$, we can get that
		$$
		U^{(k)}(\phi_1)(s,\lambda)=U(U^{k-1}(\phi_1))(s,\lambda)\le \lambda (1-e^{-\beta\chi(\theta)s})^kE(\psi_1(s,\lambda Z_{s,1}\cdots Z_{s,k})),
		$$
where
$\psi_1(s,\lambda)
:=\phi_1(s,\lambda)/\lambda.$
Since $0\le \psi_1(s,\lambda)\le1$, we have
$$
(1-e^{-\beta\chi(\theta)s})^{k}E(\psi_1(s,\lambda Z_{s,1}Z_{s,2}\cdots Z_{s,k}))\to 0,\quad k\to\infty.
$$
Thus the claim is valid and consequently
		\begin{align}\label{phi1-s-lambda}
			\phi_1(s,\lambda)=\sum_{k=1}^\infty \Delta_k(s,\lambda).
		\end{align}

		Fix an arbitrary $t>0$. Choose $\delta_t$ as in Remark \ref{remark:delta-t}.
		By Remark \ref{remark:delta-t}, we can assume $L\in\mathbb{L}_{\delta_t}$ and $L$ is differentiable.
		Since $L^*$ is slowly varying, there exists $x_0>1$ such that $x^{-\delta_t}L^*(x)$ is decreasing on $[x_0,\infty)$.
		Let $\bar{J}(x)=L^*(x\vee x_0)$. Thus $x^{-\delta_t}\bar{J}(x)$ is decreasing on $(0,\infty)$ and $\bar{J}$ is increasing on $(0,\infty)$.
		Thus for any $\lambda>0$ and $y>0$, $$\bar{J}(1/\lambda)(1\wedge (\lambda y))\le \bar{J}(y),$$
		which is easily seen to be true by dealing with the cases $\lambda y>1$ or $\lambda y\le 1$ separately.
Thus by Lemma \ref{lemma5}, \eqref{W-t-tau-2} and the fact that $e^{-x}-1+x\le x\wedge x^2$ for all $x\ge 0$, we get
		\begin{align}
			\lambda^{-1}\bar{J}(1/\lambda)|\Delta_1(s,\lambda)|\le &\bar{J}(1/\lambda)\E( W_{s\wedge \tau_1}(\theta)(1\wedge\lambda W_{s\wedge \tau_1}(\theta))
			\label{7.44}\\
			\le &\bE(W_{s\wedge \tau_1}(\theta)\bar{J}(W_{s\wedge \tau_1}(\theta))):=c(s)<\infty.\label{7.32}
		\end{align}
				 Combining this with the dominated convergence theorem, we get that
		for $s>0$,
		\begin{equation}\label{7.31}
			\lim_{\lambda\to 0}\lambda^{-1}\bar{J}(1/\lambda)|\Delta_1(s,\lambda)|=0.
		\end{equation}
		
		We claim that for any $k\ge 1$, and $s\in[0,t]$,
		\begin{align}\label{7.34}
			\lim_{\lambda\to 0}\lambda^{-1}\bar{J}(1/\lambda)|\Delta_k(s,\lambda)|=0
		\end{align}
		and
		\begin{equation}\label{7.33}
			\lambda^{-1}\bar{J}(1/\lambda)|\Delta_k(s,\lambda)|\le c(s)A_t(s)^{k-1},
		\end{equation}
		where $$A_t(s):=(1-e^{-\beta\chi(\theta)s})E(1\vee Z_s^{\delta_t})=\bE \left(\sum_{i=1}^N (1\vee Y_\theta e^{\theta S_i})^{\delta_t} Y_\theta e^{\theta S_i};\tau<s\right)\le A_t(t)<1.$$
		
		We will prove \eqref{7.34} and \eqref{7.33} by induction.
		It is clear  \eqref{7.34} and \eqref{7.33} hold for $k=1$. Assume that  \eqref{7.34}  and \eqref{7.33}  hold for $k$. Recall that
		\begin{align}\label{recall-delta}
			\Delta_{k+1}(s,\lambda)=U(\Delta_{k})(s,\lambda)=E\left(\sum_{i=1}^N\Delta_k(s-\tau,\lambda Y_\theta e^{\theta S_i}){\bf 1}_{\{\tau<s\}}\right).
		\end{align}
		Note that  $x\bar{J}(x)$ is still convex.
		Hence $\{W_{s\wedge\tau_1}(\theta)\bar{J}(W_{s\wedge\tau_1}(\theta)),0\le s\}$ is a submartingale, which implies that $c(s)$ is increasing on $[0,\infty)$.
		By \eqref{7.33} for $k$
		and the fact that $c(s)$ and $ A_t(s)$ are  increasing,  we have for $t\in(0,t]$,
		\begin{align}\label{7.38}
			&\sum_{i=1}^N\Delta_k(s-\tau,\lambda Y_\theta e^{\theta S_i}){\bf 1}_{\{\tau<s\}}\le c(s)A_t(s)^{k-1}\sum_{i=1}^N ( \lambda Y_\theta e^{\theta S_i})\bar{J}^{-1}(\lambda^{-1} Y_\theta^{-1} e^{-\theta S_i}){\bf 1}_{\{\tau<s\}}\nonumber\\
			&\le \lambda \bar{J}^{-1}(\lambda^{-1} )c(s)A_t(s)^{k-1}\sum_{i=1}^N (1\vee Y_\theta e^{\theta S_i})^{\delta_t} Y_\theta e^{\theta S_i}{\bf 1}_{\{\tau<s\}},
		\end{align}
		where we used the fact that $\frac{\bar{J}(ux)}{\bar{J}(x)}\le 1\vee u^{\delta_t}$. Thus by \eqref{recall-delta} we have
		\begin{align*}
			\Delta_{k+1}(s,\lambda)
			\le &\lambda \bar{J}^{-1}(\lambda^{-1} )c(s)A_t(s)^{k-1}\bE \left(\sum_{i=1}^N (1\vee Y_\theta e^{\theta S_i})^{\delta_t} Y_\theta e^{\theta S_i};\tau<s\right)\\
			=&\lambda \bar{J}^{-1}(\lambda^{-1} )c(s)A_t(s)^{k}.
		\end{align*}
		That is,  \eqref{7.33} holds for $k+1$.
		Now combining \eqref{recall-delta},  \eqref{7.34}, \eqref{7.38} and the dominated convergence theorem, we get   that \eqref{7.34} is valid for $k+1$.
		
		Recall the decomposition of $\phi_1(s,\lambda)$ in \eqref{phi1-s-lambda}. By \eqref{7.34}, \eqref{7.33}  and  the dominated convergence theorem, we get
		\begin{align*}
			\lim_{\lambda\to 0}\lambda^{-1}\bar{J}(1/\lambda)\phi_1(t,\lambda) =  \lim_{\lambda\to 0}\lambda^{-1}\bar{J}(1/\lambda)\sum_{k=1}^\infty\Delta_k(t,\lambda) =0.
		\end{align*}
		The desired result follows immediately from the fact $L^*(x)\sim \bar{J}(x)$ as $x\to\infty.$
		\hfill$\Box$
		
		\bigskip

\noindent{\bf Proof of Theorem \ref{prop:moment7}}\quad
Fix an arbitrary $t>0$.
		If $L^*(\infty)<\infty,$ then Theorem \ref{prop:moment7} holds trivially. In the remainder of this proof, we assume that $L^*(\infty)=\infty$. Without loss of generality, we assume that $L\in\mathbb{L}_{\delta_t}$
		
		We first prove the sufficiency.
		Assume that
		$\E(\langle e_\theta,\mathcal{P}\rangle L^*(\langle e_\theta,\mathcal{P}\rangle))<\infty,$
which, by \eqref{W-t-tau-2},
is equivalent to
\begin{equation}\label{finite-moment}
\bE(W_{t\wedge \tau_1}(\theta)L^*(W_{t\wedge \tau_1}(\theta)))<\infty.\end{equation}
		To prove $\E(W_t(\theta)L^*(W_t(\theta)))<\infty,$
by Lemma \ref{lemma:tb}, it suffices to show that
		$$\int_{0}^1 \phi_1(t,\lambda)L(1/\lambda)\lambda^{-2}\,d\lambda<\infty.$$
Recall the definition of $\Delta_1(t, x)$  in \eqref{def-Delta1}.
By Lemma \ref{lemma5},
$$0\le \Delta_1(t,\lambda)\le \E \Big(e^{-\lambda W_{t\wedge \tau_1}(\theta)}-1+\lambda W_{t\wedge \tau_1}(\theta)\Big)=\E (\mathcal{T}_1(\lambda W_{t\wedge\tau_1})).$$
Hence by Lemma \ref{lemma:tb} and \eqref{finite-moment}, we have
		\begin{align}\label{4.1.1}
			\int_{0}^1 \Delta_1(t,\lambda)L(1/\lambda)\lambda^{-2}\,d\lambda<\infty.
		\end{align}	
Since $\phi_1(t,\lambda)\le \lambda$, we have  for any $a\in(0,1)$,
$$\int_a^1 \phi_1(t,\lambda)L(1/\lambda)\lambda^{-2}\,d\lambda <\infty,\quad \mbox{and }\int_a^1 U( \phi_1)(t,\lambda)L(1/\lambda)\lambda^{-2}\,d\lambda <\infty.$$
				Note that for any $a\in(0,1)$,
		\begin{align}\label{11.11}
			&\int_{a}^1 \phi_1(t,\lambda)L(1/\lambda)\lambda^{-2}\,d\lambda\nonumber\\
\le &\int_{0}^1 \Delta_1(t,\lambda)L(1/\lambda)\lambda^{-2}\,d\lambda+\int_{a}^1 U(\phi_1)(t,\lambda)L(1/\lambda)\lambda^{-2}\,d\lambda\nonumber\\
			\le &\int_{0}^1\!\! \Delta_1(t,\lambda)L(1/\lambda)\lambda^{-2}\,d\lambda+(1-e^{-\beta\chi(\theta)t})\!\!\int_a^1\!\! E(\psi_1(t,\lambda Z_t))L(1/\lambda)\lambda^{-1}d\lambda,
		\end{align}
where in the last inequality we used
		$$U( \phi_1)(t,\lambda)\le \lambda(1-e^{-\beta\chi(\theta)t})E(\psi_1(t,\lambda Z_t)),$$
		which follows from \eqref{11.9}.
By a change of variables, we have
		\begin{align}\label{7.7}
&(1-e^{-\beta\chi(\theta)t})\int_a^1 E(\psi_1(t,\lambda Z_t))L(1/\lambda)\lambda^{-1}d\lambda\nonumber\\
=&(1-e^{-\beta\chi(\theta)t})E \int_{aZ_t}^{Z_t} \psi_1(t,\lambda )L(Z_t/\lambda)\lambda^{-1}d\lambda\nonumber\\
\le &(1-e^{-\beta\chi(\theta)t})E \int_{a}^{1} \psi_1(t,\lambda )L(Z_t/\lambda)\lambda^{-1}d\lambda\nonumber\\
&+(1-e^{-\beta\chi(\theta)t})E\Big[\int_{aZ_t}^{a} \psi_1(t,\lambda )L(Z_t/\lambda)\lambda^{-1}d\lambda; Z_t\le 1\Big]\nonumber\\
&+(1-e^{-\beta\chi(\theta)t})E\Big[\int_{1}^{Z_t} \psi_1(t,\lambda )L(Z_t/\lambda)\lambda^{-1}d\lambda;Z_t>1\Big].
		\end{align}
Since $L\in\mathbb{L}_{\delta_t}$,  $\psi_1(t,\lambda)\le 1$ and $\psi_1(t,\cdot)$ is increasing, we have
\begin{align*}
&(1-e^{-\beta\chi(\theta)t})E \int_{a}^{1} \psi_1(t,\lambda )L(Z_t/\lambda)\lambda^{-1}d\lambda\\
\le &(1-e^{-\beta\chi(\theta)t})E \Big(Z_t^{\delta_t}\vee Z_t^{-\delta_t}\Big) \int_{a}^1 \phi_1(t,\lambda)L(1/\lambda)\lambda^{-2}\,d\lambda,\end{align*}
\begin{align*}
&(1-e^{-\beta\chi(\theta)t})E\Big[\int_{aZ_t}^{a} \psi_1(t,\lambda )L(Z_t/\lambda)\lambda^{-1}d\lambda; Z_t\le 1\Big]\\
\le&(1-e^{-\beta\chi(\theta)t}) \psi_1(t,a)L(a^{-1})E(aZ_t)^{-\delta_t}\int_{0}^a \lambda^{\delta_t-1}d\lambda\\=&\psi_1(t,a)L(a^{-1})\delta_t^{-1}(1-e^{-\beta\chi(\theta)t})E \Big(Z_t^{-\delta_t}\Big)\le a^{-1}\phi_1(t,a)L(a^{-1})\delta_t^{-1},
\end{align*}
and
\begin{align*}
&(1-e^{-\beta\chi(\theta)t})E\Big[\int_{1}^{Z_t} \psi_1(t,\lambda )L(Z_t/\lambda)\lambda^{-1}d\lambda;Z_t>1\Big]\\
\le &(1-e^{-\beta\chi(\theta)t})L(1)E\int_1^\infty (Z_t/\lambda)^{\delta_t}\lambda^{-1}d\lambda\\
=&L(1)\delta_t^{-1}(1-e^{-\beta\chi(\theta)t})E \Big(Z_t^{\delta_t}\Big)\le L(1)\delta_t^{-1}.
\end{align*}
Combining these three displays with  \eqref{11.11} and \eqref{7.7}, we have that
		\begin{align*}
			0<& \Big(1-(1-e^{-\beta\chi(\theta)t})E \Big(Z_t^{\delta_t}\vee Z_t^{-\delta_t}\Big)\Big)
			\int_{a}^1 \phi_1(t,\lambda)L(1/\lambda)\lambda^{-2}\,d\lambda\\
			\le &\int_{0}^1 \Delta_1(t,\lambda)L(1/\lambda)\lambda^{-2}\,d\lambda+a^{-1}\phi_1(t,a)L(a^{-1})\delta_t^{-1}+\delta_t^{-1}L(1).
		\end{align*}
		Since $L(a^{-1})/L^*(a^{-1})\to0$,
		we have by Lemma \ref{lem:psito0-6} that  $a^{-1}L(a^{-1})\phi_1(t,a)\to 0$ as $a\to0$.
		Thus
		$$\int_{0}^1 \phi_1(t,\lambda)L(1/\lambda)\lambda^{-2}\,d\lambda<\infty,$$
		which implies that $\E(W_t(\theta)L^*(W_t(\theta)))<\infty.$
		
		Now we prove the necessity.
		Assume that $\E(W_t(\theta) L^*(W_t(\theta)))<\infty$.
		By Remark \ref{remark:delta-t}, we can assume, without loss of generality, that $xL(x)$ is increasing and $L$ is differentiable..
		It follows that  $\{W_s(\theta)L^*(W_s(\theta)),0\le s\le t\}$ is a submartingale.
		Using the bounded optional stopping theorem, we have that
		$$\E( W_{t\wedge\tau_1}(\theta)L^*(W_{t\wedge\tau_1}(\theta)))\le \E(W_t(\theta) L^*(W_t(\theta)))<\infty.$$
Now we can apply \eqref{W-t-tau-2} and complete the proof of the necessity.
		\hfill$\Box$

		\bigskip

In Theorem \ref{proposition-log}
below, we will show that  when $L^*(x)=(\log_+x)^r$ for some $r>0$,
the condition $\theta\in\Theta_0$ in Theorem \ref{prop:moment7}
can be replaced by some weaker conditions.

\begin{theorem}\label{proposition-log}
			Let $\theta\in \Theta $
			and $r>0$. Assume that \eqref{kappa'} holds
			and
			\begin{align}\label{condition 9}
				\bE\left(\sum_{i=1}^N {\bf 1}_{S_i>0}S_i^r e^{\theta S_i}\right)<\infty, \mbox{ and } \int_1^\infty y^r e^{\theta y}n(dy)<\infty.
			\end{align}
			Then $\E(W_t(\theta)(\log_+W_t(\theta))^r)<\infty$ for some $t>0$ if and only if
			\begin{align}\label{conditon10}
				\bE(\langle e_\theta,\cP\rangle (\log_+\langle e_\theta,\cP\rangle)^r)<\infty.
			\end{align}
			In particular, if $\E(W_t(\theta)(\log_+W_t(\theta))^r)<\infty$ for some $t>0$, then $\E(W_t(\theta)(\log_+W_t(\theta))^r)<\infty$ for all $t>0$.
\end{theorem}

		\noindent{\bf Proof:}
		Fix an arbitrary $t>0$. Recall that $Z_t$ is a random variable with distribution $\mu_t$.
		Note that \eqref{kappa'} and \eqref{condition 9} imply that $E[|\log Z_t|]<\infty$ and  $E[(\log_+Z_t)^r]<\infty$.
		In this proof, $t>0$ is a fixed number.	
		By \eqref{W-t-tau} we have
		\begin{align}\label{e:rsnew}
		\E\big(W_{t\wedge \tau_1}(\theta)(\log_+(W_{t\wedge \tau_1}(\theta)))^r\big)=&e^{-\beta t}\bE(e^{-\kappa(\theta)t}e^{\theta\xi_t}(\theta\xi_t-\kappa(\theta)t)_+^r)\nonumber\\
		&+\!\bE(Y_\theta\langle e_\theta,\mathcal{P}\rangle(\log Y_\theta\!+\!\log \langle e_\theta,\mathcal{P}\rangle)_+^r;\tau<t).
		\end{align}
	Since $\int_1^\infty y^r e^{\theta y}n(dy)<\infty$, we have $\bE((\xi_s\vee 0)^re^{\theta\xi_s})<\infty$  for all $s\ge0$. It follows that
	$$\bE(Y_\theta(\log_+ Y_\theta)^r;\tau<t)=\int_0^t \beta e^{-\beta s}e^{-\kappa(\theta)s}\bE(e^{\theta\xi_s}(\theta\xi_s-\kappa(\theta)s)_+^r)\,ds<\infty.$$
Hence by \eqref{e:rsnew}, and the independence of $\xi, \tau$ and $\cP$,
\begin{align}\label{W-t-tau-3}
	\bE(\langle e_\theta,\cP\rangle (\log_+\langle e_\theta,\cP\rangle)^r)<\infty\Leftrightarrow \E\big(W_{t\wedge \tau_1}(\theta)[\log_+(W_{t\wedge \tau_1}(\theta))]^r\big)<\infty.
	\end{align}

 (1) We first consider the case $r\ge1$.
		We choose $C_t>1$  large enough so that $x^{-1}[\log(x\vee C_t)]^{r}$ is decreasing and  $$
		(1-e^{-\beta\chi(\theta)t})E\left(\left(1+\frac{\log_+(Z_t)}{\log C_t}\right)^r\right)<1.
		$$
		Let $L(x)=[\log(x\vee C_t)]^{r-1}$.
		Then
\begin{equation}\label{L*(x)}
L^*(x)=\left(1-\frac{1}{r}\right)(\log C_t)^r+\frac{1}{r}(\log x)^r,\quad x>C_t.
\end{equation}	
 We claim that
		 \begin{align}\label{inequality2}
			\left(1+\frac{\log_+(u^{-1})}{\log C_t}\right)^{-1}\le \frac{\log((ux)\vee C_t)}{\log(x\vee C_t)}\le \left(1+\frac{\log_+u}{\log C_t}\right).
		\end{align}
		In fact, if $0\le u\le1$, then $\frac{\log((ux)\vee C_t)}{\log(x\vee C_t)}\le 1.$ If $u>1$, then
		$$\frac{\log((xu)\vee C_t)}{\log(x\vee C_t)}\le \frac{\log u+\log(x\vee C_t)}{\log(x\vee C_t)}=1+\frac{\log u}{\log(x\vee C_t)}\le 1+\frac{\log u}{\log(C_t)}.
		$$
	
		The proof of necessity is similar as that in the proof of
Theorem \ref{prop:moment7}.
We omit the details. We only need to prove sufficiency. By \eqref{W-t-tau-3} and \eqref{L*(x)},
it suffices to show that  if
 $\E\big(W_{t\wedge \tau_1}(\theta)L^*(W_{t\wedge \tau_1}(\theta))\big)<\infty$,
then
		$\E(W_t(\theta)(L^*(W_t(\theta))<\infty .$

		Recall that $\psi_1(t,\lambda)=\phi_1(t,\lambda)/\lambda$.
Using \eqref{7.7} first and then  \eqref{inequality2} with $u=Z_t$,
we have that  for any $a\in(0,1)$,
		\begin{align}\label{7.7'}
			&\int_a^1 U( \phi_1)(t,\lambda)L(1/\lambda)\lambda^{-2}\,d\lambda
			\le (1-e^{-\beta\chi(\theta)t})\int_a^1 E(\psi_1(t,\lambda Z_t))L(1/\lambda)\lambda^{-1}d\lambda\nonumber\\
		   \le& (1-e^{-\beta\chi(\theta)t}) E((1+\log_+(Z_t)/\log C_t)^{r-1}) \int_{a}^1 \phi_1(t,\lambda)L(1/\lambda)\lambda^{-2}\,d\lambda\nonumber\\
			&+(1-e^{-\beta\chi(\theta)t})E\Big({\bf1}_{Z_t>1}\int_1^{Z_t}
{\psi_1}(t,\lambda)
L(Z_t/\lambda)\lambda^{-1}\,d\lambda\Big)\nonumber\\
			&+(1-e^{-\beta\chi(\theta)t})E\Big({\bf1}_{Z_t<1}\int_{aZ_t}^a
{\psi_1}(t,\lambda)
L(Z_t/\lambda)\lambda^{-1}\,d\lambda\Big).
		\end{align}
		For $1<\lambda<Z_t$,  we have $L(Z_t/\lambda)\le L(Z_t)\le ((\log Z_t)^{r-1}+(\log C_t)^{r-1})$. Thus
		\begin{align*}
			&E\Big({\bf1}_{Z_t>1}\int_1^{Z_t}{\psi_1}(t,\lambda)L(Z_t/\lambda)\lambda^{-1}\,d\lambda\Big)\\
			\le& E\Big({\bf1}_{Z_t>1} ((\log Z_t)^{r-1}+(\log C_t)^{r-1})\int_1^{Z_t}\lambda^{-1}\,d\lambda\Big)\\
			=&E\Big({\bf1}_{Z_t>1} ((\log Z_t)^r+(\log C_t)^{r-1}\log Z_t)\Big)<\infty.
		\end{align*}
		If $\lambda>aZ_t$, then $L(Z_t/\lambda)\le L(1/a)$. Note that ${\psi_1}(t,\lambda)$ is increasing in $\lambda$. Thus we have
		$$E\Big({\bf1}_{Z_t<1}\int_{aZ_t}^a\psi_1(t,\lambda)L(Z_t/\lambda)\lambda^{-1}\,d\lambda\Big)\le \psi_1(t,a)L(1/a)E(|\log Z_t|).$$
		We claim that
\begin{equation}\label{limit-psi_1}
\lim_{a\to 0}{\psi_1}(t,a)L(1/a)=0.
\end{equation}
We will prove the claim later.
		Since
		$$(1-e^{-\beta\chi(\theta)t}) E((1+\log_+(Z_t)/\log C_t)^{r-1})<1,
		$$
		using an argument similar to that  in the proof of sufficiency of
Theorem \ref{prop:moment7}, we can get
		$$ \int_{0}^1 \phi_1(t,\lambda)L(1/\lambda)\lambda^{-2}\,d\lambda<\infty,$$
which, by Lemma \ref{lemma:tb}, implies that
$\E(W_t(\theta)L^*(W_t(\theta)))<\infty.$
		
Now we prove \eqref{limit-psi_1}.
		Put $\bar{J}(x)=[\log (x\vee C_t)]^r.$  Note that $\bar{J}(x)\sim L^*(x)$ as $x\to\infty$.
		Since $x^{-1}(\log(x\vee C_t))^r$ is decreasing, we can verify that for any $u>0$ and $x>0$,
		$$\bar{J}(x)(1\wedge u)\le \bar{J}(ux),$$
		and by \eqref{inequality2},
		$$\bar{J}(ux)/\bar{J}(x)\le \left(1+\frac{\log^+u}{\log C_t}\right)^r.$$
		Using arguments similar to those
		leading to \eqref{7.34} and \eqref{7.33}, we have that  for any $k\ge 1$, and $s\in[0,t]$,
		\begin{align}\label{7.34'}
			\lim_{\lambda\to 0}\lambda^{-1}\bar{J}(1/\lambda)|\Delta_k(s,\lambda)|=0,
		\end{align}
		and
		\begin{equation}\label{7.33'}
			\lambda^{-1}\bar{J}(1/\lambda)|\Delta_k(s,\lambda)|\le c(s)A_t(s)^{k-1},
		\end{equation}
		where $$A_t(s)=(1-e^{-\beta\chi(\theta)s})E((1+\log_+(Z_s)/\log C_t)^r).$$
		Note that $A_t(s)$ is increasing on $s$ and $A_t(s)\le A_t(t)<1$.
		Then using an argument similar to that the proof of Lemma \ref{lem:psito0-6}, we can get that
		\begin{align*}
			\lim_{\lambda\to 0}\lambda^{-1}\bar{J}(1/\lambda)\phi_1(t,\lambda)=0,
		\end{align*}
		which implies that
		${\psi_1}(t,\lambda)L(1/\lambda)=\lambda^{-1}\phi_1(t,\lambda)L(1/\lambda)\to 0,$ as $\lambda\to0$.
		
    (2) Now we
 consider the case when $r\in(0,1)$. We choose $C_t>1$  large enough such that
 $$
 \tilde{C}_t:=(1-e^{-\beta\chi(\theta)t})E\left(1+\frac{\log_+(Z_t^{-1})}{\log C_t}\right)<1.
 $$
 By \eqref{inequality2}, we have that for any $u\ge 0,x\ge 0$,
		\begin{align}\label{inequality3}
			\frac{L(ux)}{L(x)}=\left(  \frac{\log_+(x\vee C_t)}{\log_+((ux)\vee C_t)}\right)^{1-r}\le 1+\frac{\log_+(u^{-1})}{\log C_t}.
		\end{align}
	Using an argument similar to that  leading to \eqref{7.7'}, we can get
	that  for any $a\in(0,1)$,
	\begin{align*}
		&\int_a^1 U( \phi_1)(t,\lambda)L(1/\lambda)\lambda^{-2}\,d\lambda\\
		\le& (1-e^{-\beta\chi(\theta)t}) E(1+\log_+(Z_t^{-1})/\log C_t) \int_{a}^1 \phi_1(t,\lambda)L(1/\lambda)\lambda^{-2}\,d\lambda\nonumber\\
		&+(1-e^{-\beta\chi(\theta)t})E\Big({\bf1}_{Z_t>1}\int_1^{Z_t}{\psi_1}(t,\lambda)L(Z_t/\lambda)\lambda^{-1}\,d\lambda\Big)\nonumber\\
		&+(1-e^{-\beta\chi(\theta)t})E\Big({\bf1}_{Z_t<1}\int_{aZ_t}^a{\psi_1}(t,\lambda)L(Z_t/\lambda)\lambda^{-1}\,d\lambda\Big)\\
		\le& \tilde{C}_t\int_{a}^1 \phi_1(t,\lambda)L(1/\lambda)\lambda^{-2}\,d\lambda
		+E\Big({\bf1}_{Z_t>1}\int_1^{Z_t}\lambda^{-1}\,d\lambda\Big)+E\Big({\bf1}_{Z_t<1}\int_{aZ_t}^a\lambda^{-1}\,d\lambda\Big)\\
		= & \tilde{C}_t\int_{a}^1 \phi_1(t,\lambda)L(1/\lambda)\lambda^{-2}\,d\lambda+E|\log Z_t|,
	\end{align*}
where in the second inequality, we used the fact ${\psi_1}(t,\lambda)\le 1$ and $L(\lambda)\le 1$.
Then using an argument similar to that  in the case of $r\ge1$,
we can show that the desired result holds for $r\in(0,1)$.
		\hfill$\Box$

			\begin{remark} \label{rek:addit-mar}
			Some necessary and sufficient conditions for the non-degeneracy of the limits of the  additive martingales of branching random walks were given in \cite{BK97}.
			For branching L\'evy process, assuming that \eqref{kappa'} holds,
			we can apply the results of \cite{BK97} to $\{W_{n}(\theta),n\ge0\}$ to get that
			$\{W_{n}(\theta),n\ge0\}$ is uniformly integrable if and only if
			$$\theta\kappa'(\theta)<\kappa(\theta), \mbox{ and } \E(W_1(\theta)\log_+(W_1(\theta)))<\infty.$$
By Theorem  \ref{proposition-log},
the condition above is equivalent to
			$$\theta\kappa'(\theta)<\kappa(\theta), \mbox{ and } \bE(\langle e_\theta,\cP\rangle \log_+\langle e_\theta,\cP\rangle)<\infty.$$
			This gives an alternative proof of the necessary and sufficient condition for the non-degeneracy  of $W_\infty(\theta)$ given in \cite{BM18}.
		\end{remark}

		\section{Properties of $W_\infty$}\label{moment-W}
		\subsection{Moments of $W_\infty$}
	Recall that for our branching process $(X_t)_{t\ge 0}$, the skeleton $\{X_{n},n=0,1,2,\cdots\}$ is a branching  random walk.
	Applying  \cite[Theorem 2.1]{Liu2000} with $Y_1=W_1(\theta)$ and $\{A_j,j=1,\cdots,N\}=\{e^{-\kappa(\theta)}e^{\theta z_u(1)}, u\in\mathcal{L}_1\}$, we get that, for any  $p>1$,  the martingale limit $W_\infty(\theta)=\lim_{n\to\infty}W_n(\theta)$ has finite $p$-th moment if and only if
	$$\kappa(p\theta)<p\kappa(\theta), \mbox{ and }  \E(W_1(\theta)^p)<\infty.$$
	Now applying Theorem \ref{wtp}, we immediately get the following result.
	\begin{theorem}\label{main-Lp}
		Let $\theta\in\Theta$ and $p>1.$ Assume that \eqref{kappa'} and \eqref{LLogL} hold. Then
		$\E (W_\infty (\theta))^p<\infty$
		if and only if
		\begin{align}\label{Hp}
			\kappa(p\theta)< p\kappa(\theta), \mbox{ and } \bE(\langle e_{\theta},\cP\rangle^p)<\infty.
		\end{align}
	\end{theorem}
	
		In the following theorem, we
		give
		a sufficient condition for $\E(W_\infty(\theta)^pH(W_\infty(\theta)))<\infty$, where $H$ is a slowly varying function.
		
		\begin{theorem}\label{Theorem:WH(W)}
			Let $\theta\in\Theta$ and $p>1$. Suppose that  \eqref{kappa'} hold.
			Assume that  $H\in\mathbb{L}$, and if
			$p$ is an integer, we further assume $H\in \mathbb{L}^*\cup \tilde{\mathbb{L}}$.
			Then $\E(W_\infty(\theta)^pH(W_\infty(\theta)))<\infty$
			if (1) $p\theta\in\Theta_0$ and $\kappa(p\theta)<p\kappa(\theta)$; and (2)
			$\bE(\langle e_\theta,\mathcal{P}\rangle^pH(\langle e_\theta,\mathcal{P}\rangle)<\infty$.
		\end{theorem}

Before we prove Theorem \ref{Theorem:WH(W)}, we first give some estimates on the Laplace transform of $W_\infty(\theta)$.  Letting $t\to\infty$ in \eqref{decomp:Wt}, we get that
		\begin{align}\label{decomp:W2}
			W_\infty(\theta)=e^{-\kappa(\theta)\tau_1}e^{\theta z_\varnothing(\tau_1) }\sum_{j=1}^{N^\varnothing}e^{\theta S_i^\varnothing}W_{\infty}^j(\theta),
		\end{align}
		where $W_\infty^i(\theta)$, $j=1,\cdots,N^\varnothing$, are iid with the same law as $W_\infty(\theta)$ under $\P$, and are independent of $X_{\tau_1}$. Let  $\phi(\lambda)$  be the Laplace transform of  $W_\infty(\theta)$:
		$$
		\phi(\lambda):=\E(e^{-\lambda W_\infty(\theta)}),\quad \lambda\ge0.$$
		By \eqref{decomp:W2}, we have that
		\begin{align*}
			\phi(\lambda)&=\E\left(\prod_{i=1}^{N^\varnothing}\phi\left(\lambda  e^{-\kappa(\theta)\tau_1}e^{\theta z_\varnothing(\tau_1)} e^{\theta S_i^\varnothing}\right)\right)=\bE\left(\prod_{i=1}^N\phi\left(\lambda e^{-\kappa(\theta)\tau}e^{\theta \xi_\tau} e^{\theta S_i}\right)\right).
					\end{align*}
 Recall that $Y_\theta=e^{-\kappa(\theta)\tau}e^{\theta \xi_\tau}$, see \eqref{def-Xi'}.
		For $\lambda\ge 0$, define
		\begin{align}\label{def-I}
			I(\lambda):=\sum_{i=1}^N -\log \phi \left(\lambda Y_\theta e^{\theta S_i}\right).
		\end{align}
		Then we have
		\begin{align}\label{phi}
			\phi(\lambda)=\bE\Big(\exp\left\{-I(\lambda )\right\}\Big), \quad \lambda\ge0.
		\end{align}

		If $\E(W_\infty(\theta))^n<\infty$, define
		\begin{align}\label{def-phi-n}
			\phi_n(\lambda):=\E \mathcal{T}_n(\lambda W_\infty(\theta))=(-1)^{n+1}\Big(\phi(\lambda)-1-\sum_{k=1}^{n}\mu_k\lambda^k\Big),
		\end{align}
		where $\mu_k:=(-1)^k\E[(W_\infty(\theta))^k]/k!$.
\eqref{def-phi-n} can be rewritten as
\begin{align}\label{def-phi-n'}
			\phi(\lambda)=1+\sum_{k=1}^{n}\mu_k\lambda^k+(-1)^{n+1}\phi_n(\lambda).
		\end{align}
		By \eqref{est:T} and the dominated convergence theorem, we  have
		\begin{align}\label{est-phi-n}
			0\le \phi_n(\lambda)\le \E\left(\frac{2(W_\infty(\theta))^n}{n!}\wedge\frac{(W_\infty(\theta))^{n+1}\lambda}{(n+1)!}\right)\lambda^n=o(\lambda^n),\quad \lambda\to0.
		\end{align}
		
		Recall that
		$	\Xi_\theta=Y_\theta\langle e_{\theta},\cP\rangle,$ see \eqref{def-Xi'-2}.
			Since $\tau $ and $\xi$ are independent, we get that, if $\beta+r\kappa(\theta)-\varphi(r\theta)>0$, then
		\begin{align}\label{y-Ly}
			\bE(Y_\theta^r)&=\int_0^\infty \beta e^{-\beta s}e^{-r\kappa(\theta) s}e^{\varphi(r\theta)s}ds=\frac{\beta}{\beta+r\kappa(\theta)-\varphi(r\theta)}\\
			&=\frac{\beta}{(r-1)\kappa(\theta)+\beta\gamma(\theta)},
		\end{align}
		where in the last equality we used \eqref{def-kappa}.
		Now we  give a crucial lemma that will play an important role in the proof  of Theorems \ref{Theorem:WH(W)}.

		\begin{lemma}\label{lem:psi2}
			Suppose $\theta\in\Theta$ and that \eqref{kappa'}--\eqref{LLogL} hold. Let $n\ge1$ be an integer. If $\E(W_\infty(\theta))^n<\infty$, then there exists a constant $C_n$ such that for all $\lambda>0$,
			\begin{align}\label{psi2}
				\left|\phi_n(\lambda)-\bE \Big(\sum_{i=1}^{N} \phi_n(\lambda  Y_\theta e^{\theta S_i})\Big)\right|
				\le C_{n}\bE\left[\Big(\Xi_\theta^n\lambda^n\Big)\wedge \Big(\Xi_\theta ^{n+1} \lambda^{n+1}\Big)\right]<\infty.
			\end{align}
		\end{lemma}
		{\bf Proof:}
		We claim that
		\begin{equation}\label{4.19}
			\bE(\Xi_\theta^n)=\bE(Y_\theta^n)\bE[\langle e_{\theta},\cP\rangle^n]<\infty.
		\end{equation}
		In fact, for $n=1$, using $\beta+\kappa(\theta)-\varphi(\theta)=\beta\chi(\theta)>0$ and \eqref{y-Ly}, we get
		\begin{align}
			\bE(\Xi_\theta)=\E(Y_\theta)\bE\langle e_{\theta},\cP\rangle=1.
		\end{align}
		For $n\ge 2$, since $\E(W_\infty(\theta))^n<\infty$,
we have by Theorem \ref{main-Lp} that
$n\kappa(\theta)>\kappa(n\theta)$ and $\bE\langle e_{\theta},\cP\rangle^n<\infty$.  Note that $n\kappa(\theta)-\beta-\varphi(n\theta)=n\kappa(\theta)-\kappa(n\theta)+\beta\chi(n\theta)>0$. Thus we have by \eqref{y-Ly} that $\E(Y_\theta^n)<\infty$.
Therefore \eqref{4.19} holds.
		
		Recall the definition of $I(\lambda)$ in \eqref{def-I}.
		By Jensen's inequality,  we have $-\log \phi(\lambda)\le\E( \lambda W_\infty(\theta))=\lambda$.
		Consequently,
		\begin{align}\label{4.1}
			I(\lambda)\le \lambda \Xi_\theta, \quad \lambda\ge0.
		\end{align}
		Thus $\bE(I(\lambda)^n)<\infty$ by \eqref{4.19}.
		It follows from \eqref{phi} that
		\begin{align}\label{phi2}
			\phi(\lambda)=\bE(e^{-I(\lambda)})=(-1)^{n+1}\bE(\mathcal{T}_n(I(\lambda)))+1+\sum_{l=1}^n (-1)^l\frac{\bE(I(\lambda)^l)}{l!}.
		\end{align}
Since $\mathcal{T}_n$ is an increasing function on $(0, \infty)$
and $I(\lambda)\le \lambda \Xi_\theta,$
		we get by \eqref{est:T} that
		\begin{align}\label{4.22}
			0\!\le \!\bE(\mathcal{T}_n(I(\lambda)))\!\le\! \bE(\mathcal{T}_n(\lambda \Xi_\theta))\!\le \!\bE \Big(\!\frac{2\Xi_\theta^n\lambda^n}{n!}\!\wedge \! \frac{\Xi_\theta^{n+1}\lambda^{n+1}}{(n+1)!}\!\Big)=
					o(\lambda^n), \quad \lambda\to0.
		\end{align}

		Now we deal with $\bE(I(\lambda)^l)$,
for $1\le l\le n$.
We claim that there exist constants $a_k, k=0,\cdots, n$, such that
		\begin{align}\label{L-phi}
			-\log\phi(\lambda)=\sum_{l=1}^n a_l\lambda^l
			-(-1)^{n+1}\phi_n(\lambda) +g_n(\lambda),
		\end{align}
			where
		\begin{equation}\label{4.9}
			|g_n(\lambda)|\le c_n(\lambda^n\wedge \lambda^{n+1}).
		\end{equation}
In fact, by Taylor's formula, we have
		$$
		-\log(1-x)=\sum_{k=1}^n \frac{x^k}{k}+r_n(x), \quad 0\le x<1,
		$$
		where \begin{align}\label{est-R}
			0\le r_n(x)=\sum_{k=n+1}^\infty \frac{x^k}{k}\le \sum_{k=n+1}^\infty x^k=\frac{x^{n+1}}{1-x}.
		\end{align}
By \eqref{def-phi-n'},
 we have that  $1-\phi(\lambda)=-\sum_{l=1}^n\mu_l \lambda^l-(-1)^{n+1}\phi_n(\lambda).$ Thus
 we have for all $\lambda>0$,
		\begin{align*}
			-\log\phi(\lambda)=&-\log(1-(1-\phi(\lambda))=\sum_{k=1}^n \frac{(1- \phi(\lambda))^k}{k}+r_n(1-\phi(\lambda))\\
			=&-\sum_{r=1}^n\mu_r \lambda^r-(-1)^{n+1}\phi_n(\lambda)  +\sum_{k=2}^n\frac{(-\sum_{r=1}^n\mu_r \lambda^r-(-1)^{n+1}\phi_n(\lambda)))^k}{k}\\
&+r_n(1-\phi(\lambda)).
		\end{align*}
		It follows from \eqref{est-R} that for $\lambda<1$,
		\begin{align}\label{r-n}
			r_n(1-\phi(\lambda))\le \frac{(1-\phi(\lambda))^{n+1}}{\phi(\lambda)}\le \frac{\lambda^{n+1}}{\phi(1)},
		\end{align}
		where in the last inequality, we used the fact that $1-\phi(\lambda)\le \lambda$ and $\phi(\lambda)\ge \phi(1)>0$ for $\lambda<1$.
By \eqref{est-phi-n}, we have $0<\phi_n(\lambda)\le 2|\mu_n|\lambda^n$. Thus by \eqref{r-n},
there exist constants $b_2,\cdots, b_n$ and $g_n(\lambda)\sim o(\lambda^n)$ as $\lambda\to0$ such that
		\begin{align*}
			\sum_{k=2}^n\frac{(-\sum_{r=1}^n\mu_r \lambda^r-(-1)^{n+1}\phi_n(\lambda)))^k}{k}+r_n(1-\phi(\lambda))=\sum_{r=2}^n b_r\lambda^r+g_n(\lambda).
		\end{align*}
		Note that for $\lambda < 1$,
		$|g_n(\lambda)|\le c_n\lambda^{n+1}$.
		Since $-\log \phi(\lambda)\le \lambda$ and $\phi_n(\lambda)\le 2|\mu_n|\lambda^n$,
		we have for $\lambda\ge 1$,
		\begin{equation*}
			|g_n(\lambda)|=\left|-\log\phi(\lambda)-\sum_{l=1}^n a_l \lambda^l
			+(-1)^{n+1}\phi_n(\lambda) \right|\le c_n\lambda^n,
		\end{equation*}
which implies \eqref{4.9}. Hence \eqref{L-phi} holds.
		
		(1) For $l=1$, by the definition of $I(\lambda)$ and \eqref{L-phi}, we have
	\begin{align*}
	I(\lambda)=&\sum_{i=1}^N -\log\phi\left(\lambda Y_\theta e^{\theta S_i}\right)\\
=&\sum_{r=1}^n a_r Y_\theta^r\langle e_{r\theta},\cP\rangle \lambda^r -(-1)^{n+1}\sum_{i=1}^N \phi_n(\lambda Y_\theta e^{\theta S_i})+\sum_{i=1}^N g_n(\lambda Y_\theta e^{\theta S_i}).
\end{align*}
		Thus
		\begin{align}\label{E-I}
			\bE(I(\lambda))
=&\sum_{r=1}^n a_r \bE(Y_\theta^r)\bE(\langle e_{r\theta},\cP\rangle) \lambda^r -(-1)^{n+1}\bE\left[\sum_{i=1}^N \phi_n(\lambda Y_\theta e^{\theta S_i})\right]\nonumber\\
&+\bE\left[\sum_{i=1}^N g_n(\lambda Y_\theta e^{\theta S_i})\right].
		\end{align}
		By \eqref{4.9} we have
		\begin{align}\label{g_n}
			\left|\bE\left[\sum_{i=1}^N g_n(\lambda Y_\theta e^{\theta S_i})\right]\right|\le c_n\bE \left[\Big(\Xi_\theta^{n} \lambda^{n}\Big)\wedge \Big(\Xi_\theta^{n+1} \lambda^{n+1}\Big)\right].
		\end{align}
		
		(2) For $l=2,3,\cdots,n$,
		\begin{align}\label{4.10}
			I(\lambda)^l &=\sum_{m=l}^{n} \sum_{r_1+r_2+\cdots+r_l=m}a_{r_1}\cdots a_{r_l}\langle e_{r_1\theta},\cP\rangle\cdots \langle e_{r_l\theta},\cP\rangle Y_\theta^m\lambda^m+ I_{l,n}(\lambda),
		\end{align}
		where
		\begin{align*}
			I_{l,n}(\lambda)=&\sum_{m=n+1}^{nl} \sum_{r_1+r_2+\cdots+r_l=m}a_{r_1}\cdots a_{r_l}\langle e_{r_1\theta},\cP\rangle\cdots \langle e_{r_l\theta},\cP\rangle Y_\theta^m\lambda^m\\
			&+\sum_{k=1}^{l}{n\choose k} \Big(\sum_{i=1}^N g_n^*(\lambda Y_\theta e^{\theta S_i})\Big)^k\Big(\sum_{r=1}^n a_r Y_\theta^r\langle e_{r\theta},\cP\rangle \lambda^r\Big)^{l-k},
		\end{align*}
		and $g_n^*(\lambda )=	-(-1)^{n+1}\phi_n(\lambda) +g_n(\lambda)$.
		Note that  for $m\le n$ and $r_1+r_2+\cdots+r_l=m$,
		$$\bE(\langle e_{r_1\theta},\cP\rangle\cdots \langle e_{r_l\theta},\cP\rangle Y_\theta^m)\le \bE(\langle e_{\theta},\cP\rangle^m Y_\theta^m)=\bE(\Xi_\theta^m)<\infty.$$
		Then by \eqref{4.10} we get that
		\begin{align}\label{equa-I}
			\bE\left[I(\lambda)^l\right]= &\sum_{m=l}^{n} \sum_{r_1+r_2+\cdots+r_l=m}a_{r_1}\cdots a_{r_l}\bE\left(\langle e_{r_1\theta},\cP\rangle\cdots \langle e_{r_l\theta},\cP\rangle\right)\bE( Y_\theta^m )\lambda^m\nonumber\\
&+ \bE [I_{l,n}(\lambda)].
		\end{align}
		We claim that,  for $2\le l\le n$, there exists $c_{l,n}>0$ such that
		\begin{align}\label{4.7}
			|I_{l,n}(\lambda)|
			\le c_{l,n}\left[\Big(\Xi_\theta^{n} \lambda^{n}\Big)\wedge \Big(\Xi_\theta^{n+1} \lambda^{n+1}\Big)\right].
		\end{align}
		In fact, by \eqref{est-phi-n} and \eqref{4.9}, $|g^*(\lambda)|\le c_n\lambda^n$. Thus if $\lambda \Xi_\theta<1$,  for $l=2,\cdots,n$, we have
		\begin{align*}
			|I_{l,n}(\lambda)|\le &c_{l,n}\Xi_\theta^{n+1} \lambda^{n+1}
			+\sum_{k=1}^{l}c_{l,n,k} (\lambda\Xi_\theta)^{nk} (\lambda\Xi_\theta)^{l-k}
			\le c_{l,n}\Xi_\theta^{n+1} \lambda^{n+1},
		\end{align*}
		where in the last inequality we used the inequality $nk+l-k\ge n-1+l\ge n+1$.
		On the other hand, if $\lambda \Xi\ge 1$, by \eqref{4.10} and \eqref{4.1}, we have that for $2\le l\le n$,
		\begin{align*}
			|I_{l,n}(\lambda)|&\le I(\lambda)^l+\sum_{m=l}^{n} \sum_{r_1+r_2+\cdots+r_l=m}|a_{r_1}\cdots a_{r_l}|\langle e_{r_1\theta},\cP\rangle\cdots \langle e_{r_l\theta},\cP\rangle Y_\theta^m\lambda^m\\
			&\le \Xi_\theta^{l} \lambda^{l}+c_{l,n}\Xi_\theta^{n} \lambda^{n}\le c_{l,n}\Xi_\theta^{n} \lambda^{n}.
		\end{align*}
		Thus \eqref{4.7} is valid.

		Combining \eqref{phi2}, \eqref{4.22},  \eqref{E-I},\eqref{g_n}, \eqref{equa-I} and \eqref{4.7}, we get that there exist  $\mu_l', l=1,\cdots,n$, such that
		\begin{align}\label{phi3}
			\phi(\lambda)=1+\sum_{l=1}^n \mu_l' \lambda^l+(-1)^{n+1}\bE\left[\sum_{i=1}^N \phi_n(\lambda Y_\theta e^{\theta S_i})\right]+R_{n}(\lambda),
		\end{align}
		where
		\begin{align} \label{Rn1}
			|R_{n}(\lambda)|\le c_{n}\bE\left[\Big(\Xi_\theta^{n} \lambda^{n}\Big)\wedge \Big(\Xi_\theta^{n+1} \lambda^{n+1}\Big)\right]=o(\lambda^n), \quad \mbox{ as } \lambda\to0.
		\end{align}
		Since $\phi_n(\lambda)\le 2|\mu_n|\lambda^n$ and $\phi_n(\lambda)=o(\lambda^n)$, we have by the dominated convergence theorem that
		\begin{align}\label{lim-phi-n}
			\bE\left[\sum_{i=1}^N \phi_n(\lambda Y_\theta e^{\theta S_i})\right]=o(\lambda^n), \quad \mbox{ as } \lambda\to0.
		\end{align}
Now comparing  \eqref{phi3} with  \eqref{def-phi-n'}, and using \eqref{Rn1} and \eqref{lim-phi-n},  we  deduce
 that $\mu_l=\mu_l'$, $l=1,\cdots,n$, and
		\begin{align*}
\phi_n(\lambda)=
\bE\left[\sum_{i=1}^N \phi_n(\lambda Y_\theta e^{\theta S_i})\right]+(-1)^{n+1}R_{n}(\lambda).
		\end{align*}
The desired result now follows immediately from \eqref{Rn1}.
		
		\hfill$\Box$

Now we define a  measure $\mu$ on $\R^+$ as follows:
for any Borel subset of $\R^+$,
		\begin{align}\label{def:mu}
			\mu(B):=\bE \Big(\sum_{i=1}^N{\bf 1}_B(Y_\theta e^{\theta S_i})Y_\theta e^{\theta S_i}\Big).
		\end{align}
		Note that by \eqref{y-Ly},
		$\mu(\R^+)
		=\bE(Y_\theta)\chi(\theta)=1$.  Thus $\mu$ is a probability measure on $\R^+$.  In the remainder of this section, we assume that $Z,Z_1,Z_2,\cdots$ are i.i.d. random variables with common distribution $\mu$.
		
		Put
$$\psi_n(\lambda)=\phi_n(\lambda)/\lambda.$$
 By the definition  of $\mu$,  we have that
		$$\bE \Big(\sum_{i=1}^N\phi_n(\lambda Y_\theta e^{\theta S_i})\Big)=\lambda \bE \Big(\sum_{i=1}^N\psi_n(\lambda Y_\theta e^{\theta S_i})Y_\theta e^{\theta S_i}\Big)=\lambda E(\psi_n(\lambda Z)).$$
		It follows from  Lemma \ref{lem:psi2} that
		\begin{align*}
			|\psi_n(\lambda)-E(\psi_n(\lambda Z))|&=\lambda^{-1}|\phi_n(\lambda)-\bE \Big(\sum_{i=1}^N\phi_n(\lambda Y_\theta e^{\theta S_i})\Big)|\\
			&\le C_{n}\bE\left[\Big(\Xi_\theta^{n} \lambda^{n-1}\Big)\wedge \Big(\Xi_\theta^{n+1} \lambda^{n}\Big)\right] .
		\end{align*}

		Using \eqref{y-Ly}, we can easily get the following result.
		
		\begin{lemma}\label{lemma:Z}
			Let $\theta\in\Theta$  and  $r\ge  0$.
If $\kappa(r\theta)<\infty$ and $r\kappa(\theta)-\kappa(r\theta)+\beta\chi(r\theta)>0$, then
			\begin{align}\label{Z-r}
				E(Z^{r-1})=\bE(Y_\theta^{r})\chi(r\theta)=\frac{\beta\chi(r\theta)}{r\kappa(\theta)-\kappa(r\theta)+\beta\chi(r\theta)}.
			\end{align}
			Moreover if $r\kappa(\theta)>\kappa(r\theta)$, then $E(Z^{r-1})<1$.
		\end{lemma}

		\begin{lemma}\label{exp:psi-n}
	Suppose $\theta\in\Theta$ and that \eqref{kappa'}--\eqref{LLogL} hold. Assume  $\E(W_\infty(\theta))^n<\infty$ for an integer $n\ge 1$. When $n=1$, we further assume  that there exists $\theta_0>\theta$ such that $\kappa(\theta_0)<\infty$.
		Then for any $k\ge1$,
			$E\psi_n(\lambda Z_1Z_2\cdots Z_k)<\infty$  and
			\begin{align}\label{phi6}
				\psi_n(\lambda)=\sum_{k=1}^\infty\delta_{n,k}(\lambda),
			\end{align}
			where
			\begin{align}\label{def:deltan1}
			\delta_{n,1}(\lambda)=\psi_n(\lambda)-E(\psi_n(\lambda Z_1))
			\end{align}
			and for $k\ge 2$,
			\begin{align}\label{def:deltank}
				\delta_{n,k}(\lambda):=E(\psi_n(\lambda Z_1Z_2\cdots Z_{k-1}))-E(\psi_n(\lambda Z_1Z_2\cdots Z_k))=E(\delta_{n,k-1}(\lambda Z)).
			\end{align}
		\end{lemma}
		{\bf Proof}
		For  $n\ge2, $ since $0<\E(W_\infty(\theta))^n<\infty$, we have $\kappa(n\theta)<n\kappa(\theta)$. It follows from  Lemma \ref{lemma:Z}  that $E(Z^{n-1})<1$.  Using the fact $0\le \psi_n(\lambda)\le 2|\mu_n|\lambda^{n-1}$,  we get that
		\begin{align}\label{2.1}
			E\psi_n(\lambda  Z_1Z_2\cdots Z_k)\le 2|\mu_n|\lambda^{n-1}[ E(Z^{n-1})] ^k.
		\end{align}
		Thus, for any $k\ge 1$, $\delta_{n,k}$ is well-defined and
		\begin{align}\label{phi5}
			\psi_n(\lambda)=E(\psi_n(\lambda Z_1Z_2\cdots Z_k))+\sum_{j=1}^k\delta_{n,j}(\lambda).
		\end{align}
		Letting $k\to\infty$, using \eqref{2.1} and the fact $E(Z^{n-1})<1$, we get the desired result immediately.
		
		Now we deal with the case $n=1$. Note that
		$0\le \psi_1(\lambda)\le 2\E(W_\infty(\theta))=2$.
Thus $E\psi_1(\lambda Z_1Z_2\cdots Z_k)<\infty$ for any $k\geq 1$.
To finish the proof, it suffices to show that
		\begin{align}\label{claim3}
			\lim_{k\to\infty}E(\psi_1(\lambda Z_1Z_2\cdots Z_k))= 0.
		\end{align}
		Since  $\theta\kappa'(\theta)<\kappa(\theta)$,
		we know that $\kappa(s)/s$ is decreasing on $[\theta,\theta_0]$. Thus there exists $r>1$ such that  $\kappa(r\theta)<r\kappa(\theta)$.
		Now using Lemma \ref{lemma:Z} we get $E(Z^{r-1})<1$. Applying Jensen's inequality,  we have
		$$
		E(\log Z)=(r-1)^{-1}E(\log Z^{r-1})\le (r-1)^{-1}\log E(Z^{r-1})<0.
		$$
		Hence using the law of large numbers, we get that $\lim_{k\to\infty}\sum_{j=1}^k \log Z_j=-\infty$  which implies that
		\begin{align}\label{limz2}
			\lim_{k\to\infty}Z_1Z_2\cdots Z_k=0.
		\end{align}
		Note that  $\psi_1(\lambda)$ is bounded and $\psi(0+)=0$. Using \eqref{limz2} and the dominated convergence theorem,  we get  that $E(\psi_1(\lambda Z_1Z_2\cdots Z_k))\to 0$, as $k\to\infty$.
		
		The proof is now complete.
		\hfill $\Box$

		\bigskip

		If
		$\theta p\in\Theta_0 $ and $\kappa(p\theta)<p\kappa(\theta)
		$, then  by the continuity of $\kappa$ on $\Theta_0$, there exists $r>p$ such that $\kappa(r\theta)<r\kappa(\theta)$.
		Thus $E(Z^{r-1})<1$ by Lemma \ref{lemma:Z}.
By the dominated convergence theorem, we have that
		$$\lim_{\delta\to0} E(Z^{p+\delta-1} \vee Z^{p-\delta-1})=E(Z^{p-1})<1.$$
		So we can always choose $\delta\in(0, (n+1-p)\wedge (p-n))$ if $p\in(n,n+1)$ and  $\delta\in(0,1)$ if $p$ is an integer, such that $p+\delta<r$
		and
		\begin{align}\label{3.8'}
			&E(Z^{p+\delta-1} \vee Z^{p-\delta-1})\nonumber\\
			=&\bE\Big(\sum_{i=1}^N Y_\theta^{p-\delta}e ^{(p-\delta)\theta S_i}{\bf 1}_{Y_\theta e ^{\theta S_i}<1}\Big)+\bE\Big(\sum_{i=1}^N Y_\theta^{p+\delta}e ^{(p+\delta)\theta S_i}{\bf 1}_{Y_\theta e ^{\theta S_i}\ge 1}\Big)<1.
		\end{align}
		It is clear  that
		\begin{align}\label{3.1}
			\bE(Y_\theta^{p+\delta}) \chi((p+\delta)\theta)=E(Z^{p+\delta-1})<1
\end{align}
and
\begin{align}\label{3.1'}
\bE(Y_\theta^{p-\delta}) \chi((p-\delta)\theta)=E(Z^{p-\delta-1})<1.
		\end{align}
Thus, $\bE(Y_\theta^{p+\delta})<\infty$ and $\bE(Y_\theta^{p-\delta})<\infty$.
In the remainder of this section, $\delta$ always stands for such a constant.

		\begin{lemma}\label{lemma2}
			Let $\theta\in\Theta$,   $p>1$ and  $h: \R_+\to\R$ be a bounded Borel function.
			Assume  that $p\theta\in\Theta_0$ and $\kappa(p\theta)<p\kappa(\theta)$.
			If  there exists a slowly varying function $L\in\mathbb{L_\delta}$ such that
			$|h(\lambda)|\le \lambda^{p-1}L(\lambda^{-1})$ for all $\lambda>0$, then
			\begin{align}\label{3.12}
				E(|h|(\lambda Z))\le A\lambda^{p-1}L(\lambda^{-1}), \quad \lambda>0,
			\end{align}
			for some constant $A\in(0,1)$.
			In addition, if the limit $\lim_{\lambda\to0}\frac{h(\lambda)}{\lambda^{p-1}L(\lambda^{-1})}=a$ exists,  then
			\begin{align}\label{3.11}
				\lim_{\lambda\to0}\frac{E(h(\lambda Z))}{\lambda^{p-1}L(\lambda^{-1})}=a \frac{\beta \,\chi(p\theta)}{p\kappa(\theta)-\kappa(p\theta)+\beta\chi(p\theta)}.
			\end{align}
		\end{lemma}
		{\bf Proof:}
		Assume that  $L\in\mathbb{L_\delta}$ satisfies the assumption of the lemma. Then we have
		\begin{align}\label{3.81}
			\frac{|h|(\lambda Z)}{\lambda^{p-1}L(\lambda^{-1})}\le Z^{p-1}\frac{L(\lambda^{-1} Z^{-1})}{L(\lambda^{-1})}\le (Z^{p-1+\delta}\vee Z^{p-1-\delta}).
		\end{align}
By \eqref{3.8'},
$A=E(Z^{p-1+\delta}\vee Z^{p-1-\delta})<1$.  Thus \eqref{3.12} is valid.
		
		If $\lim_{\lambda\to0}\frac{h(\lambda)}{\lambda^{p-1}L(\lambda^{-1})}=a$, then by \eqref{3.81} and the dominated convergence theorem, we have
		$$\lim_{\lambda\to0}\frac{E(h(\lambda Z))}{\lambda^{p-1}L(\lambda^{-1})}=aE(Z^{p-1})=a\frac{\beta \,\chi(p\theta)}{p\kappa(\theta)-\kappa(p\theta)+\beta\chi(p\theta)},$$
where in the last equality we used \eqref{Z-r}.
The proof is now complete.
		\hfill$\Box$
		
		\bigskip

		\begin{lemma}\label{lem:psito0-5}
			Suppose $\theta\in\Theta$ and \eqref{kappa'} holds. Let  $n\ge1$ be an integer and
			$p\in[n,n+1]\cap (1, \infty)$.
			Let $L\in\mathbb{L}$  and let $H$ defined by \eqref{def:H}. 	
			If (1) $p\theta\in\Theta_0$  and $\kappa(p\theta)<p\kappa(\theta)$; and (2)
			$\bE(\langle e_\theta,\mathcal{P}\rangle^pH(\langle e_\theta,\mathcal{P}\rangle))<\infty$,
			then $$\lim_{\lambda\to 0}\lambda^{1-p}H(1/\lambda)\psi_n(\lambda)=0.$$
		\end{lemma}
		{\bf Proof:}  Since $p>1$ and  $\kappa(p\theta)<p\kappa(\theta)$,
by Remark \ref{rem:kappa} (1),
we have
		$\theta\kappa'(\theta)<\kappa(\theta)$.
Since 	$\bE(\langle e_\theta,\mathcal{P}\rangle^pH(\langle e_\theta,\mathcal{P}\rangle))<\infty$, we know 	$\bE(\langle e_\theta,\mathcal{P}\rangle^q)<\infty$, for any $1<q\le p$.
Therefore \eqref{LLogL} holds, and consequently $W_\infty(\theta)$ is non-degenerate.
		
		If $n=1$, then $\E(W_\infty(\theta))=1$. Since $p>1$ and  $\kappa(p\theta)<p\kappa(\theta)$,
by Remark \ref{rem:kappa} (1),
we have $\kappa(r\theta)<r\kappa(\theta)$ for
all $1<r< p$.
Thus, if   $n\ge 2$,  then $\kappa(n\theta)<n\kappa(\theta)$.
	Consequently $\E(W_\infty(\theta))^n<\infty$.

	Since $ \bE\left( \langle e_\theta,\mathcal{P}\rangle^{p}H(\langle e_\theta,\mathcal{P}\rangle)\right)<\infty$ and $\bE(Y_\theta^{p+\delta})<\infty$, we have that
	\begin{align}\label{moment-xi-finite}
		\bE\left( \Xi_\theta^{p}H(\Xi_\theta)\right) <\infty.
	\end{align}
		
		We first consider the case $p\in(n,n+1)$.
		Without loss of generality, we can assume that $H\in \mathbb{L}_\delta$. Then we have that
		By Lemma \ref{exp:psi-n},
		\begin{align}\label{expansion-psi}
			\psi_n(\lambda)=\sum_{k=1}^\infty\delta_{n,k}(\lambda),
		\end{align}
		where $\delta_{n, k}$ are defined in \eqref{def:deltan1} and \eqref{def:deltank}.
		We have by Lemma \ref{lem:psi2} that
		\begin{align} \label{5.71}
			|\delta_{n,1}(\lambda)|\le C_{n}\bE\Big((\lambda^{n}\Xi_\theta^{n+1})\wedge (\lambda^{n-1}\Xi_\theta^{n})\Big).
		\end{align}
		Since $H\in\mathbb{L}_\delta$, we have by Lemma \ref{inequality} in the  Appendix that
		\begin{align}\label{5.33}
			\lambda^{n}\Xi_\theta^{n+1}\wedge
			\lambda^{n-1}\Xi_\theta^{n}&\le \lambda^{p-1}\Xi_\theta^{p}H(\Xi_\theta)/H(1/\lambda).
		\end{align}		
		Thus, by \eqref{5.71},\eqref{5.33} and the dominated convergence theorem, we have
		\begin{equation}\label{5.31}
			\lim_{\lambda\to 0}\lambda^{1-p}H(1/\lambda)|\delta_{n,1}(\lambda)|=0
		\end{equation}
		and
		\begin{equation}\label{5.32}
			\lambda^{1-p}H(1/\lambda)|\delta_{n,1}(\lambda)|\le C_{n}\bE(\Xi_\theta^{p}H(\Xi_\theta)).
		\end{equation}
		Using Lemma \ref{lemma2} and  induction, we get that
		\begin{align}\label{deltan-5}
			\lim_{\lambda\to 0}\lambda^{1-p}H(1/\lambda)|\delta_{n,k}(\lambda)|=0
		\end{align}
		and
		\begin{equation}\label{5.2}
			\lambda^{1-p}H(1/\lambda)|\delta_{n,k}(\lambda)|\le C_{n}\bE(\Xi_\theta^{p}L(\Xi_\theta))A^{k-1}<\infty,
		\end{equation}
		where $A\in(0,1)$ is a constant.
		Thus, by \eqref{expansion-psi},  \eqref{deltan-5},  \eqref{5.2} and  the dominated convergence theorem, we have
		\begin{align*}
			\lim_{\lambda\to 0}\lambda^{1-p}H(1/\lambda)\psi_n(\lambda) = 0.
		\end{align*}
		
		For $p=n$, $H(x)=\int_1^{x\vee 1} L(t)/t dt$.
		Since $H$ is slowly varying at $\infty$, there exists $x_0>1$ such that $x^{-1}H(x)$ is decreasing on $[x_0,\infty)$.
		Let $\hat{H}(x)=H(x\vee x_0)$. Then
		$\hat{H}$ is non-decreasing and $\hat{H}(x)/x$ is non-increasing on $(0,\infty)$. Applying Lemma \ref{inequality} in Appendix with $p=n$ and $L$ replaced by $\hat{H}$, we get \eqref{5.33} with $H$ replaced by $\tilde{H}$.
		Using an argument similar to the case $p\in(n,n+1)$ above, we can show that the desired result holds for $p=n$ as well.

		The case  $p=n+1$ can be dealt with similarly and we omit the details.
		\hfill$\Box$
		\bigskip

		Recall that $\phi_n(\lambda)=\E \mathcal{T}_n(\lambda W_\infty(\theta))$ and $\psi_n(\lambda)=\phi_n(\lambda)/\lambda$.

		\begin{proposition}\label{prop:WH(W)}
			Under the conditions of Lemma \ref{lem:psito0-5}, we have
			$$\int_0^1 \psi_n(\lambda)L(\lambda^{-1})\lambda^{-p}\,d\lambda<\infty.$$
		\end{proposition}
		{\bf Proof}
		Without loss of generality, we assume that $L\in\mathbb{L}_\delta$. It follows from  \eqref{3.8'}  that
		$$
		E \Big(Z^{p-1}(Z^\delta\vee Z^{-\delta}))<1.
		$$
		By Lemma \ref{lem:psi2}, we have that
		\begin{align}\label{5.21}
			&\int_{0}^{1}|\psi_n(\lambda)-E(\psi_n(\lambda Z))|L(1/\lambda)\lambda^{-p}\,d\lambda\nonumber\\
			\le &C_{n}\left(\bE\left[ \Xi_\theta^{n+1}\int_{0}^{\Xi_\theta^{-1}\wedge 1}L(1/\lambda)\lambda^{-(p-n)}\,d\lambda \right]+\bE\left[ \Xi_\theta^{n}{\bf 1}_{\Xi_\theta>1}\int_{\Xi_\theta^{-1}}^1 L(1/\lambda)\lambda^{-(1+p-n)}\,d\lambda\right]\right).
		\end{align}
		For $p\in[n,n+1)$, we have that, for any $x>0$,
		\begin{align*}
			\int_0^x L(1/\lambda)\lambda^{-(p-n)}\,d\lambda\le \int_0^x L(1/x)(x/\lambda)^\delta \lambda^{-(p-n)}\,d\lambda=(1+n-p-\delta)^{-1}L(1/x) x^{1+n-p},
		\end{align*}
		where in the inequality,
		we use the fact that $L(1/\lambda)\le L(1/x) (x/\lambda)^\delta$ for $\lambda<x$.
		For $p=n+1$, we have
		\begin{align*}
			\int_0^x L(1/\lambda)\lambda^{-(p-n)}\,d\lambda=\int_{x^{-1}}^\infty L(\lambda)\lambda^{-1}\,d\lambda=\tilde{L}(1/x).
		\end{align*}
		Similar as above, we get that for $0<x<1$,
		\begin{align}\label{5.5}
			\int_{x}^1 L(1/\lambda)\lambda^{-(1+p-n)}\,d\lambda\le \left\{\begin{array}{cc}
				(p-n-\delta)^{-1}L(1/x)x^{n-p}, & \mbox{ if }  p\in(n,n+1];\\
				L^*(1/x),&\mbox{ if }  p=n.
			\end{array} \right.
		\end{align}
		Plugging the last three displays into \eqref{5.21}, and using  \eqref{moment-xi-finite} and  \eqref{compare2}, we get
		\begin{align}\label{5.9}
			&\int_{0}^{1} |\psi_n(\lambda)-E\psi_n(\lambda Z)|L(1/\lambda)\lambda^{-p}\,d\lambda
			\le c_{n,p}\bE(\Xi_\theta^{p}(L(\Xi_\theta)+H(\Xi_\theta))< \infty.
		\end{align}
		Since $\psi_n(\lambda)\le 2|\mu_n|\lambda^{n-1}$, we have  for any $a\in(0,1)$,
		\begin{align*}
			\int_a^1  \psi_n(\lambda)L(1/\lambda)\lambda^{-p}\,d\lambda\le 2|\mu_n|\int_a^1 L(1/\lambda)\lambda^{n-p-1}\,d\lambda<\infty,
		\end{align*}
		where in the last inequality we used \eqref{5.5}.
		Similarly, $\int_a^1 E \psi_n(\lambda Z)L(1/\lambda)\lambda^{-p}\,d\lambda <\infty$. It is clear that
		\begin{align}\label{4.81}
		&\int_a^1  \psi_n(\lambda)L(1/\lambda)\lambda^{-p}\,d\lambda
		\nonumber\\
		\le&\int_a^1 E \psi_n(\lambda Z)L(1/\lambda)\lambda^{-p}\,d\lambda+\int_{0}^{1} |\psi_n(\lambda)-E\psi_n(\lambda Z)|L(1/\lambda)\lambda^{-p}\,d\lambda.
		\end{align}
		Using Fubini's theorem, we have that
		\begin{align}\label{5.7}
			&\int_a^1 E \psi_n(\lambda Z)L(1/\lambda)\lambda^{-p}\,d\lambda =E \Big(Z^{p-1}\int_{aZ}^Z  \psi_n(\lambda)L(Z/\lambda)\lambda^{-p}\,d\lambda\Big)\nonumber\\
			&\le E \Big(Z^{p-1}\int_{a}^1 \psi_n(\lambda)L(Z/\lambda)\lambda^{-p}\,d\lambda\Big)
			+E \Big({\bf 1}_{Z\le1}Z^{p-1}\int_{aZ}^a \psi_n(\lambda)L(Z/\lambda)\lambda^{-p}\,d\lambda\Big)\nonumber\\
			&\quad +E \Big({\bf 1}_{Z>1}Z^{p-1}\int_{1}^Z \psi_n(\lambda)L(Z/\lambda)\lambda^{-p}\,d\lambda\Big)\nonumber\\
			&=:J_1(\lambda)+J_2(\lambda)+J_3(\lambda).
		\end{align}
		It follows from \eqref{L-epsilon} that
		\begin{align}\label{5.8}
			J_1(\lambda)\le E \Big(Z^{p-1}(Z^\delta\vee Z^{-\delta})\Big)\int_{a}^1  \psi_n(\lambda)L(1/\lambda)\lambda^{-p}\,d\lambda.
		\end{align}
		Since $\psi_n$ is increasing and $L(Z/\lambda)\le L(1/a)(\lambda/aZ)^\delta$, we have
		\begin{align}\label{5.19}
			J_2(\lambda)&\le 	E \Big({\bf 1}_{Z\le1}\psi_n(a)L(1/a)Z^{p-1}\int_{aZ}^{\infty} (\lambda/(aZ))^{\delta}\lambda^{-p}\,d\lambda\Big)\nonumber\\
			&\le \frac{1}{p-1-\delta}\psi_n(a)L(1/a)a^{1-p}.
		\end{align}
		Combining \eqref{L-epsilon} with the fact $\psi_n(\lambda)\le 2|\mu_n|\lambda^{n-1}$, we get
		\begin{align}\label{5.10}
			J_3(\lambda)
			\le &2|\mu_n|L(1)E \Big({\bf 1}_{Z>1}Z^{p-1}\int_{1}^Z \lambda^{n-1} (Z/\lambda)^\delta \lambda^{-p}\,d\lambda\Big)\nonumber\\
			\le &2|\mu_n|L(1)E(Z^{p+\delta-1})\int_1^\infty \lambda^{-(1+p-n+\delta)}d\lambda\nonumber\\
			\le &2|\mu_n|L(1)\frac {1}{p+\delta-n}.
		\end{align}
		Combining \eqref{5.7}-\eqref{5.10} we get that
		\begin{align}\label{5.30}
			\int_a^1 E \psi_n(\lambda Z)&L(1/\lambda)\lambda^{-p}\,d\lambda\le E \Big(Z^{p-1}(Z^\delta\vee Z^{-\delta})\Big)\int_{a}^1 \psi_n(\lambda)L(1/\lambda)\lambda^{-p}\,d\lambda\nonumber\\
			&+\frac {2|\mu_n|L(1)}{p+\delta-n}+ \frac{1}{p-1-\delta}\psi_n(a)L(1/a)a^{1-p}.
		\end{align}
Combining this with \eqref{4.81}, we have
		\begin{align*}
			0&<\left(1-E \Big(Z^{p-1}(Z^\delta\vee Z^{-\delta})\Big)\right)\int_a^1\psi_n(\lambda)L(\lambda^{-1})\lambda^{-p}\,d\lambda\\
			&\le  \frac{1}{p-\delta-1}a^{1-p}L(a^{-1})\psi_n(a)+\frac {2|\mu_n|L(1)}{p+\delta-n}+ \int_0^1 |\psi_n(\lambda)-E(\psi_n(\lambda Z))|L(\lambda^{-1}) \lambda^{-p}\,d\lambda.
		\end{align*}
		Letting $a\to0$ and applying Lemma \ref{lem:psito0-5}, we get
		\begin{align*}
			&0<\left(1-E \Big(Z^{p-1}(Z^\delta\vee Z^{-\delta})\Big)\right)
			\int_0^1\psi_n(\lambda)L(\lambda^{-1})\lambda^{-p}\,d\lambda\\
			\le &\frac {2|\mu_n|L(1)}{p+\delta-n}+ \int_0^1 |\psi_n(\lambda)-E(\psi_n(\lambda Z))|L(\lambda^{-1}) \lambda^{-p}\,d\lambda<\infty.
		\end{align*}
		The desired result now follows immediately.
		
		\hfill$\Box$

		\noindent{\bf Proof of Theorem \ref{Theorem:WH(W)}:} Combining Lemma \ref{lemma:tb} with Proposition \ref{prop:WH(W)}, we immediately get the desired assertion.  \hfill$\Box$

		\subsection{Tail probability of $W_\infty(\theta)$}\label{sub:tail}
Now we study  the tail behavior of $W_\infty(\theta)$. More precisely, we will givee conditions for
\begin{align*}
	\P(W_\infty(\theta)>x) & \sim x^{-p}L(x),\quad x\to\infty,
\end{align*}	
with   $L\in\mathbb{L}$.

	If there exists $p_*>1$ such that $\kappa(\theta p_*)=p_*\kappa(\theta)$,
then $\E(W_\infty(\theta)^{p_*})=\infty$. Applying \cite[Theorem 2.2]{Liu2000} to the  branching random walk $\{X_{n}(\theta),n\ge 0\}$, we get  that  if
\begin{align}\label{condition:tail}
	\E\left(\sum_{u\in\mathcal{L}_1} (z_u(1)\vee0)e^{p_*\theta z_u(1)}\right)<\infty,\mbox{ and } \E(W_1(\theta)^{p_*})<\infty,
\end{align}
then the tail of
$W_\infty(\theta)$ is   of order $x^{-p_*}$.   Applying Theorem \ref{wtp}  and Lemma \ref{f-moment},  we can
rewrite conditions \eqref{condition:tail} in terms of $(\tau,\xi,\mathcal{P})$. As a consequence, we immediately get the following theorem.

\begin{proposition}\label{boundry case}
	Let $\theta\in\Theta$. Assume that there exists $p_*>1$ such that $\kappa(\theta p_*)=p_*\kappa(\theta)$.
	If
	$$
	\E(\langle e_\theta,\mathcal{P}\rangle^{p_*})<\infty,\,
	\int_{y>1} ye^{p_*\theta y}n(dy)<\infty \mbox{ and } \bE\left(\sum_{i=1}^N (S_i\vee 0)e^{p_*\theta S_i}\right)<\infty,$$
	then
	$$\P(W_\infty(\theta)>x)\sim c x^{-p_*}, \quad x\to\infty,$$
	for some constant $c\in(0,\infty)$.
\end{proposition}

\bigskip

 Now we deal with the case
$\kappa(p\theta)< p\kappa(\theta)$.

\begin{theorem}\label{Theorem:tailW}
Let $\theta\in\Theta$ and  $p\in(1,2)$
satisfy $p\theta\in\Theta_0$ and  $\kappa(p\theta)< p\kappa(\theta)$.
Assume that  \eqref{kappa'} holds and $L\in\mathbb{L}$.  Then
	\begin{align}\label{tail-W}
		\P(W_\infty(\theta)>x) & \sim \frac{\beta+p\kappa(\theta)-\varphi(p\theta)}{p\kappa(\theta)-\kappa(p\theta)}x^{-p}L(x),\quad x\to\infty,
	\end{align}
	if and only if  $\bP(\Xi_\theta >x)\sim x^{-p}L(x)$ as $x\to\infty$.
\end{theorem}
\begin{remark}
		If $p\theta\in\Theta_0$ and  $\kappa(p\theta)< p\kappa(\theta)$, then by \eqref{y-Ly}, there exists $r>p$ such that $\bE(Y_\theta^r)<\infty$. As a consequence of Lemma \ref{lemma:A1} in the Appendix, we get that if
	$$\bP(\langle e_\theta,\mathcal{P}\rangle >x)\sim x^{-p}L(x), \quad x\to\infty,$$
	then
	$$\bP(\Xi_\theta >x)\sim \bE(Y_\theta^p)x^{-p}L(x) ,\quad x\to\infty.$$

	Suppose  $\cP=\sum_{i=1}^{N}\delta_{S_i}$, where $S_i,i\ge1$, are iid with the same law as $S$, and are also independent of $N$.  Then
	$$\kappa(\theta)=\beta\Big(m\bE(e^{\theta S})-1\Big)+\varphi(\theta).$$
	\begin{itemize}
		\item[(1)] For $r>1$, $\bE(\langle e_\theta,\cP\rangle^r)<\infty$ if and only if $\bE(N^r)<\infty,$ and $\bE(e^{r\theta S})<\infty.$
		\item[(2)] Assume that $\lim_{x\to\infty}x^{p}L^{-1}\bP(N>x)=c_1$ and $\lim_{x\to\infty}x^{p}L^{-1}\bP(e^{\theta S}>x)=c_2$, where $c_1+c_2>0$.
		Then by Lemma \ref{lemma:A1} in Appendix,
		$$\P(\langle e_\theta,\mathcal{P}\rangle >x)\sim (c_1(Ee^{\theta S})^p+c_2EN)x^{-p}L(x),\quad x\to\infty.$$
		
	\end{itemize}
	\end{remark}

Recall that $Z$ is a non-negative random variable with distribution $\mu$ defined in \eqref{def:mu}.
For $p\in(1,2)$, let $\delta\in(p-1,2-p)$ be
 a constant small enough such that  \eqref{3.8'} holds.

		The following lemma will play an important role in the proof of  Theorem \ref{Theorem:tailW}. Recall that  $\psi_1(\lambda)=\frac{\phi_1(\lambda)}{\lambda}=\frac{\phi(\lambda)-1+\lambda}{\lambda}$.
Since $\phi(\lambda)\leq 1$, we have $\psi_1\in[0,1]$.
		
		\begin{lemma}\label{lemma4}
			Suppose $\theta\in\Theta$ and \eqref{kappa'} holds. Let $p\in(1,2)$ and $L\in\mathbb{L}$. If $p\theta\in\Theta_0$, $\kappa(p\theta)< p\kappa(\theta)$  and $\bP(\Xi_\theta >x)\sim x^{-p}L(x)$ as $x\to\infty$,
			then
			\begin{align}\label{limpsi1}
				\lim_{\lambda\to0}\frac{\psi_1(\lambda)-E(\psi_1(\lambda Z)) }{\lambda^{p-1}L(\lambda^{-1})}
				=\frac{\Gamma(2-p)}{p-1}.
			\end{align}
		\end{lemma}
		{\bf Proof:} Note that $W_\infty(\theta)$ is non-degenerate under the conditions in this lemma. By the definition of $Z$, we have
		\begin{align*}
			E(\psi_1(\lambda Z))
			=\lambda^{-1}\bE \Big(\sum_{i=1}^N\phi_1(\lambda Y_\theta e^{\theta S_i})\Big).
		\end{align*}
Therefore, we only need to prove that
$$\lim_{\lambda\to 0}\frac{\phi_1(\lambda)-\bE \Big(\sum_{i=1}^N\phi_1(\lambda Y_\tau e^{\theta S_i})\Big)}{\lambda^{p}L(\lambda^{-1})}=\frac{\Gamma(2-p)}{p-1}.$$
		By \eqref{phi2} and \eqref{E-I} with $n=1$, we have that
		\begin{align}\label{decomp:phi1}
			\phi(\lambda)&=\bE(\mathcal{T}_1(I(\lambda)))+1-\bE(I(\lambda))\nonumber\\
			&=1-\lambda+\bE(\mathcal{T}_1(I(\lambda)))+\bE\Big(\sum_{i=1}^N\phi_1(\lambda Y_\theta e^{\theta S_i})\Big)-\bE\Big(\sum_{i=1}^N g_1(\lambda Y_\theta e^{\theta S_i})\Big),
		\end{align}
		where $|g_1(\lambda)|\le c(\lambda\wedge \lambda^2)$,
see \eqref{4.9}.
		Thus
		\begin{align}\label{phi4}
			\phi_1(\lambda)-\bE \Big(\sum_{i=1}^N\phi_1(\lambda Y_\theta e^{\theta S_i})\Big) =\bE(\mathcal{T}_1(I(\lambda)))-\bE\Big(\sum_{i=1}^N g_1(\lambda Y_\theta e^{\theta S_i})\Big).
		\end{align}
		Using the fact $|g_1(\lambda)|\le c(\lambda\wedge \lambda^2)\le c\lambda^{p+\delta}$, we have
		\begin{align}\label{est:g1}
			\left|\bE\Big(\sum_{i=1}^N g_1(\lambda Y_\theta e^{\theta S_i})\Big)\right|\le c\lambda^{p+\delta} \bE Y_\theta^{p+\delta}\chi((p+\delta)\theta)\le c\lambda^{p+\delta},
		\end{align}
		where in the last inequality we used
\eqref{3.1} and \eqref{3.1'}.
Then we  only need to prove that
		\begin{align}\label{toprove-limit}
			\lim_{\lambda\to0}\frac{\bE(\mathcal{T}_1(I(\lambda)))}{\lambda^{p}L(\lambda^{-1})}=
			\frac{\Gamma(2-p)}{p-1}.
		\end{align}

		By Lemma \ref{lem:taub} with $n=1$, we have
		\begin{align}\label{J11}
			\lim_{\lambda\to0}\frac{\bE(\mathcal{T}_1(\lambda \Xi_\theta))}{ \lambda^{p}L(\lambda^{-1})}=
			\frac{\Gamma(2-p)}{p-1}.
		\end{align}
		Since $I(\lambda)\le \lambda \Xi_\theta$ and $\mathcal{T}_1$ is increasing, we have
		\begin{align}\label{limitsup}
			\limsup_{\lambda\to0}\frac{\bE(\mathcal{T}_1(I(\lambda)))}{\lambda^{p}L(\lambda^{-1}})\le \lim_{\lambda\to0}\frac{\bE(\mathcal{T}_1(\lambda \Xi_\theta))}{ \lambda^{p}L(\lambda^{-1})}=
			\frac{\Gamma(2-p)}{p-1}.
		\end{align}
		On the other hand, since $\lim_{\lambda\to0}\frac{-\log \phi(\lambda)}{\lambda}=\lim_{\lambda\to0}\frac{1-\phi(\lambda)}{\lambda}=1$,
for any $\epsilon>0$, there exists $\lambda_0>0$ such that for any $\lambda\in(0,\lambda_0)$,
		\begin{align}\label{1.5}
			-\log \phi(\lambda)\ge (1-\epsilon)\lambda.
		\end{align}
		It follows that for any $\lambda>0$,
		\begin{align*}
			I(\lambda )&=\sum_{i=1}^N -\log \phi(\lambda Y_\theta e^{\theta S_i})\ge (1-\epsilon)\lambda Y_\theta\sum_{i=1}^Ne^{\theta S_i}{\bf 1}_{\lambda Y_\theta e^{\theta S_i}<\lambda_0}\\
			&= (1-\epsilon)\lambda Y_\theta\langle e_\theta,\cP\rangle-(1-\epsilon)\lambda Y_\theta\sum_{i=1}^Ne^{\theta S_i}{\bf 1}_{\lambda Y_\theta e^{\theta S_i}\ge \lambda_0}.
		\end{align*}
		Thus, we have
		\begin{align}\label{lowerbound}
			\bE(\mathcal{T}_1(I(\lambda)))\ge &\bE \Big(\mathcal{T}_1((1-\epsilon)\lambda \Xi_\theta)\Big)
			-(1-\epsilon)\lambda \bE\Big(Y_\theta\sum_{i=1}^Ne^{\theta S_i}{\bf 1}_{\lambda Y_\theta e^{\theta S_i}\ge \lambda_0}\Big),
		\end{align}
		where we used the inequality that
		for any $a>b>0$,
		$$\mathcal{T}_1(a-b)\ge \mathcal{T}_1(a)-b.$$
		By  Lemma \ref{lem:taub} again, we have
		\begin{align}\label{J13}
			\lim_{\epsilon\to0}\lim_{\lambda\to0}\frac{\bE \Big(\mathcal{T}_1((1-\epsilon)\lambda \Xi_\theta)\Big)}{ \lambda^{p}L(\lambda^{-1})}=\lim_{\epsilon\to0}(1-\epsilon)^p \frac{\Gamma(2-p)}{p-1}=\frac{\Gamma(2-p)}{p-1}.
		\end{align}
		Note that
		\begin{align*}
			&\bE\Big( \lambda Y_\theta\sum_{i=1}^Ne^{\theta S_i}{\bf 1}_{\lambda Y_\theta e^{\theta S_i}\ge \lambda_0}\Big)
			\le \lambda_0^{-(p+\delta-1)}\bE \big(\lambda^{p+\delta} Y_\theta^{p+\delta}\sum_{i=1}^Ne^{(p+\delta)\theta S_i}\Big)\\
			&=\lambda_0^{-(p+\delta-1)}\lambda^{p+\delta} \bE(Y_\theta^{p+\delta})\chi((p+\delta)\theta)\le \lambda_0^{-(p+\delta-1)}\lambda^{p+\delta}.
		\end{align*}
		Plugging this into \eqref{lowerbound}, we get that
		\begin{align}\label{limitinf}
			\liminf_{\lambda\to0}\frac{\bE(\mathcal{T}_1(I(\lambda)))}{\lambda^{p}L(\lambda^{-1})}\ge
			\frac{\Gamma(2-p)}{p-1}.
		\end{align}
Combining \eqref{limitsup} and \eqref{limitinf}, we
get \eqref{toprove-limit}. The proof is now complete.
		\hfill$\Box$
		
		\bigskip

		Now we are ready to prove Theorem \ref{Theorem:tailW}.
		
		\bigskip
		
		\noindent{\bf Proof of Theorem \ref{Theorem:tailW}:} We first prove sufficiency.
		By Lemma \ref{lem:taub}, \eqref{tail-W} is equivalent to
		\begin{align}\label{limphi1}
			&E(\mathcal{T}_1(\lambda W_\infty(\theta)))
=\phi(\lambda)-1+\lambda=\phi_1(\lambda)\nonumber\\
 &\sim \frac{\Gamma(2-p)}{p-1}\frac{\beta+p\kappa(\theta)-\varphi(p\theta)}{p\kappa(\theta)-\kappa(p\theta)}\lambda^{p}L(\lambda^{-1}),\quad \lambda\to0.
		\end{align}
		We now prove \eqref{limphi1}.

		It follows from \eqref{phi6} with $n=1$ that
		\begin{align}\label{phi7}
			\psi_1(\lambda)=\sum_{k=1}^\infty\delta_{1,k}(\lambda),
		\end{align}
		where where $\delta_{1, k}$ are defined in \eqref{def:deltan1} and \eqref{def:deltank}.
		Since  $0\le \psi_1(\lambda)\le 1$, we have
$$\frac{|\delta_{1,1}(\lambda)|}{\lambda^{p-1}L(1/\lambda)}\le \frac{2}{\lambda^{p-1}L(1/\lambda)}\to0,\quad \lambda\to \infty.$$
By Lemma \ref{lemma4} , we have that  $c:=\sup_{\lambda>0}\frac{|\delta_{1,1}(\lambda)|}{\lambda^{p-1}L(\lambda^{-1})}\in(0,\infty)$ and
		\begin{align*}
			\delta_{1,1}(\lambda)\sim \frac{\Gamma(2-p)}{p-1} \lambda^{p-1}L(\lambda^{-1}),\quad \lambda\to0.
		\end{align*}
		Using  Lemma \ref{lemma2} and  induction on $k$, we get that
		\begin{align}\label{deltan}
			\delta_{1,k}(\lambda)\sim \frac{\Gamma(2-p)}{p-1}\left(\frac{\beta \, \chi(p\theta)}{\beta+p\kappa(\theta)-\varphi(p\theta)}\right)^{k-1} \lambda^{p-1}L(\lambda^{-1}),\quad \lambda\to0
		\end{align}
		and
		\begin{align}\label{deltan2}
			|\delta_{1,k}(\lambda)|\le c A^{k-1} \lambda^{p-1}L(\lambda^{-1}).
		\end{align}
		Combining \eqref{phi7}-\eqref{deltan2} and  the dominated convergence theorem, we have
		\begin{align*}
			\lim_{\lambda\to0}\frac{\psi_1(\lambda)}{\lambda^{p-1}L(\lambda^{-1})}
			=\lim_{\lambda\to0}\frac{\sum_{k=1}^\infty\delta_{1,k}(\lambda)}{\lambda^{p-1}L(\lambda^{-1})}
			=\frac{\Gamma(2-p)}{p-1}\frac{p\kappa(\theta)-\kappa(p\theta)+\beta\chi(p\theta)}{p\kappa(\theta)-\kappa(p\theta)}.
		\end{align*}
Thus \eqref{limphi1} is valid. The proof of sufficiency is now complete.
		
		Now we prove necessity.
		Without loss of generality, we assume that $L\in\mathbb{L}_{p-1}$.  Suppose that
		$$
		\P(W_\infty(\theta>x)\sim cx^{-p}L(x),\quad x\to\infty,
		$$
		where $c=\frac{p\kappa(\theta)-\kappa(p\theta)+\beta\chi(p\theta)}{p\kappa(\theta)-\kappa(p\theta)}.$
		By Lemma \ref{lem:taub} ,  we have
		$$\psi_1(\lambda)\sim c \frac{\Gamma(2-p)}{p-1}\lambda^{p-1}L(\lambda^{-1}),\quad \lambda\to0.$$
		Thus there exists $\lambda_0>0$ such that $\frac{\psi_1(\lambda)}{\lambda^{p-1}L(\lambda^{-1})}\le 2c\frac{\Gamma(2-p)}{p-1}$ for any  $\lambda\in (0,\lambda_0]$.
		Since $L\in\mathbb{L}_{p-1}$,
$\lambda^{p-1}L(\lambda^{-1})$ is increasing. Thus for any $\lambda>\lambda_0$, we have that   $\frac{\psi_1(\lambda)}{\lambda^{p-1}L(\lambda^{-1})}\le \frac{1}{\lambda_0^{p-1}L(\lambda_0^{-1})}$ . Thus $\frac{\psi_1(\lambda)}{\lambda^{p-1}L(\lambda^{-1})}$ is bounded on $[0,\infty).$
		Then by Lemma  \ref{lemma2}, we have that
		$$E(\psi_1(\lambda Z))\sim \frac{\Gamma(2-p)}{p-1}
		\frac{\beta \,\chi(p\theta)}{\kappa(\theta)p-\kappa(p\theta)}\lambda^{p-1}L(\lambda^{-1}),\quad \lambda\to0.$$
		Thus, using   \eqref{phi4} and \eqref{est:g1}, we get that
		\begin{align*}
			\lim_{\lambda\to0}\frac{\bE(\mathcal{T}_1(I(\lambda)))}{\lambda^{p}L(\lambda^{-1})}=\lim_{\lambda\to0}\frac{\lambda(\psi_1(\lambda)-E(\psi_1(\lambda Z)))}{\lambda^{p}L(\lambda^{-1})}=
			\frac{\Gamma(2-p)}{p-1}.
		\end{align*}
		Now using arguments similar to that in the proof of Lemma \ref{lemma4}, we can get that
		$$\lim_{\lambda\to0}\frac{\bE \Big(\mathcal{T}_1(\lambda \Xi_\theta)\Big)}{\lambda^{p}L(\lambda^{-1})}=\frac{\Gamma(2-p)}{p-1}
		,$$
		which implies that
		$$\bP(\Xi_\theta>x)\sim x^{-p}L(x),\quad x\to\infty.$$
		\hfill$\Box$

\section{Applications}\label{sec:application}

		\subsection{
Derivative martingale}\label{subsection: derivative martingale}
		
	In this subsection, we study the derivative martingale in the borderline case $\theta=\theta_*$.
 	
		The following lemma is a key tool for deriving certain  limit results for branching L\'evy processes $\{X_t, t\ge0\}$ by means of their discrete skeletons $\{X_{nt},n=0,1,2\cdots\}$.

		\begin{lemma}\label{Croft-Kingman}
			Let $\{\mathcal{Y}_t, t\ge0\}$ be a real-valued cadlag process
			and let $\mathcal{Y}_\infty$ be a random variable.
			If, for any fixed $s>0$,  $\lim_{n\to\infty}\mathcal{Y}_{ns}=\mathcal{Y}_\infty$  in distribution or in probability, then $\lim_{t\to\infty}\mathcal{Y}_{t}=\mathcal{Y}_\infty$ in  distribution or in probability, respectively.
		\end{lemma}
		{\bf Proof}
Assume that $\lim_{n\to\infty}\mathcal{Y}_{ns}=\mathcal{Y}_\infty$  in distribution.
Define $g(t,\lambda):=E(e^{i\lambda \mathcal{Y}_t})$. Since $\{\mathcal{Y}_t, t\ge0\}$ is cadlag, $g(t,\lambda)$ is right continuous in $t$. Since $\lim_{n\to\infty}\mathcal{Y}_{ns}=\mathcal{Y}_\infty$ in distribution, we have $\lim_{n\to\infty}g(ns,\lambda)=E(e^{i\lambda \mathcal{Y}_\infty})$. According to the Croft-Kingman lemma (see \cite[Corollary 2]{King63}, $\lim_{t\to\infty}g(t,\lambda)=E(e^{i\lambda \mathcal{Y}_\infty})$.
		This implies  that $\lim_{t\to\infty}\mathcal{Y}_t=\mathcal{Y}_\infty$ in distribution.

		If $\lim_{n\to\infty}\mathcal{Y}_{ns}=\mathcal{Y}_\infty$  in probability,
		then $\lim_{n\to\infty}(\mathcal{Y}_{ns}-\mathcal{Y}_\infty)\to 0$  in probability. By the paragraph above, $\lim_{t\to\infty}(\mathcal{Y}_t-\mathcal{Y}_\infty)\to0$ in distribution. Since convergence in distribution to a constant implies convergence in probability, we have  that $\lim_{t\to\infty}\mathcal{Y}_{t}=\mathcal{Y}_\infty$ in probability.
		
		\hfill$\Box$

		Recall the definition of $\kappa'$ and $\kappa''$ in \eqref{def-kappa'} and  \eqref{kappa''}.
		If $\theta\in\Theta_0$, then by the dominated convergence theorem, we have
		$$
		\E\left(\sum_{u\in\mathcal{L}_t}z_u(t)^ke^{\theta z_u(t)}\right)=\frac{d^k}{d\theta^k}\E\left(\sum_{u\in\mathcal{L}_t}e^{\theta z_u(t)}\right)=\frac{d^k}{d\theta^k} e^{\kappa(\theta)t},
\quad k\geq 1.
$$
		
		The following lemma says that it remains valid under certain weaker conditions.

		\begin{lemma}\label{f-moment-14}
Let $\theta\in \Theta$.
			\begin{itemize}
				\item [(1)] If
				\begin{align}\label{4.51}
					\bE\Big(\sum_{i=1}^N |S_i|e^{\theta S_i}\Big)+\int_{|y|>1} |y| e^{\theta y} n(dy)<\infty,
				\end{align}
				then
				$$\E\left(\sum_{u\in\mathcal{L}_t}z_u(t)e^{\theta z_u(t)}\right)=\kappa'(\theta )te^{\kappa(\theta )t}.$$
				\item[(2)] If
				$$
				\bE\Big(\sum_{i=1}^N |S_i|^2e^{\theta  S_i}\Big)+\int_{\R} |y|^2 e^{\theta  y} n(dy)<\infty,
				$$
				then
				$$\E\left(\sum_{u\in\mathcal{L}_t}z_u(t)^2e^{\theta z_u(t)}\right)=(\kappa''(\theta )t+\kappa'(\theta )^2t^2)e^{\kappa(\theta )t}.
				$$
			\end{itemize}
			
		\end{lemma}
		{\bf Proof:}
(1) Assume that \eqref{4.51} holds.
		Since $\int_{|y|>1}|y|e^{\theta y}n(dy)<\infty$,
by the dominated convergence theorem, we have
		$$\varphi'(\theta-)=\lim_{\epsilon\downarrow0}\frac{\varphi(\theta)-\varphi(\theta-\epsilon)}{\epsilon}=a+b^2\theta+\int_\R y(e^{\theta  y}-{\bf 1}_{\{|y|\le 1\}})n(dy).$$
		By
		\cite[Proposition 25.4]{Sato},
		we have $\bE(|\xi_t|e^{\theta \xi_t})<\infty$ for all $t>0$.
		Note that
		\begin{align}\label{4.56}
			\bE(\xi_t e^{\theta\xi_t})=\bE\left(\lim_{\epsilon\downarrow0}\frac{1-e^{-\epsilon \xi_t}}{\epsilon}e^{\theta\xi_t}\right).
		\end{align}
		For $0<\epsilon<\theta/4$ and $x<0$,  we have
		$$
		\frac{e^{-\epsilon x}-1}{\epsilon}\le (-x)e^{-\epsilon x}\le (-x)e^{-\theta x/4}\le \frac{4}{\theta}e^{-\theta x/2},$$
		where in the last inequality we use the inequality $|x|\le \frac{4}{\theta}e^{\theta |x|/4}$.
		Thus for $0<\epsilon<\theta/4$, we have
		\begin{align}\label{4.57}
			\frac{|1-e^{-\epsilon \xi_t}|}{\epsilon}e^{\theta\xi_t}\le 1_{\xi_t>0}\xi_t e^{\theta\xi_t}+1_{\xi_t\le 0}\frac{4}{\theta}e^{\frac{\theta}{2} \xi_t}\le |\xi_t| e^{\theta\xi_t}+
			\frac{4}{\theta}.
		\end{align}
		By \eqref{4.56} and \eqref{4.57}, applying  the dominated convergence theorem, we get
		\begin{align}\label{4.54}
			\bE(\xi_t e^{\theta\xi_t})=\lim_{\epsilon\downarrow0}\bE\left(\frac{1-e^{-\epsilon \xi_t}}{\epsilon}e^{\theta\xi_t}\right)=\lim_{\epsilon\downarrow0}\frac{e^{\varphi(\theta)t}-e^{\varphi(\theta-\epsilon)t}}{\epsilon}=\varphi'(\theta-)te^{\varphi(\theta)t}.
		\end{align}
		By Lemma \ref{f-moment} we have that  $\E(\sum_{u\in\mathcal{L}_t}|z_u(t)|e^{\theta z_u(t)})<\infty.$
		Using the Markov property and the branching property, we get that for any $t>0$,
		\begin{align}\label{4.52}
			U(t)&:=\E(\sum_{u\in\mathcal{L}_t}z_u(t)e^{\theta z_u(t)})\nonumber\\
			&=e^{-\beta t}\bE(\xi_t e^{\theta \xi_t})+\bE\int_0^t \beta e^{-\beta r}\sum_{k=1}^N \E_{\delta_{\xi_r+S_k}}\!\!\left(\sum_{u\in\mathcal{L}_t}z_u(t-r)e^{\theta z_u(t-r)}\right)\!dr.
		\end{align}
		Since
		$$\E_{\delta_x}\left(\sum_{u\in\mathcal{L}_t}z_u(t)e^{\theta z_u(t)}\right)=\E\left(\sum_{u\in\mathcal{L}_t}(x+z_u(t))e^{\theta(x+z_u(t))}\right)=e^{\theta x}(xe^{\kappa(\theta)t}+U(t)),
		$$
		we have
		\begin{align}\label{4.55}
			&\bE  \left[\sum_{k=1}^N\E_{\delta_{\xi_r+S_k}}\left(\sum_{u\in\mathcal{L}_t}z_u(t-r)e^{\theta z_u(t-r)}\right) \right]\nonumber\\
			=&\bE\left[\sum_{k=1}^N(\xi_r+S_k)e^{\theta(\xi_r+S_k)}\right]e^{\kappa(\theta)(t-r)}+\bE\left[\sum_{k=1}^Ne^{s(\xi_r+S_k)}\right]U(t-r)\nonumber\\
			=&\left[\bE(\xi_r
e^{\theta\xi_r})
\chi(\theta)+A(\theta)e^{\varphi(s)r}\right]e^{\kappa(\theta)(t-r)}+e^{\varphi(\theta)r}\chi(\theta)U(t-r)\nonumber\\
			=&\left[\varphi'(\theta-)\chi(\theta)r+A(\theta)\right]e^{\varphi(\theta)r}e^{\kappa(\theta)(t-r)}+e^{\varphi(\theta)r}\chi(\theta)U(t-r),
		\end{align}
		where $A(\theta)=\bE\Big(\sum_{i=1}^N S_ie^{\theta S_i}\Big)$.
		By \eqref{4.52}  and \eqref{4.55} we have
		\begin{align}\label{4.53}
			U(t)=&e^{(\varphi(\theta)-\beta) t}\varphi'(\theta-) t+e^{\kappa(\theta)t}\int_0^t \beta\chi(\theta)\varphi'(\theta-) r e^{-\beta\chi(\theta)r}\,dr\nonumber\\
			&+e^{\kappa(\theta)t}\int_0^t \beta A(\theta)e^{-\beta\chi(\theta)r}\,dr
			+\int_0^t \beta\chi(\theta)e^{(\varphi(\theta)-\beta) r}U(t-r)\,dr.
		\end{align}
		Since
		\begin{align*}
			&e^{\kappa(\theta)t}\int_0^t \beta\chi(\theta)\varphi'(\theta-) r e^{-\beta\chi(\theta)r}\,dr=e^{\kappa(\theta)t}\varphi'(\theta-) \int_0^t r d(-e^{-\beta\chi(\theta)r})\\
			=&-e^{(\varphi(\theta)-\beta) t}\varphi'(\theta-) t+e^{\kappa(\theta)t}\varphi'(\theta-) \int_0^t e^{-\beta\chi(\theta)r} dr,
		\end{align*}
		we have by \eqref{4.53} that
		\begin{align*}
			U(t)=\kappa'(\theta)e^{\kappa(\theta)t}\ \int_0^t e^{-\beta\chi(\theta)r} dr+e^{(\varphi(\theta)-\beta) t}\int_0^t \beta\chi(\theta)e^{-(\varphi(\theta)-\beta) r}U(r)\,dr,
		\end{align*}
		which implies that
		\begin{align*}
			V(t):=e^{-\kappa(\theta)t}U(t)=\kappa'(\theta)\ \int_0^t e^{-\beta\chi(\theta)r} dr+e^{-\beta\chi(\theta) t}\int_0^t \beta\chi(\theta)e^{\beta\chi(\theta) r}V(r)\,dr.
		\end{align*}
		Differentiating with respect to
		$t$, we get
		\begin{align*}
			V'(t)&=\kappa'(\theta)e^{-\beta\chi(\theta)t} -\beta\chi(\theta)e^{-\beta\chi(\theta)t}\int_0^t \beta\chi(\theta)e^{\beta\chi(\theta) r}V(r)\,dr+\beta\chi(\theta)V(t)\nonumber\\
			&=\kappa'(\theta)e^{-\beta\chi(\theta)t} +\beta\chi(\theta)\kappa'(\theta)\ \int_0^t e^{-\beta\chi(\theta)r} dr=\kappa'(\theta).
		\end{align*}
		Since $V(0)=0,$ then $V(t)=\kappa'(\theta)t$, this implies that $U(t)=\kappa'(\theta)te^{\kappa(\theta)t}.$
		
The proof of (2) is similar.
Here we omit the details.
		\hfill$\Box$
		
		\bigskip

		Suppose there exists $\theta_*\in\Theta'$ such that $\theta_*\kappa'(\theta_*)=\kappa(\theta_*)$.
		Notice  that   $W_\infty(\theta_*)=0$ almost surely.
In this boundary case
$\theta=\theta^*$,
there is a signed martingale $\{D_t\}$, called the derivative martingale, which is defined as
		$$
		D_t:=\sum_{u\in \mathcal{L}_t} (\kappa(\theta_*)t-\theta_*z_u(t))e^{\theta_* z_u(t)-\kappa(\theta_*)t},\quad t\geq 0.
		$$
		For branching random walks, A\"idekon \cite{Aid13} derived  sufficient integrability conditions for the non-degeneracy of the limit of the derivative martingale.
Later, Chen  \cite{Ch15}
showed that these conditions are  also necessary.
Recently, Mallein and Shi \cite{MSh23}
proved the corresponding results for the general branching L\'evy process model of \cite{BM18}.

		\begin{lemma}[Mallein and Shi \cite{MSh23}]\label{derivative-mart}
			Suppose there exists $\theta_*\in \Theta'$ such that
			$\theta_*\kappa'(\theta_*)=\kappa(\theta_*)$.
			The derivative martingale $D_t$ converges to a non-negative non-degenerate limit $D_\infty$ if and only if
			\begin{align}\label{condition-deriv}
				\bE\left( \langle e_{\theta_*}, \mathcal{P}\rangle \log_+^2 \langle e_{\theta_*} , \mathcal{P}\rangle \right)<\infty \mbox{ and } \bE\left( V\log_+(V)\right)<\infty,
			\end{align}
			where $V=\sum_{i=1}^N {\bf 1}_{\{S_i\le 0\}}(-S_i)e^{\theta_* S_i}.$
		\end{lemma}

		The following result says that, after appropriate scaling, $W_t(\theta_*)$ converges to $D_\infty$ in probability. This result was proved for branching random walks in \cite{ASh14}. Applying the result in \cite{ASh14} for $W_{n\delta}$ and  using Lemma \ref{Croft-Kingman}, we can get that the result holds for $W_t$.
Note that,  in the following result, our conditions are on $\mathcal{P}$ instead of on $W_\delta$.

		\begin{proposition}\label{lemma-rate}
			Suppose that there exists $\theta_*>0$ such that
		$\kappa''(\theta_*)<\infty$, $\theta_*\kappa'(\theta_*)=\kappa(\theta_*)$ and
			\begin{align}\label{H3}
				\bE\left( \langle e_{\theta_*}, \mathcal{P}\rangle \log_+^2 \langle e_{\theta_*} , \mathcal{P}\rangle \right)<\infty \mbox{ and } \bE\left( \langle e_{\theta_*-\epsilon}, \mathcal{P}\rangle \log_+ \langle e_{\theta_*-\epsilon} , \mathcal{P}\rangle \right)<\infty,
			\end{align}
			for some $\epsilon\in(0, \theta_*)$.
			Then we have
			$$
			\sqrt{t} W_t(\theta_*)\overset{P}{\rightarrow} \sqrt{\frac{2}{\pi \theta_*\kappa''(\theta_*)}}D_\infty,\quad \mbox{ as } t\to\infty.$$
		\end{proposition}

		\noindent{\bf Proof:} Since $\kappa''(\theta_*)<\infty$,   applying Lemma \ref{f-moment} with $r=2$, we get that for all $t>0$,
		$$
		\E\Big(\sum_{u\in\mathcal{L}_t} (z_u(t))^2e^{\theta_*z_u(t)}\Big)<\infty.
		$$
		By \eqref{H3} and
Theorem \ref{proposition-log},
we have
		$$
		\bE\left( \langle e_{\theta_*}, X_t\rangle \log_+^2 \langle e_{\theta_*} , X_t\rangle \right)<\infty\mbox{ and } \bE\left( \langle e_{\theta_*-\epsilon}, X_t\rangle \log_+ \langle e_{\theta_*-\epsilon} , X_t\rangle \right)<\infty.
		$$
		Since
		$$
		V_t:=\sum_{u\in \mathcal{L}_t} {\bf 1}_{\{z_u(t)\le 0\}}|z_u(t)|e^{\theta_* z_u(t)}\le\epsilon^{-1} \sum_{u\in \mathcal{L}_t}e^{(\theta_*-\epsilon)z_u(t)},
		$$
		we have $\E(V_t \log_+(V_t))<\infty.$
 Now applying \cite[Theorem 1.1]{ASh14}
to  $\{W_{nt},n=0,1,\cdots\}$, we get that  as $n\to\infty,$
		\begin{align}\label{discrete}
			\sqrt{n } W_{n t}(\theta_*)\overset{P}{\rightarrow} \sqrt{\frac{2}{\pi \sigma_t^2}}D_\infty,
		\end{align}
		where
		$$
		\sigma_t^2=\E\left(\sum_{u\in \mathcal{L}_t} (\kappa(\theta_*)t-\theta_*z_u(t))^2e^{\theta_* z_u(t)-\kappa(\theta_*)t}\right).
$$
	
		Using Lemma \ref{f-moment-14} and the fact that $\theta_*\kappa'(\theta_*)=\kappa(\theta_*)$, we get
		\begin{align*}
		\sigma_t^2=\kappa(\theta_*)^2t^2-2\kappa(\theta_*)\theta_*\kappa'(\theta_*)t^2+\theta_*^2(\kappa''(\theta_*)t+\kappa'(\theta_*)^2t^2)=\theta_*^2\kappa''(\theta_*)t.
		\end{align*}
		Thus
		\begin{align*}
			\sqrt{n t} W_{n t}(\theta_*)\overset{P}{\rightarrow}\sqrt{\frac{2}{\pi \theta_*^2\kappa''(\theta_*)}}D_\infty,\quad \mbox{ as } t\to\infty.
		\end{align*}
		Now the desired result follows from Lemma \ref{Croft-Kingman} immediately.
		\hfill$\Box$

\begin{remark}
Dadoun \cite{Dadoun} used the same technique  to study the limit behavior of growth-fragmentation processes, i.e.,  by  transferring the result  for branching random walks in \cite{ASh14} to discrete skeletons of the growth-fragmentation, and then infer the behavior of the continuous time process with the help of Lemma \ref{Croft-Kingman}.
\end{remark}

		\subsection{
Extremal process of branching L\'evy process}\label{ss:4.3}
For any $t>0$, define $M_t:=\max\{u\in\mathcal{L}_t:z_u(t)\}$.
		We first give a result about $M_t$.
		\begin{lemma}\label{lemma:rightmost}
			If $\theta_0\in \Theta$, then
			\begin{equation*}
				\lim_{t\to\infty}e^{\theta_0M_t-\kappa(\theta_0)t}=0,\quad \P\mbox{-a.s.}
			\end{equation*}
		\end{lemma}
		{\bf Proof:}
		Since
		\begin{align*}
			e^{\theta_0 M_t-\kappa(\theta_0)t}\le e^{-\kappa(\theta_0)t}\sum_{u\in\cL_t}e^{\theta_0 z_u(t)}=W_t(\theta_0),
		\end{align*}
		we have
		$$
		U:=\limsup_{t\to\infty}e^{\theta_0 M_t-\kappa(\theta_0)t}\le W_\infty(\theta_0)\in[0,\infty).
		$$
		Thus, $\E(U)\le \E(W_\infty(\theta_0))\le 1$.
		To prove the assertion of the lemma,
		it suffices to show $\P(U>0)=0$.
		For any $0\le s<t$,
 by the branching property, we have
		\begin{align*}
			M_t=\sup_{u\in\cL_s}\{z_u(s)+M_{t-s}^u\},
		\end{align*}
		where $M_{t-s}^u$, $u\in\cL_s$, are independent copies of $M_{t-s}$, and  independent of $\{z_u(s),u\in\cL_s\}$.
		Thus
		\begin{align*}
			U=\sup_{u\in\cL_s} e^{\theta_0z_u(s)}e^{-\kappa(\theta_0)s}U^u,
		\end{align*}
		where $U^u$, $u\in\cL_s$ are independent copies of $U$, and  independent of $\{z_u(s),u\in\cL_s\}$.
		Since $\E(\|X_s\|)=e^{\beta(\E(N)-1)s}>1$, $\P(\|X_s\|>1)>0$.
		If $\P(U>0)>0$, then
		\begin{align*}
			\E(U)&=\E\Big(\sup_{u\in\cL_s} e^{\theta_0z_u(s)}e^{-\kappa(\theta_0)s}U^u\Big)
			<\E\Big(\sum_{u\in\cL_s} e^{\theta_0z_u(s)}e^{-\kappa(\theta_0)s}U^u\Big)\\
			&=\E\Big(\sum_{u\in\cL_s} e^{\theta_0z_u(s)}e^{-\kappa(\theta_0)s}\Big)\E(U)=\E(U),
		\end{align*}
		which is   a contradiction.
		Thus $\P(U>0)=0$. The proof is complete.
		
		\hfill$\Box$
		
		In this subsection, we assume that there exists $\theta_*>0$ such that  $\theta_*\kappa'(\theta_*)=\kappa(\theta_*)$.
		Note that for any $t>0,$ $\{\widetilde{X}_n:=-\theta_*X_{nt}+\kappa(\theta_*)tn,n=0,1,\cdots\}$ is a branching random walk with $ \E(\langle e^{-x},\widetilde{X}_1\rangle)=1$ and $\E(\langle xe^{-x},\widetilde{X}_1\rangle)=0$.
		Combining \cite[Theorem 1.2]{HSh09}, \cite[Theorem 1.1]{Aid13} and \cite[Theorems 1.1 and 2.3]{Madaule} for $\{\widetilde{X}_n\}$,  with our Lemma \ref{Croft-Kingman}, Lemma \ref{f-moment},
Theorem \ref{proposition-log},
 we immediately get the following Propositions \ref{prop:Maximum}, \ref{maximum-weak} and  \ref{prop:extremal}.

		\begin{proposition}\label{prop:Maximum}
				Assume that there exists $\theta_*\in\Theta_0$ such that  $\theta_*\kappa'(\theta_*)=\kappa(\theta_*)$.
			If there exist $\delta>0$ and $\delta_->0$ such that
			\begin{align}\label{maximum-1}
				\bE(N^{1+\delta})<\infty,\mbox{ and } \kappa(-\delta_-)<\infty,
			\end{align}
			then conditioned on the event of survival, we have
			\begin{align}\label{maximum-2}
				\lim_{t\to\infty}\frac{M_t-c_*t}{\log t}=-\frac{3}{2\theta_*},  \quad \mbox{  in probability.}
			\end{align}
			where $c_*=\kappa(\theta_*)/\theta_*$.
		\end{proposition}

\begin{remark}
Dadoun \cite[Corollary 2.10] {Dadoun}
got a similar result for  the largest fragment of growth-fragmentation  using similar arguments.
\end{remark}

		\bigskip
		
		\begin{proposition}\label{maximum-weak}
					Suppose there exists $\theta_*>0$ such that
			$\kappa''(\theta_*)<\infty$ and   $\theta_*\kappa'(\theta_*)=\kappa(\theta_*)$.
						If $N$ is  finite, $\mathcal{P}$ is non-lattice and \eqref{H3} holds, then there exists a constant $C^*\in(0,\infty)$ such that
			$$\lim_{t\to\infty}\P(M_t-m(t)\le x)=\E(e^{-C^* e^{-\theta_*x}D_\infty}),$$
			where $m(t)=c_*t-\frac{3}{2\theta_*}\log t.$
		\end{proposition}
		
		Now  we introduce the point processes:
		$$\mathcal{E}_t^*:=\sum_{u\in\mathcal{L}_t}\delta_{z_u(t)-m(t)-\frac{\log D_\infty}{\theta_*}},$$
and
		$$\mathcal{E}_t:=\sum_{u\in\mathcal{L}_t}\delta_{z_u(t)-m(t)},$$
		where $m(t)=c_*t-\frac{3}{2\theta_*}\log t.$
		
		\begin{proposition}\label{prop:extremal}
			Suppose there exists $\theta_*>0$ such that
			$\kappa''(\theta_*)<\infty$ and   $\theta_*\kappa'(\theta_*)=\kappa(\theta_*)$.
						If  $\mathcal{P}$ is non-lattice and  \eqref{H3} holds, then, conditioned on the event of non-extinction, $(\mathcal{E}_t^*,D_t)$
			converges jointly in law to ($\mathcal{E}_\infty^*,D_\infty)$ as $t\to\infty$, where $\mathcal{E}_\infty^*$ and $D_\infty$ are independent. Moreover, for $f\in C_c^+(\R)$, $C(f):=-\log E(e^{-\langle f,\mathcal{E}_\infty^*\rangle})$ satisfies
			$$
			C(f(x+\cdot))=e^{\theta_*x}C(f), \quad  \mbox{ for any } x.
			$$
		\end{proposition}

\bigskip

		\begin{corollary}
			Suppose there exists $\theta_*>0$ such that
			$\kappa''(\theta_*)<\infty$ and   $\theta_*\kappa'(\theta_*)=\kappa(\theta_*)$.
			If $\mathcal{P}$ is non-lattice and  \eqref{H3} holds, then as $t\to\infty$,
			$\mathcal{E}_t$ converges in law to a random measure $\mathcal{E}_\infty$ with  Laplace transform
			$$E(e^{-\langle f,\mathcal{E}_\infty\rangle})=\E(e^{-C(f)D_\infty}).$$
		\end{corollary}

\begin{remark}
Comparing with results on the extreme processes of supercritical branching L\'evy processes with local branching mechanism, we see that the non-local nature of the branching mechanism does not
affect the form of the extreme process, it only affects the parameters involved in  the extreme process.
\end{remark}

	\subsection{The central limit theorem of $W_t(\theta)-W_\infty(\theta)$}\label{clt}
		In this section, we apply
		the results in earlier sections
to prove some
central limit theorems for $W_t(\theta)-W_\infty(\theta)$.

If $\mathcal{Y}, \mathcal{Y}_t, t\ge0$, are real-valued random variables, we write
$$\mathbf{L}(\mathcal{Y}_t|\mathcal{F}_t)\overset{w}{\to }\mathbf{L}(\mathcal{Y}|\mathcal{F}_\infty) \mbox{ in probability},$$
if $\E(f(\mathcal{Y}_t)|\mathcal{F}_t)\to \E(f(\mathcal{Y})|\mathcal{F}_\infty)$ in probability for every $f\in C_b(\R)$. By the dominated convergence theorem, we see that  $\mathbf{L}(\mathcal{Y}_t|\mathcal{F}_t)\overset{w}{\to }\mathbf{L}(\mathcal{Y}|\mathcal{F}_\infty) \mbox{ in probability}$ implies that $ \mathcal{Y}_t\overset{d}{\to}\mathcal{Y}$ as $t\to\infty$.

Similarly, we write
$$\mathbf{L}(\mathcal{Y}_t|\mathcal{F}_t)\overset{w}{\to }\mathbf{L}(\mathcal{Y}|\mathcal{F}_\infty) ,\quad \P\mbox{-a.s., }$$
if $\E(f(\mathcal{Y}_t)|\mathcal{F}_t)\to \E(f(\mathcal{Y})|\mathcal{F}_\infty), \P$-a.s.  for every $f\in C_b(\R)$. Note that $\mathbf{L}(\mathcal{Y}_t|\mathcal{F}_t)\overset{w}{\to }\mathbf{L}(\mathcal{Y}|\mathcal{F}_\infty),  \P$-a.s. implies that $ \mathcal{Y}_t\overset{d}{\to}\mathcal{Y}$ as $t\to\infty$.
		
		Recall that
		\begin{align*}
			W_\infty(\theta)
			=e^{-\kappa(\theta)t}\sum_{u\in\cL_t}e^{\theta z_u(t)} W^u_\infty(\theta),
		\end{align*}
		where $\{W^u_\infty(\theta), u\in\cL_t\}$ are i.i.d., each having  the same law as $W_\infty(\theta)$, and are also independent of  $X_t$.  Observe that
		\begin{align}\label{decomp2-a}
			e^{\kappa(\theta)t}(W_t(\theta)-W_\infty(\theta)) &= \sum_{u\in\cL_t}e^{\theta z_u(t)} (1-W^u_\infty(\theta)).
		\end{align}

\subsubsection{Normal central limit theorems}

It follows from Theorem \ref{main-Lp}, under the condition \eqref{H2}, $W_\infty(\theta)$ is square integrable. The next lemma gives an explicit expression for its variance.

		\begin{lemma}
			Let $\theta\in\Theta$.
			If
			\begin{align}\label{H2}
				\kappa(2\theta)< 2\kappa(\theta) \mbox{ and } \bE(\langle e_{\theta},\cP\rangle^2)<\infty,
			\end{align}
			then
			\begin{equation}\label{var:W}
				\sigma_\theta^2:=\mathbb{V}ar(W_\infty(\theta))=\frac{\beta\left(\bE \langle e_{\theta},\cP\rangle^2-\chi(2\theta)\right)}{2\kappa(\theta)-\kappa(2\theta)}-1.
			\end{equation}
		\end{lemma}
		{\bf Proof:}
		By \eqref{decomp:W2}, we have
		\begin{align*}
			\E(W_\infty(\theta)^2) &=\E\left[ e^{-2\kappa(\theta)\tau_1} e^{2\theta z_\varnothing(\tau_1 )}\left(\sum_{j=1}^{N^\varnothing}e^{\theta S_i^\varnothing}W_{\infty}^j(\theta)\right)^2\right]\\
			&=\bE\left[ Y_\theta^2\right]\E\left[ \left(\sum_{j=1}^{N^\varnothing}e^{\theta S_i^\varnothing}W_{\infty}^j(\theta)\right)^2\right].
		\end{align*}
		It follows from \eqref{y-Ly} that
		$$\bE\left[ Y_\theta^2\right]=\frac{\beta}{\beta+2\kappa(\theta)-\varphi(2\theta)}.$$
		Since $\{W_{\infty}^j(\theta),j\ge 1\}$ are   i.i.d. with mean 1, we have
		\begin{align*}
			&\E\left[ \left(\sum_{j=1}^{N^\varnothing}e^{\theta S_i^\varnothing}W_{\infty}^j(\theta)\right)^2\right]\\
			=&\E\left[ \sum_{j=1}^{N^\varnothing}e^{2\theta S_i^\varnothing}\E(W_{\infty}^j(\theta)^2)\right]+2\E\left[ \sum_{1\le i<j\le N^\varnothing }e^{\theta S_i^\varnothing}e^{\theta S_j^\varnothing}\E(W_{\infty}^i(\theta))\E(W_{\infty}^j(\theta))\right]\\
			&=\bE\langle e_{2\theta}, \mathcal{P}\rangle	\E(W_\infty(\theta)^2)+\bE(\langle e_\theta,\mathcal{P}\rangle^2-\langle e_{2\theta}, \mathcal{P}\rangle).
		\end{align*}
		Thus
		\begin{align*}
			\E(W_\infty(\theta)^2)=\frac{\beta\chi(2\theta) }{\beta+2\kappa(\theta)-\varphi(2\theta)}	\E(W_\infty(\theta)^2)+\frac{\beta}{\beta+2\kappa(\theta)-\varphi(2\theta)}\bE(\langle e_\theta,\mathcal{P}\rangle^2-\langle e_{2\theta}, \mathcal{P}\rangle).
		\end{align*}
		Since $\beta+2\kappa(\theta)-\varphi(2\theta)=2\kappa(\theta)-\kappa(2\theta)+\beta\chi(2\theta)$,  the desired result follows immediately.
		
		\hfill$\Box$

		\begin{theorem}\label{them-clt2}
			Let $\theta\in\Theta$. Assume that \eqref{H2} holds.
			\begin{itemize}
				\item [(1)] If 	there exists a positive function $b(t)$ and a random variable $\Lambda_\theta$ such that
				\begin{align}\label{clt-1}
					\lim_{t\to \infty}b(t)^2 W_{t}(2\theta)= \Lambda_\theta,
				\end{align}
				and
				\begin{align}\label{clt-2}
					\lim_{t\to \infty}b(t)e^{-\kappa(2\theta)t/2}e^{\theta M_t}=0
				\end{align}
				in probability, then as $t\to\infty$,
				$$
				\mathbf{L}(b(t)e^{(\kappa(\theta)-\kappa(2\theta)/2)t}	\left(W_t(\theta)-W_\infty(\theta)\right)|
				\mathcal{F}_t)\overset{w}{\to }\mathbf{L}(\sigma_\theta\sqrt{\Lambda_\theta}Z|\mathcal{F}_\infty)$$   in probability with respect to $\P$,
				where $\sigma_\theta^2=\mathbb{V}ar(W_\infty(\theta))$, and  $Z$ is a standard normal random variable independent of  $\mathcal{F}_\infty$.
				\item[(2)] If the limits in \eqref{clt-1} and \eqref{clt-2} hold almost surely, then
				$$\mathbf{L}(b(t)e^{(\kappa(\theta)-\kappa(2\theta)/2)t}	\left(W_t(\theta)-W_\infty(\theta)\right)|
				\mathcal{F}_t)\overset{w}{\to }\mathbf{L}(\sigma_\theta\sqrt{\Lambda_\theta}Z|\mathcal{F}_\infty),  \quad \P\mbox{-a.s.}.$$
			\end{itemize}
		\end{theorem}
{\bf Proof:}
(1) Note that
		\begin{align}\label{decomp2-b}
			b(t)e^{(\kappa(\theta)-\kappa(2\theta)/2)t}	\left(W_t(\theta)-W_\infty(\theta)\right)&= b(t)e^{-\kappa(2\theta)t/2}\sum_{u\in\cL_t}e^{\theta z_u(t)} (1-W^u_\infty(\theta))\nonumber\\
			&=:\sum_{u\in\cL_t}\mathcal{V}(t,z_u(t))(1-W^u_\infty(\theta)),
		\end{align}
		where $$\mathcal{V}(t,x):=b(t)e^{-\kappa(2\theta){t}/2}e^{\theta x} .$$
		Put $\Psi(\lambda):=\E(e^{i\lambda(1-W_\infty(\theta))}).$ Then  for any $\lambda\in\R$,
		\begin{align}\label{5.61}
			\E\left[ \exp\{i\lambda b(t)e^{(\kappa(\theta)-\kappa(2\theta)/2)t}	\left(W_t(\theta)-W_\infty(\theta)\right)\}|\mathcal{F}_t\right]=\prod_{u\in\mathcal{L}_t}\Psi(\lambda \mathcal{V}(t,z_u(t))).
		\end{align}
		By \cite[Section 2.4, Lemma 4.3]{Durrett}, we have
		\begin{align}\label{5.62}
			&\left|\prod_{u\in\mathcal{L}_t}\Psi(\lambda \mathcal{V}(t,z_u(t)))-e^{-\frac{1}{2}\lambda^2\sigma_\theta^2b(t)^2W_{t}(2\theta)}\right|
			\le \sum_{u\in\mathcal{L}_t}\left| \Psi(\lambda \mathcal{V}(t,z_u(t)))
			-e^{-\frac{1}{2}\lambda^2\sigma_\theta^2\mathcal{V}(t,z_u(t))^2}\right|.
		\end{align}
		Recall that $\sigma_\theta^2=\E(1-W_\infty(\theta))^2$. For any $x\in \R$, by Taylor's formula,  we have that
		\begin{align}\label{5.63}
			\left|\Psi(x)-e^{-\frac{1}{2}\sigma_\theta^2x^2}\right|&\le \left|\Psi(x)-1+\frac{1}{2}\sigma_\theta^2x^2\right|+ \left|e^{-\frac{1}{2}\sigma_\theta^2x^2}-1+\frac{1}{2}\sigma_\theta^2x^2\right|\nonumber\\
			&\le \E\left(x^2(1-W_\infty(\theta))^2\wedge\frac{|x|^3|1-W_\infty(\theta)|^3}{6}\right)+\frac{1}{8}\sigma_\theta^4x^4.
		\end{align}
		Note that for any $\epsilon>0$,
		\begin{align}\label{5.64}
			&\E\left(x^2(1-W_\infty(\theta))^2\wedge\frac{|x|^3|1-W_\infty(\theta)|^3}{6}\right)\nonumber\\
			\le &x^2\E\left((1-W_\infty(\theta))^2;|1-W_\infty(\theta)|>\epsilon |x|^{-1} \right)\nonumber\\
			&+|x|^3\E\left(|1-W_\infty(\theta)|^3;|1-W_\infty(\theta)|\le \epsilon |x|^{-1} \right)\nonumber\\
			\le &x^2G(\epsilon |x|^{-1} )+\epsilon \sigma_\theta^2x^2,
		\end{align}
		where
		$$G(x)=\E( (1-W_\infty(\theta))^2;|1-W_\infty(\theta)|>x)\to 0, \quad \mbox{ as} x\to\infty.$$
		Combining \eqref{5.62}-\eqref{5.64}, we get that  for any $\epsilon>0$,
		\begin{align}\label{5.65}
			&\left|\prod_{u\in\mathcal{L}_t}\Psi(\lambda \mathcal{V}(t,z_u(t)))-e^{-\frac{1}{2}\lambda^2\sigma_\theta^2b(t)^2W_{t}(2\theta)}\right|\nonumber\\
			\le &\sum_{u\in\mathcal{L}_t} \lambda^2 \mathcal{V}(t,z_u(t))^2G(\epsilon \lambda^{-1}\mathcal{V}(t,z_u(t))^{-1})+\epsilon \sigma_\theta^2\sum_{u\in\mathcal{L}_t} \lambda^2 \mathcal{V}(t,z_u(t))^2+\frac{1}{8}\sigma_\theta^4\sum_{u\in\mathcal{L}_t} \lambda^4 \mathcal{V}(t,z_u(t))^4\nonumber\\
			\le&b(t)^2W_t(2\theta)\left[\lambda^2G(\epsilon\lambda^{-1}b(t)^{-1}e^{\kappa(2\theta)t/2}e^{-\theta M_t})+\epsilon \sigma_\theta^2\lambda^2+\frac{1}{8}\sigma_\theta^4\lambda^4b(t)^{2}e^{-\kappa(2\theta)t}e^{2\theta M_t}\right],
		\end{align}
		where in the last inequality we used the facts $\mathcal{V}(t,z_u(t))\le b(t)e^{-\kappa(2\theta)t/2}e^{\theta M_t}$ and $\sum_{u\in\mathcal{L}_t}\mathcal{V}(t,z_u(t))^2=b(t)^2W_t(2\theta)$.
		If  \eqref{clt-1} and \eqref{clt-2} hold,  then we have
		\begin{align*}
			\left|\prod_{u\in\mathcal{L}_t}\Psi(\lambda \mathcal{V}(t,z_u(t)))-e^{-\frac{1}{2}\lambda^2\sigma_\theta^2b(t)^2W_{t}(2\theta)}\right|
			\overset{\P}{\to}0.
		\end{align*}
		Now applying  \eqref{clt-1}, we get
		$$\prod_{u\in\mathcal{L}_t}\Psi(\lambda \mathcal{V}(t,z_u(t)))
		\overset{\P}{\to}
		e^{-\frac{1}{2}\lambda^2\sigma_\theta^2\Lambda_\theta}=\E\left(e^{i\lambda \sigma_\theta\sqrt{\Lambda_\theta}Z}|\mathcal{F}_\infty\right).$$
		The desired result then follows from \eqref{5.61}.

		(2) The  second part of the theorem can be proved similarly. We omit the details.
		
		\hfill$\Box$

Assume that
	there exists $\theta_*\in\Theta'$ such that $\theta_*\kappa'(\theta_*)=\kappa(\theta_*)$.
Since $\kappa(s)/s$ is increasing
on $(\theta_*, \theta_+)$, then
$\kappa(2\theta)<2\kappa(\theta)$ implies $\theta<\theta_*$.
The rates of $W_t(\theta)$  converges to $W_\infty(\theta)$ as $t\to\infty$ are different in  the cases $\theta<( =, >) \theta_*/2$.
		In the following two results, we will give central limit theorems for $W_t(\infty)-W_\infty(\theta)$
for the cases $\theta<\theta_*/2$ and $\theta=\theta_*/2$ respectively.
The case $\theta>\theta_*/2$
will be dealt with in Proposition \ref{CLT-extremal case}.
		
		\begin{corollary}\label{them-clt}
			Let $\theta\in\Theta$.
			If \eqref{H2} holds, then
			\begin{align*}
				\mathbf{L}(e^{(\kappa(\theta)-\kappa(2\theta)/2)t}	\left(W_t(\theta)-W_\infty(\theta)\right)|\mathcal{F}_t) \overset{w}{\to}\mathbf{L}(\sigma_\theta\sqrt{W_\infty(2\theta)}Z|\mathcal{F}_\infty),\mbox{ $\P$-a.s.,}
			\end{align*}
			where $\sigma_\theta^2=\mathbb{V}ar(W_\infty(\theta)),$ and $Z$ is a standard normal random variable independent of  $\mathcal{F}_\infty$.
		\end{corollary}
		{\bf Proof:} Since $2\theta\in\Theta$,
		$W_t(2\theta)\to W_\infty(2\theta)$, a.s., and thus by Lemma \ref{lemma:rightmost},
		$$e^{-\kappa(2\theta)t/2}e^{\theta M_t}\to 0, \quad \P\mbox{-a.s.}$$
		Applying  Theorem \ref{them-clt2} with   $b(t)=1$, we get the desired result.
		
		\hfill$\Box$

Note that  when $2\theta\ge \theta_*$,  $2\theta\kappa'(2\theta)\ge\kappa(2\theta)$ which implies that $W_{\infty}(2\theta)=0$.
		Thus, when $2\theta\ge \theta_*$, the limit in Corollary \ref{them-clt} is $0$.
Consequently,
Corollary \ref{them-clt} is a central limit theorem only
in the case $\theta< \theta_*/2$. In the next corollary, we will deal with the case $\theta=\theta_*/2$.
		
		\begin{corollary}\label{cor:CLT3}
			Assume that there exists $\theta_*>0$ such that the conditions in Proposition \ref{maximum-weak} hold.
			If \eqref{H2} holds for $\theta=\theta_*/2$, then
			\begin{align*}
				\mathbf{L}(t^{1/4}e^{(\kappa(\theta)-\kappa(2\theta)/2)t}	\left(W_t(\theta)-W_\infty(\theta)\right)|\mathcal{F}_t) \overset{w}{\to}\mathbf{L}(\sigma_\theta\sqrt{CD_\infty}Z|\mathcal{F}_\infty),
			\end{align*}
			in probability with respect  to $\P$,
			where $C=\sqrt{\frac{2}{\pi \theta_*\kappa''(\theta_*)}},$ $Z\sim N(0,1)$  is independent of  $\mathcal{F}_\infty$.
		\end{corollary}
		{\bf Proof:}
		By  Proposition \ref{lemma-rate},  the limit
		\eqref{clt-1} holds in probability
		with $b(t)=t^{1/4}$ and $\Lambda_\theta=\sqrt{\frac{2}{\pi \theta_*\kappa''(\theta_*)}}D_\infty$.  By Proposition \ref{maximum-weak},  we have
		\begin{align*}
			t^{1/4}e^{-\kappa(\theta_*)t/2}e^{\theta_*M_t/2}=e^{\frac{\theta_*}{2}(M_t-m(t))} t^{-1/2}\overset{\P}{\to}0.
		\end{align*}
		Now the desired result follows from Theorem \ref{them-clt2}.
		\hfill$\Box$
		
		\bigskip
		
		For branching random walks, assuming  the second moment of $W_\infty(\theta)$ is finite,   Iksanov and Kabluchko \cite{IK16} established a central limit theorem  for  $W_\infty(\theta)-W_n(\theta)$ when $\kappa(2\theta)<2\kappa(\theta)$.
		Iksanov, Kolesko and Meiners \cite[Theorems 2.2 and 2.3]{IKM} are the counterparts of  Corollaries  \ref{them-clt} and \ref{cor:CLT3}
		for the fluctuations of the Biggins martingales with complex parameters for branching random walks.
		Note that since Corollary \ref{them-clt} is an almost sure result, we could not transfer the result about the branching random walks  to branching L\'evy processes
		directly using Lemma \ref{Croft-Kingman}.

		\subsubsection{Stable central limit theorems}

		First we give a  classical central limit theorem for some independent random variables. This is basically \cite[Section 25, Theorem 1]{GK}.
		
		\begin{lemma}\label{CLT}
			Assume that, for any $n\ge1$, $\{\xi_{n,k}, 1\le k\le k_n\}$ are independent random variables. Assume the following hold:
			\begin{itemize}
				\item[(1)] There exist $\alpha\in(1,2)$  and $c>0$ such that for any $x>0$,
				$$\lim_{n\to\infty} \sum_{k=1}^{k_n} P(\xi_{n,k}\le -x)=c x^{-\alpha},$$
				and
				$$\lim_{n\to\infty} \sum_{k=1}^{k_n} P(\xi_{n,k}> x)=0.$$
				\item[(2)] For any $\epsilon>0$,
				$$\lim_{\epsilon\to0} \lim_{n\to\infty}\sum_{k=1}^{k_n}E(\xi_{n,k}^2;|\xi_{n,k}|<\epsilon)=0.$$
				\item[(3)]
				$$\lim_{n\to\infty}\sum_{k=1}^{k_n}E(\xi_{n,k};|\xi_{n,k}|<1)=a_0.$$
			\end{itemize}
			Then
			$\sum_{k=1}^{k_n}\xi_{n,k}$ converges in distribution  and the characteristic function of the limit law is
			given by
			$$\exp\left\{i\lambda\left(a_0-\frac{c\alpha}{\alpha-1}\right)-c\Gamma(1-\alpha)e^{i\pi\alpha/2} \lambda^\alpha\right\},\quad \lambda>0.$$
		\end{lemma}
		{\bf Proof:} By \cite[Section 25, Theorem 1]{GK},
		$\sum_{k=1}^{k_n}\xi_{n,k}$ converges in distribution  and the characteristic function of the limit law is given by
		$$\exp\left\{ia_0\lambda+\int_{-\infty}^0 (e^{i\lambda x}-1-i\lambda x {\bf 1}_{\{x>-1\}})c\alpha |x|^{-1-\alpha}\,dx \right\}.$$
		Since $\alpha\in(1,2)$,  by  \cite[Lemma 14.11]{Sato}, we have that for $\lambda>0$,
		\begin{align*}
			&\int_{-\infty}^0 (e^{i\lambda x}-1-i\lambda x {\bf 1}_{\{x>-1\}}) |x|^{-1-\alpha}\,dx\\
			=&\int_{-\infty}^0 (e^{i\lambda x}-1-i\lambda x ) |x|^{-1-\alpha}\,dx-i\lambda\int_{-\infty}^{-1}|x|^{-\alpha} dx\\
			=&\lambda^{\alpha}\int_0^\infty (e^{-ir}-1+ir) r^{-1-\alpha} dr-i\frac{\lambda}{\alpha-1}\\
			=&\Gamma(-\alpha)e^{i\pi\alpha/2} \lambda^{\alpha}-i\frac{\lambda}{\alpha-1}.
		\end{align*}
		The desired result now follows from the fact $-\alpha\Gamma(-\alpha)=-\Gamma(1-\alpha).$
		\hfill$\Box$
		
		\bigskip

		\begin{theorem}\label{main-stableCLT}
			Let $\theta\in\Theta$, $p\in(1,2)$ and $L\in\mathbb{L}$.
			Assume that $p\theta\in\Theta_0$, $\kappa(p\theta)<p\kappa(\theta)$ and that
			$\P(\Xi_\theta>x)\sim  x^{-p}L(x)$.
			\begin{itemize}
				\item [(1)]
				If there exist a positive function  $b(t)$ and a random variable $\Lambda_{\theta,p}$ such that
				\begin{align}\label{con:bt}
					\lim_{t\to\infty} \mathcal{V}(t,M_t)=0, \quad \P\mbox{-a.s.}
				\end{align}
				and
				\begin{align}\label{con:bt2}
					\lim_{t\to\infty}\sum_{u\in\cL_t} \mathcal{V}(t,z_u(t))^pL(\mathcal{V}(t,z_u(t))^{-1})=\Lambda_{\theta,p}, \quad \P\mbox{-a.s.},
				\end{align}
				where $\mathcal{V}(t,x):=b(t)e^{-(\kappa(p\theta)/p) t}e^{\theta x}$,
				then
				\begin{align*}
					\mathbf{L}\left(b(t)(e^{(\kappa(\theta)-\kappa(p\theta)/p)t}( W_t(\theta)- W_\infty(\theta))|\mathcal{F}_t\right)\overset{w}{\to} \mathbf{L}(c_pU_p \Lambda_{\theta,p}^{1/p}|\mathcal{F}_\infty),
					\quad \P\mbox{-a.s., }
				\end{align*}
				where $c_p=\left(-\Gamma(1-p)\frac{\beta+p\kappa(\theta)-\varphi(p\theta)}{p\kappa(\theta)-\kappa(p\theta)}\right)^{1/p}>0$, $U_p$ is a stable random variable with characteristic function
				\begin{align}\label{Laplace-Up}
					E(e^{i\lambda U_p }) &=\exp\left\{
					e^{i\pi p/2}
					\lambda^p \right\},\quad \lambda>0.
				\end{align}
				Moreover, $U_p$ is independent of $\mathcal{F}_\infty$.
				
				\item[(2)]
				If \eqref{con:bt} and \eqref{con:bt2} hold in probability instead of almost surely, then
				\begin{align*}
					\mathbf{L}\left(b(t)(e^{(\kappa(\theta)-\kappa(p\theta)/p)t}( W_t(\theta)- W_\infty(\theta))|\mathcal{F}_t\right)\overset{w}{\to} \mathbf{L}(c_pU_p \Lambda_{\theta,p}^{1/p}|\mathcal{F}_\infty),
					\quad \mbox{in probability.}
				\end{align*}	
			\end{itemize}
		\end{theorem}
		
		\noindent{\bf Proof:}
		(1) Note that
		\begin{align*}
			b(t)(e^{\kappa(\theta)-\kappa(p\theta)/p)t}( W_t(\theta)- W_\infty(\theta))&=\sum_{u\in\cL_t}b(t)e^{-(\kappa(p\theta)/p) t}e^{\theta z_u(t)}(1- W^u_\infty(\theta))\\
			&:=\sum_{u\in\cL_t} \mathcal{V}(t,z_u(t))(1- W^u_\infty(\theta)).
		\end{align*}
		We now show that, conditioned on $\mathcal{F}_t$,
		the sequence $\{\mathcal{V}(t,z_u(t))(1- W^u_\infty(\theta)), u\in\mathcal{L}_t\}$ satisfies the conditions of Lemma \ref{CLT}.
		
		(i)
		Put $G_1(y);=\P(W_\infty(\theta)>y).$  Then for any $x>0$, we have
		\begin{align}\label{5.51}
			\sum_{u\in\cL_t}\P\left(\mathcal{V}(t,z_u(t))(1-W^u_\infty(\theta))<-x|\mathcal{F}_t\right)
			=\sum_{u\in\cL_t}G_1(1+x\mathcal{V}(t,z_u(t))^{-1} ).
		\end{align}
		By Theorem  \ref{Theorem:tailW}, $G_1(y)\sim  C_0y^{-p}L(y)$ as $y\to\infty$, where $C_0=\frac{\beta+p\kappa(\theta)-\varphi(p\theta)}{p\kappa(\theta)-\kappa(p\theta)}$. Thus for any $x>0$ and $\epsilon>0$, there exist $y_0=y_0(x,\epsilon)>0$ such that for any $y>y_0$,
		\begin{align*}
			(1-\epsilon)C_0x^{-p} &\le  \frac{G_1(1+xy)}{y^{-p}L(y)}
			\le (1+\epsilon) C_0x^{-p}.
		\end{align*}
		Since $\lim_{t\to\infty} b(t) e^{-(\kappa(p\theta)/p) t})e^{\theta M_t}\to 0$ $\P$-a.s., there exitss $\Omega_0$ with  $\P(\Omega_0)=1$ such that for any $\omega\in\Omega_0$, there exist $t_0=t_0(x,\epsilon,\omega)>0$ such that for all $t>t_0$,
		$$\inf_{u\in\cL_t}\mathcal{V}(t,z_u(t))^{-1} =b_t^{-1}e^{(\kappa(p\theta)/p) t} e^{-\theta M_t}>y_0.$$
		Thus, by \eqref{5.51}, for $\omega\in\Omega_0$ and $t>t_0$,
		\begin{align*}
			&(1-\epsilon)C_0\sum_{u\in\cL_t} \mathcal{V}(t,z_u(t))^pL(\mathcal{V}(t,z_u(t))^{-1}) x^{-p}\\
			\le& \sum_{u\in\cL_t}\P\left(
			\mathcal{V}(t,z_u(t))
			(1-W^u_\infty(\theta))<-x|\mathcal{F}_t\right)\\
			\le &(1+\epsilon)C_0\sum_{u\in\cL_t}  \mathcal{V}(t,z_u(t))^pL(\mathcal{V}(t,z_u(t))^{-1})  x^{-p}.
		\end{align*}
		Letting $t\to\infty$ and then $\epsilon\to0$, we get
		\begin{align}\label{CLT1}
			\lim_{t\to\infty}\sum_{u\in\cL_t}\P\left(\mathcal{V}(t,z_u(t)) (1-W^u_\infty(\theta))<-x|\mathcal{F}_t\right)
			=C_0\Lambda_{\theta,p} x^{-p},\quad \P\mbox{-a.s.}
		\end{align}

		(ii)
		Since $W_\infty(\theta)$ is non-negative and $\lim_{t\to\infty} \mathcal{V}(t,M_t)\to 0,$ $\P$-a.s.,
		we have
		\begin{align}\label{CLT2}
			\lim_{t\to\infty}\sum_{u\in\cL_t}\P\left(\mathcal{V}(t,z_u(t))(1-W^u_\infty(\theta))>x|\mathcal{F}_t\right)=0, \quad\mbox{ for any } x>0,  \P\mbox{ -a.s.}
		\end{align}

		(iii)
		For any $\epsilon>0$, we have
		\begin{align}\label{3.6'}
			&\sum_{u\in\cL_t}\E\left(\mathcal{V}(t,z_u(t))^2(1-W^u_\infty(\theta))^2{\bf 1}_{\{\mathcal{V}(t,z_u(t))|1-W^u_\infty(\theta)|<\epsilon\}}|\mathcal{F}_t\right)\nonumber\\
			=&\sum_{u\in\cL_t}\mathcal{V}(t,z_u(t))^2G_2(\epsilon \mathcal{V}(t,z_u(t))^{-1}),
		\end{align}
		where
		$G_2(y):=\E(|W_\infty(\theta)-1|^2;|W_\infty(\theta)-1|<y)$.
		Note that for $y>1$,
		\begin{align}\label{3.7'}
			\P(|W_\infty(\theta)-1|>y)= \P(W_\infty(\theta)>y+1)\sim C_0y^{-p}L(y),\quad y\to\infty.
		\end{align}
		Thus by Karamata's theorem, we have
		\begin{align*}
			G_2(y)\sim \frac{p}{2-p}C_0y^{2-p}L(y),\quad y\to\infty,
		\end{align*}
		Thus $y^{-2}T(\epsilon y)\sim \frac{p}{2-p}C_0\epsilon^{2-p}y^{-p}L(y)$ as $ y\to\infty.$
		Using an argument similar to that leading to \eqref{CLT1}, we get
		\begin{align}\label{CLT3}
			&\lim_{t\to\infty} \sum_{u\in\cL_t}\E\left(\mathcal{V}(t,z_u(t))^2(1-W^u_\infty(\theta))^2{\bf 1}_{\{\mathcal{V}(t,z_u(t))|1-W^u_\infty(\theta)|<\epsilon\}}|\mathcal{F}_t\right)\nonumber \\
			=&\frac{p}{2-p}C_0\epsilon^{2-p}\Lambda_{\theta,p}, \quad \P\mbox{-a.s.}
		\end{align}
		Note that the right hand side of \eqref{CLT3} tends to $0$ as $\epsilon\to0$.

		(iv) Note that
		\begin{align*}
			&\sum_{u\in\cL_t} \E\left(\mathcal{V}(t,z_u(t))(1-W^u_\infty(\theta));\mathcal{V}(t,z_u(t))|1-W^u_\infty(\theta)|\le1|\mathcal{F}_t\right)\\
			=&\sum_{u\in\cL_t}\mathcal{V}(t,z_u(t))G_3(\mathcal{V}(t,z_u(t))^{-1}),
		\end{align*}
		where
		$G_3(y)=\E(1-W_\infty(\theta);|W_\infty(\theta)-1|<y)=\E(W_\infty(\theta)-1;|W_\infty(\theta)-1|>y)$.
		Note that for $y>1$, by Karamata's theorem and \eqref{3.7'},
		\begin{align*}
			G_3(y) & = \E(|W_\infty(\theta)-1|;|W_\infty(\theta)-1|>y)\sim C_0\frac{p}{p-1}y^{1-p}L(y),\quad y\to\infty.
		\end{align*}
		Then, using an argument similar to that leading to \eqref{CLT1}, we have
		\begin{align}\label{CLT4}
			&\lim_{t\to\infty} 	\sum_{u\in\cL_t} \E\left(\mathcal{V}(t,z_u(t))(1-W^u_\infty(\theta));\mathcal{V}(t,z_u(t))|1-W^u_\infty(\theta)|\le1|\mathcal{F}_t\right)\nonumber\\
			=&C_0\frac{p}{p-1}\Lambda_{\theta,p},\quad \P\mbox{-a.s.}.
		\end{align}

		Combining \eqref{CLT1}, \eqref{CLT2}, \eqref{CLT3} and \eqref{CLT4}, we have by Lemma \ref{CLT} that for $\lambda>0$,
        \begin{align}\label{clt-513}
			&\lim_{t\to\infty}\E\left(\exp\Big\{i\lambda b(t)(e^{\kappa(\theta)-\kappa(p\theta)/p)t}( W_t(\theta)- W_\infty(\theta))\Big\}|\mathcal{F}_t\right)\nonumber\\
			=&\exp\{-C_0\Gamma(1-p)\Lambda_{\theta,p}
			{e^{i\pi p/2}}
			\lambda^p\},\quad \P\mbox{-a.s.}
\end{align}
		The desired result now follows.
		
		(2) If \eqref{con:bt} and \eqref{con:bt2} hold in probability instead of almost surely, then for any subsequence $\{t_n\}$, there exists a further subsequence $\{t_n'\}\subset\{t_n\}$, such that
		the limits in \eqref{con:bt} and \eqref{con:bt2} hold almost surely. Thus using the preceding analysis, we get the limits in \eqref{clt-513} for the subsequence $\{t_n'\}$ holds almost surely. This conclusion, in turn, implies that \eqref{clt-513} holds in probability as $t\to\infty$.
		
		\hfill$\Box$

		\begin{corollary}\label{cor:stable1}
			Let $\theta\in\Theta$, $p\in(1,2)$. Assume that $p\theta\in\Theta_0$ and $\kappa(p\theta)<p\kappa(\theta)$.
			If $\P(\Xi_\theta>x)\sim  lx^{-p}$ for some constant $l>0$, then
			\begin{description}
				\item[(1)]
				\begin{align*}
					\mathbf{L}\left(e^{(\kappa(\theta)-\kappa(p\theta)/p)t}( W_t(\theta)- W_\infty(\theta))|\mathcal{F}_t\right)\overset{w}{\to}
\mathbf{L}(c_p(lW_\infty(p\theta))^{1/p}U_p |\mathcal{F}_\infty),
\quad \P\mbox{-a.s.,}
				\end{align*}
where $U_p$ is a stable random variable with characteristic function given by \eqref{Laplace-Up} and  independent of $\mathcal{F}_\infty$.
				\item[(2)] 	 Assume that there exists $\theta_*>0$ such that the conditions in Proposition \ref{maximum-weak} hold.  If  $p\theta=\theta_*$, then
				\begin{align*}
					\mathbf{L}(t^{\frac{1}{2p}}e^{\left(\kappa(\theta)-\frac{\kappa(p\theta)}{p}\right)t}	(W_t(\theta)- W_\infty(\theta))|\mathcal{F}_t)\overset{w}{\to}
\mathbf{L}(c_p (l\sqrt{\frac{2}{\pi \theta_*\kappa''(\theta_*)}} \,D_\infty)^{1/p}U_p|\mathcal{F}_\infty)
				\end{align*}
				in probability with respect to $\P$,	
where $U_p$ is a stable random variable with characteristic function given by \eqref{Laplace-Up} and  independent of $\mathcal{F}_\infty$.
			\end{description}
		\end{corollary}
		{\bf Proof:} (1)  One can check that $b(t)=1$  satisfies the conditions in
 Theorem \ref{main-stableCLT} (1)
with $L(x)=l$.  In fact, by Lemma \ref{lemma:rightmost},
		$\frac{\kappa(p\theta)}{p\theta}t- M_t\to\infty$, $\P${-a.s.} as $t\to\infty$, which implies that $e^{-(\kappa(p\theta)/p)t}e^{\theta M_t}\to 0$, $\P$-a.s.
		Note that
		$$\sum_{u\in\mathcal{L}_t} le^{-\kappa(p\theta)t} e^{p\theta z_u(t)}=l W_{t}(p\theta)\to  l W_\infty(p\theta), \quad \P\mbox{-a.s.}$$
		Applying Theorem \ref{main-stableCLT} with $\Lambda_{\theta,p}=l W_\infty(p\theta)$, we get the desired result.
		
		(2) One can check that $b(t)=t^{\frac{1}{2p}}$  satisfies the conditions in
Theorem \ref{main-stableCLT} (2)
with $L(x)=l$.   	In fact, by Lemma \ref{maximum-weak},    $$b(t)e^{-\frac{\kappa(\theta p)}{p}t} e^{\theta M_t}=t^{-1/p}\exp\left\{\theta\left(M_t-m(t)\right)\right\}\overset{\P}{\to} 0.$$
		Note that
		$$\sum_{u\in\mathcal{L}_t}l b(t)^{p} e^{-\kappa(\theta p)t}e^{p\theta z_u(t)}=l \sqrt{t}W_{t}(p\theta)\overset{\P}{\to}l\sqrt{\frac{2}{\pi \theta_*\kappa''(\theta_*)}}D_\infty.$$
         Applying Theorem \ref{main-stableCLT} (2) with
		 $\Lambda_{\theta,p}=l\sqrt{\frac{2}{\pi \theta_*\kappa''(\theta_*)}}D_\infty$, we get the desired result.
		
		\hfill$\Box$

			When  $W_\infty(p\theta)$ is degenerate, Corollary \ref{cor:stable1} (1) is still valid with limit 0.
		Assume that there exists $\theta_*\in\Theta'$ such that $\theta_*\kappa'(\theta_*)=\kappa(\theta_*)$.
Note that
$\kappa(p\theta)<p\kappa(\theta)$ for some $p>1$ implies $\theta<\theta_*$.
	Note that  when $p\theta\ge \theta_*$,  $W_{\infty}(p\theta)=0$.
Consequently, Corollary \ref{cor:stable1} is a central limit theorem only
for the case $\theta< \theta_*/p$. In the next corollary, we will deal with the case $\theta=\theta_*/p$. The case $\theta>\theta_*/p$
will be dealt with in Proposition \ref{CLT-extremal case}.

		\begin{corollary}\label{cor:stable2}
			Let $\theta\in\Theta$, $p\in(1,2)$ and $l\in(0,\infty)$. Assume that
			$\P(\Xi_\theta>x)\sim  lx^{-p}\log_+x$, \eqref{condition-deriv} holds and there exists $\theta_*\in\Theta_0$ such that $\theta_*\kappa'(\theta_*)=\kappa(\theta_*)$.
          If $\theta=\theta_*/p$, then
			\begin{align*}
				e^{\left(\kappa(\theta)-\frac{\kappa(p\theta)}{p}\right)t}	(W_t(\theta)- W_\infty(\theta))\overset{d}{\to} c_p\left(\frac{l}{p}\right)^{1/p}
(\,D_\infty)^{1/p}U_p.
			\end{align*}	
		\end{corollary}
		{\bf Proof:}  Put $L(x):=l\log_+(x)$.  By Lemma \ref{lemma:rightmost},   we have  $e^{-(\kappa(p\theta)/p)t}e^{\theta M_t}\to 0$, $\P$-a.s.
		Thus, for $t$ large enough  such that $e^{(\kappa(p\theta)/p )t} e^{-\theta M_t}>1$ almost surely, we have
		\begin{align*}
			\sum_{u\in\mathcal{L}_t} e^{-\kappa(p\theta)t} e^{p\theta z_u(t)}\log_+(e^{(\kappa(p\theta)/p )t} e^{-\theta z_u(t)})&=\sum_{u\in\mathcal{L}_t} e^{-\kappa(p\theta)t}e^{p\theta z_u(t)}\left(\frac{\kappa(\theta p)}{p}t-\theta z_u(t)\right)\\
			&=\frac{1}{p} D_t\to \frac{1}{p} D_\infty,\quad \P\mbox{-a.s.}.
		\end{align*}
		Applying Theorem \ref{main-stableCLT} with $b(t)=1$ and $\Lambda_{\theta,p}=\frac{l}{p} D_\infty$, we get the desired result.
		\hfill$\Box$
		
		\bigskip

		The counterpart of Corollary  \ref{cor:stable1} for Galton-Watson processes was proved in \cite{Heyde71}.
		An analogue of \ref{cor:stable1} for  branching random walks can be found in \cite{IKM, IKM20}. 	\cite[Theorem 2.9]{IKM20} is on stable fluctuation of Biggins' martingale with complex parameters for branching random walks. Note that the conditions in \cite[Theorem 2.9]{IKM20} are not applicable to real parameters.

\subsubsection{The case $\theta>\theta_*/p$}

	In the case $p\theta>\theta_*$, we study the fluctuations of $W_t(\theta)-W_\infty(\theta)$ via the extremal process of $X$.
A counterpart for branching random walks can be found in \cite[Theorem 2.1]{ IKM20} and an analogue for super-Brownian motions can be found in \cite{Yang}.

		Suppose there exists $\theta_*>0$ such that $\kappa''(\theta_*)<\infty$ and   $\theta_*\kappa'(\theta_*)=\kappa(\theta_*)$.   Additionally,  we assume that $\mathcal{P}$ is non-lattice and  \eqref{H3} holds.  Then
		$\mathcal{E}_t:=\sum_{u\in\mathcal{L}_t}\delta_{z_u(t)-m(t)}$ converges weakly to  $\mathcal{E}_\infty$. Let $e_1\ge e_2\ge e_3\ge \cdots$ be the atoms of $\mathcal{E}_\infty$.
		
		Let $\theta\in\Theta$ and assume that $\theta\in(\theta_*/2,\theta_*)$. Assume that there exists $p\in(\theta_*/\theta,2]$
		such that \eqref{Hp} holds.
		Applying \cite[Theorem 2.5]{IKM20} to $\{-X_{nt},n=0,1,\cdots\}$ and noting that $\mu_n$ in \cite{IKM20} is equal to $-s^*\mathcal{E}_{nt}+t^{\frac{3}{2}}$,  we get $\mu_\infty=-s_*\mathcal{E}_\infty+\frac{3}{2}\log t$. Let $W^{(k)}_\infty(\theta), k\ge1,$ be i.i.d. random variables with the same law as  $W_\infty(\theta)$. By \cite[Theorem 2.5]{IKM20}, we know that
		$$\lim_{n\to\infty}\sum_{k=1}^ne^{\theta e_k}(W^{(k)}_\infty(\theta)-1)$$
		converges almost surely to a non-degenerate random variable which we  denote by $\tilde{\mathcal{E}}_\theta$. Note that $X_{ext}$ in \cite{IKM20} is equal to $t^{-\frac{3\theta}{2\theta_*}}\tilde{\mathcal{E}}_\theta$.

		Applying \cite[Theorem 2.5]{IKM20} to $\{-X_{nt},n=0,1,\cdots\}$, we get that as $n\to\infty$,
		\begin{align*}
			e^{-\theta m(nt)}e^{\kappa(\theta)nt}(W_\infty(\theta)-W_{nt}(\theta))\overset{d}{\to}\tilde{\mathcal{E}}_\theta.
		\end{align*}
		Then using Lemma \ref{Croft-Kingman}, we can obtain the following result.
		\begin{proposition}\label{CLT-extremal case}
			Suppose there exists $\theta_*>0$ such that $\kappa''(\theta_*)<\infty$ and   $\theta_*\kappa'(\theta_*)=\kappa(\theta_*)$.   Assume that $\mathcal{P}$ is non-lattice and  \eqref{H3} holds.  Let
			$\theta\in(\theta_*/2,\theta_*)$.
			If there exists $p\in(\theta_*/\theta,2]$ such that \eqref{Hp} holds, then as $t\to\infty$,
			\begin{align*}
				e^{-\theta m(t)}e^{\kappa(\theta)t}(W_\infty(\theta)-W_t(\theta))\overset{d}{\to}\tilde{\mathcal{E}}_\theta.
			\end{align*}
		\end{proposition}

		\section{Appendix}
		
We first give the proof of Lemma \ref{f-moment}.

		\noindent{\bf Proof of Lemma \ref{f-moment}.}
(1) We first prove the sufficiency in the case $\bE(N)<\infty$. Let $\tau_m$ be the $m$-th branching time. Since $\bE(N)<\infty$, on the survival event, for any $m\ge1$, we have $\tau_m<\infty$ and $\tau_m\uparrow \infty$ . On the extinction event, the total number $T$ of branching events is finite and $\tau_m=\infty$ for all $m>T$. Note that $\tau_1$  is exponentially distributed with parameter $\beta$. Assume that \eqref{condition-f} holds. Let $$U_m(t,x)=\E_{\delta_x}(\langle f,X_t\rangle, t<\tau_m).$$
By \cite[Proposition 25.4]{Sato}, we have
$$
U_1(t,x)=e^{-\beta t}\bE_x(f(\xi_t))<\infty.
$$
By the Markov property and the branching property, we have that
\begin{align}\label{f-2}
U_{m+1}(t,x)&=\E_{\delta_x}(\langle f,X_t\rangle, t<\tau_1)+\E_{\delta_x}(\langle f,X_t\rangle, \tau_1\le t<\tau_{m+1})\nonumber\\
&=U_1(t,x)+\bE_x\int_0^t\beta e^{-\beta s} \E_{\xi_s+\mathcal{P}}(\langle f,X_{t-s}\rangle, t<\tau_m)\,ds\nonumber\\
&\le U_1(t,x)+\bE_x\int_0^t\beta e^{-\beta s} \sum_{i=1}^NU_m(t-s,\xi_s+S_i)\,ds.
\end{align}
Note that there exist $c_2>0$ such that $|x+y|^r\le c_2(|x)^r+|y|^r)$ and $((x+y)\vee0)^r\le c_2((x\vee 0)^r+(y\vee0)^r)$, which implies that
\begin{align}\label{f-3}
f(x+y)\le c_2(e^{\theta y}f(x)+e^{\theta x}f(y)).
\end{align}
Thus, by the spatial homogeneity of $X$, we have
\begin{align}\label{f-1}
U_m(t,x+y)&=\E_{\delta_x}(\langle f(y+\cdot),X_t\rangle,t<\tau_m)\nonumber\\
&\le c_2f(y)\E_{\delta_x}(\langle e_\theta,X_t\rangle,t<\tau_m)+c_2e^{\theta y}\E_{\delta_x}(\langle f,X_t\rangle,t<\tau_m)\nonumber\\
&\le c_2e^{\kappa(\theta)t}e^{\theta x}f(y)+c_2e^{\theta y}U_m(t,x).
\end{align}
Applying  \eqref{f-1} with $x=\xi_s$ and $y=S_i$,  we get
\begin{align}
&\bE_x\int_0^t\beta e^{-\beta s} \sum_{i=1}^NU_m(t-s,\xi_s+S_i)\,ds\nonumber\\
&\le c_2\bE_x \int_0^t \beta e^{-\beta s}e^{\kappa(\theta)(t-s)}e^{\theta \xi_s}\sum_{i=1}^N f(S_i)\,ds
+c_2\bE_x \int_0^t \beta e^{-\beta s}U_m(t-s,\xi_s)\langle e_\theta,\mathcal{P}\rangle\,ds\nonumber\\
&=c_2\bE\Big(\sum_{i=1}^N f(S_i)\Big)\bE_x \int_0^t \beta e^{-\beta s}e^{\kappa(\theta)(t-s)}e^{\theta \xi_s}\,ds+c_2\chi(\theta)\bE_x \int_0^t \beta e^{-\beta s}U_m(t-s,\xi_s)\,ds.
\end{align}
Plugging this into \eqref{f-2}, we get that
\begin{align}
U_{m+1}(t,x)\le &e^{-\beta t}\bE_x(f(\xi_t))+c_2\bE\Big(\sum_{i=1}^N f(S_i)\Big)\bE_x \int_0^t \beta e^{-\beta s}e^{\kappa(\theta)(t-s)}e^{\theta \xi_s}\,ds\nonumber\\
&\quad +c_2\chi(\theta)\bE_x \int_0^t \beta e^{-\beta s}U_m(t-s,\xi_s)\,ds.
\end{align}
Using this and induction on $m$, we can show that
$$
U_m(t,x)\le e^{\gamma t}\bE_x f(\xi_t)+c_2\bE\Big(\sum_{i=1}^N f(S_i)\Big) \int_0^t \beta e^{\gamma s}e^{\kappa(\theta)(t-s)}e^{\varphi(\theta)s}\,ds\,e^{\theta x},
$$
where $\gamma=\beta(c_2\chi(\theta)-1).$
Thus we have that for all $t>0$,
\begin{align}\label{111}
\E(\langle f,X_t\rangle)&=\sup_{m}U_m(t,0)\nonumber\\
&\le e^{\gamma t}\bE f(\xi_t)+c_2\bE\Big(\sum_{i=1}^N f(S_i)\Big) \int_0^t \beta e^{\gamma s}e^{\kappa(\theta)(t-s)}e^{\varphi(\theta)s}\,ds<\infty.
\end{align}
		
(2) We next  prove the sufficiency in the case $\bE(N)=\infty$.
Using \eqref{111} for $X^{(n)}$, we get that
\begin{align}
\E(\langle f,X_t^{(n)}\rangle)&\le e^{\gamma_n t}\bE f(\xi_t)+c_2\bE\Big(\sum_{i=1}^N {\bf 1}_{\{S_i>-n\}}f(S_i)\Big) \int_0^t \beta e^{\gamma_n s}e^{\kappa^{(n)}(\theta)(t-s)}e^{\varphi(\theta) s}\,ds\\
&\le e^{\gamma t}\bE f(\xi_t)+c_2\bE\Big(\sum_{i=1}^N f(S_i)\Big) \int_0^t \beta e^{\gamma s}e^{\kappa(\theta)(t-s)}e^{\varphi(\theta)s}\,ds<\infty,
\end{align}
where $\gamma_n=\beta(c_2 \bE(\langle e_\theta,\mathcal{P}^{(n)}\rangle-1).$
Thus
$$
\E(\langle f,X_t\rangle)=\lim_{n\to\infty}\E(\langle f,X_t^{(n)}\rangle)<\infty.
$$

(3) We now prove the necessity. Assume that $\E(\langle f,X_t\rangle)<\infty $ for some $t>0$. Note that for any $x$,
\begin{align}
\E_{\delta_x}(\langle f,X_t\rangle)&=\E(\langle f(x+\cdot),X_t\rangle)
\le c_2f(x)\E(\langle e_\theta,X_t\rangle)+c_2e^{\theta x}\E(\langle f,X_t\rangle)\nonumber\\
&=c_2e^{\kappa(\theta)t}f(x)+c_2\E(\langle f,X_t\rangle)e^{\theta x}<\infty.
\end{align}
By the  Markov property and the branching property, we get that
\begin{align*}
\infty>\E(\langle f,X_t\rangle)&=\E(\langle f,X_t\rangle;t<\tau_1)+\E(\langle f,X_t\rangle;t\ge \tau_1)\\
&=e^{-\beta t}\bE(f(\xi_t))+\bE \int_0^t \beta e^{-\beta s}\sum_{i=1}^N \E_{\delta_{\xi_s+S_i}}(\langle f,X_{t-s}\rangle)\,ds.
\end{align*}
Thus we have  $\bE(f(\xi_t))<\infty$, which implies $\int_{|y|>1}f(y)\,n(dy)<\infty$ by
\cite[Proposition 25.4]{Sato}.
From the display above we also get that there exist $x$ and $r>0$ such that
$$
\bE\left(\sum_{i=1}^N \E_{\delta_{x+S_i}}(\langle f,X_{r}\rangle)\right)<\infty.
$$
		
If $f(x)=|x|^re^{\theta x}$, using the fact that $|y|^r\le c_2(|x+y|^r+|x|^r)$, we have
$$
f(x+y)\ge c_2^{-1}f(y)e^{\theta x}-f(x)e^{\theta y}.
$$
Thus for any $y$
\begin{align*}
\E_{\delta_{x+y}}(\langle f,X_r\rangle)=\E_{\delta_{x}}(\langle f(y+\cdot),X_r\rangle)
\ge c_2^{-1}e^{\kappa(\theta)r}e^{\theta x}f(y)-\E_{\delta_x}(\langle f,X_r\rangle)e^{\theta y}.
\end{align*}
It follows that
\begin{align*}
\infty>\bE\left(\sum_{i=1}^N \E_{\delta_{x+S_i}}(\langle f,X_{r}\rangle)\right)\ge c_2^{-1}e^{\kappa(\theta)r}e^{\theta x}\bE\left(\sum_{i=1}^N f(S_i)\right) -\E_{\delta_x}(\langle f,X_r\rangle)\chi(\theta),
\end{align*}
which implies that $\bE\left(\sum_{i=1}^N f(S_i)\right)<\infty.$

If  $f(x)=(x\vee 0)^re^{\theta x}$, then for any $a>0$,
$$
f(x+y)\ge {\bf 1}_{\{x>-a\}} e^{\theta x}{\bf 1}_{\{y>a\}} (y-a)^r e^{\theta y}.
$$
Thus for any $a>0$,
\begin{align}\label{3.1.3}
\infty>\bE\left(\sum_{i=1}^N \E_{\delta_{x+S_i}}(\langle f,X_{r}\rangle)\right)\ge \bE\left(\sum_{i=1}^N {\bf 1}_{\{S_i>a\}}(S_i-a)^re^{\theta S_i}\right) \E_{\delta_x}(\langle {\bf 1}_{(-a,\infty)} e_\theta,X_r\rangle).
\end{align}
Since $ \E_{\delta_x}(\langle e_\theta,X_r\rangle)=e^{\theta x}e^{\kappa(\theta)r}>0,$
 we can choose $a$ sufficiently large
such that $ \E_{\delta_x}(\langle {\bf 1}_{(-a,\infty)} e_\theta,X_r\rangle)>0$. Hence   $\bE\left(\sum_{i=1}^N {\bf 1}_{\{S_i>a\}}(S_i-a)^re^{\theta S_i}\right)<\infty$.
It is clear that $(y\vee 0)^r\le {\bf 1}_{\{y<a\}} a^r+{\bf 1}_{\{y\ge a\}}y^r\le a^r+ c_2a^r+c_2{\bf 1}_{\{y\ge a\}}(y-a)^r$. Thus we can get
$\bE\left(\sum_{i=1}^N f(S_i)\right)<\infty.$ Hence the proof is complete.
\hfill $\Box$

		\begin{lemma}\label{lemma:A1}
			Assume that $X,Y$ are two independent positive random variables . If $E(X^{p+\delta})<\infty$ and $P(Y>y)\sim y^{-p}L(y)$ as $y\to \infty$, where $p>\delta>0$ and $L$ is slowly varying at $\infty$, then
			$$P(XY>y)\sim E(X^p)y^{-p}L(y).$$
		\end{lemma}
		{\bf Proof:} Without loose of generality,  we assume that $x^{\delta}L(x)$ is increasing and $x^{-\delta}L(x)$ is decreasing.
		Put $g(y):=P(Y>y)$. Note that $P(XY>y)=E g(yX^{-1})$.  Since $g(y)\sim y^{-p}L(y)$, there exists $y_0>0$ such that $g(y)\le 2y^{-p}L(y)$ for all $y\ge y_0$.
		Thus
		\begin{align*}
			g(yX^{-1})&\le {\bf 1}_{\{y_0X>y\}}+{\bf 1}_{\{y_0X\le y\}} 2y^{-p}(X^p L(yX^{-1}))\\
			&\le (y_0X)^{p+\delta}y^{-p-\delta}+2 y^{-p}L(y) (X^{p+\delta}\vee X^{p-\delta}),
		\end{align*}
		where in the last inequality, we used the fact
		$$L(yX^{-1})/L(y)\le X^\delta\vee X^{-\delta}.$$
		Now by the dominated convergence theorem, we have
		$$\lim_{y\to\infty}\frac{E g(yX^{-1})}{y^{-p}L(y)}=E(X^p).$$

		\hfill$\Box$
		
		\begin{lemma}\label{lemma:A2}
			Assume that $X, X_i,i\ge 1,$ are i.i.d.  non-negative random variables,
			and  $N$ is a positive integer-valued random variable independent of $\{X_i\}$. Let  $p\in(1,2)$ and $L$ be slowly varying at $\infty$. If,  as $x\to\infty$
			$$P(N>x)\sim c_1x^{-p}L(x), \quad P(X>x)\sim c_2x^{-p}L(x),$$
			for some  $c_1,c_2\ge0$, then
			$$P\left(\sum_{i=1}^N X_i>x\right)\sim \Big(c_1(EX)^p+c_2EN\Big)x^{-p}L(x), \quad x\to\infty.$$
		\end{lemma}
		{\bf Proof} Put $Z:=\sum_{i=1}^N X_i$, $J_X(\lambda):=-\log E(e^{-\lambda X})$ and $\Phi_N(\lambda):=E(e^{-\lambda N})$.
		Since $E(Z)=E(N)E(X)$, we have
		\begin{align}\label{A2}
			&E(e^{-\lambda Z} -1+\lambda Z)=\Phi_N(J_X(\lambda))-1+\lambda E(N)E(X)\nonumber\\
			&=\Phi_N(J_X(\lambda))-1+J_X(\lambda)E(N)+(\lambda E(X)-J_X(\lambda))E(N).
		\end{align}
		Since $P(N>x)\sim c_1x^{-p}L(x)$, we have
		$$\Phi(\lambda)-1+\lambda E(N)\sim \frac{\Gamma(2-p)}{p-1}c_1 \lambda^pL(1/\lambda).$$
		Note that $J_N(\lambda)\sim \lambda E(X)$ as $\lambda\to0$. It follows that
		$$\Phi_N(J_N(\lambda))-1+J_N(\lambda)E(N)\sim \frac{\Gamma(2-p)}{p-1}c_1 (E(X))^p\lambda^pL(1/\lambda), \quad \lambda\to0.$$
		Since $P(X>x)\sim c_2x^{-p}L(x)$, wew have
		$$
		e^{-J_N(\lambda)}-1+\lambda E(X)\sim \frac{\Gamma(2-p)}{p-1}c_2 \lambda^pL(1/\lambda).$$
		Noting that  as $\lambda\to0$
		$$e^{-J_N(\lambda)}-1+J_N(\lambda)\sim \frac{1}{2} J_N(\lambda)^2\sim \frac{1}{2}( E(X))^2\lambda^2,$$ we have
		$$\lambda E(X)-J_X(\lambda)=e^{-J_N(\lambda)}-1+\lambda E(X)-(e^{-J_X(\lambda)}-1+J_X(\lambda))\sim \frac{\Gamma(2-p)}{p-1}c_2 \lambda^pL(1/\lambda).$$
		Plugging this into \eqref{A2}, we get
		$$E(e^{-\lambda Z} -1+\lambda Z)\sim \frac{\Gamma(2-p)}{p-1}\Big(c_1(EX)^p+c_2EN\Big)\lambda^pL(1/\lambda).$$
		Tthe proof is complete.
		\hfill$\Box$

		\begin{lemma}\label{inequality}
			Let $n\ge 0$ and  $p\in[n,n+1]$.
			If $L$ is slowly varying at $\infty$ such that $x^{p-n}L(x)$ is increasing and $x^{-(n+1-p)}L(x)$ is decreasing, then for any $x,y>0$,
			$$(xy)^n\wedge (xy)^{n+1}\le x^pL(x) y^p L^{-1}(y^{-1}).$$
		\end{lemma}
		{\bf Proof}  Since $x^{p-n}L(x)$ is increasing, we have for $xy>1$,
		$$x^{p-n}L(x)\ge y^{n-p}L(y^{-1}),$$
		which implies that
		$$(xy)^n\le x^pL(x) y^p L^{-1}(y^{-1}).$$
		Since $x^{-(n+1-p)}L(x)$ is decreasing, we have for $xy\le 1$,
		$$x^{-(n+1-p)}L(x)\ge y^{n+1-p}L(y^{-1}),$$
		which implies that
		$$(xy)^{n+1}\le x^pL(x) y^p L^{-1}(y^{-1}).$$
		Now the proof is complete.
		\hfill$\Box$

\section*{Acknowledgments}
We thank Bastien Mallein for informing us the paper Dadoun \cite{Dadoun}.
Yan-Xia Ren  is supported by   NSFC (Grant No. 12231002) and the Fundamental Research Funds for Central Universities, Peking University LMEQF.
Renming Song's research  was supported in part by a grant from the Simons
Foundation (\#960480, Renming Song). Rui Zhang is  supported by  NSFC (Grant No.  12271374, 12371143),  and Academy for Multidisciplinary Studies, Capital Normal University.

		\end{doublespace}


\begin{thebibliography}{99}
		
		\bibitem{Aid13} E. A\"idekon. Convergence in law of the minimum of a branching random walk. \emph{Ann. Probab.}, {\bf 41}(2013), 1362--1426.
		
		\bibitem{ASh14} E. A\"idekon  and Z. Shi. The Seneta-Heyde scaling for the branching random walk. \emph{Ann. Probab.}, {\bf 42}(3) (2014) , 959--993.
		
		\bibitem{Athreya}  K. B. Athreya and  P. E. Ney. \emph{Branching Processes}. Dover, 2004.



		\bibitem{BM} J. Bertoin and B. Mallein. Infinitely ramified point measures and branching L\'evy processes. \emph{Ann. Probab.},
		{\bf 47} (2019), 1619--1652.
		
		\bibitem{BM18} J. Bertoin and B. Mallein.  Biggins' martingale convergence for branching L\'evy processes. \emph{Electron. Commun. Probab.} {\bf 23} (2018),  1--12.
		
		\bibitem{Biggins77} J. D. Biggins.  Martingale convergence in the branching random walk. \emph{ J. Appl. Probability}, {\bf 14} (1977), 25--37.
		
		\bibitem{Biggins92} J. D. Biggins.  Uniform convergence of martingales in the branching random walk. \emph{Ann. Probab.} {\bf 20} (1992), 137--151.
		
		
		\bibitem{BK97}  J. D. Biggins and  A. E. Kyprianou.  Seneta-Heyde norming in the branching random walk. \emph{Ann. Prob.} {\bf  25}(1997), 337--360.
		
		\bibitem{BD74} N. H. Bingham and R. A. Doney.  Asymptotic Properties of Supercritical Branching Processes I: The Galton-Watson Process. \emph{Adv. Appl. Probab.}, {\bf 6}(1974),  711--731.
		
		
		\bibitem{Bingham} N. H. Bingham, C. M. Goldie and J. L. Teugels. \emph{ Regular Variation}. Cambridge Univ. Press, Cambridge, 1978.
		
		\bibitem{Ch15} X. Chen. A necessary and sufficient condition for the nontrivial limit of the derivative martingale in a branching random walk. \emph {Adv.  Appl. Probab.}, {\bf 47}(3) (2015), 741--760.


\bibitem{Dadoun} B. Dadoun. Asymptotics of self-similar growth-fragmentation processes. \emph{  Electron. J. Probab.} {\bf 22} (2017), no. 27, 1--30.

		
		\bibitem{Durrett} R. Durrett. {\it Probability: Theory and Examples.} Fourth edition. Cambridge University Press, Cambridge, 2010.
		
		\bibitem{GBTG}G. Fa\"{y}, B. Gonz\'{a}lez-Ar\'{e}valo, T. Mikosch, G. Samorodnitsky. Modeling teletraffic arrivals
		by a Poisson cluster process. \emph{Queueing Syst.} {\bf 54 }(2) (2006), 121--140.
		
		\bibitem{GK} B. V. Gnedenko  and A. N. Kolmogorov. \emph{Limit distributions for sums of independent random variables}.  Addison-Wesley publishing company.


\bibitem{Hering}
H. Hering. Critical Markov branching processes with general set of types. \emph{Trans. Amer. Math. Soc.}, {\bf 160} (1971), 185--202.

		
		\bibitem{Heyde71} C. C. Heyde. Some central limit analogues for supercritical Galton-Watson processes. \emph{ J. Appl. Probab.} {\bf 8} (1971) 52--59.
		
		
		\bibitem{HSh09} Y. Hu and Z. Shi. Minimal position and critical martingale convergence in branching random walks, and directed polymers on disordered trees.
		\emph{Ann. Prob.} {\bf 37} (2009),  742-789.
		
		
		\bibitem{IK16} A. Iksanov and Z. Kabluchko. A central limit theorem and a law of the iterated logarithm for the Biggins martingale of the supercritical branching random walk. \emph{J. Appl. Probab}, {\bf 53}(2016), 1178--1192.
		
		\bibitem{IKM} A. Iksanov, K. Kolesko and M. Meiners. Stable-like fluctuations of Biggins' martingales. \emph{Stochastic Process. Appl.} {\bf 129}(2019), 4480--4499.
		
		\bibitem{IKM20} A. Iksanov, K. Kolesko and M. Meiners. Fluctuations of Biggins' martingales at complex parameters. \emph{Ann. l'Institut Henri Poincar\'e--Probab.  Statist.}, {\bf 56}(2020), 2445--2479.
		
		\bibitem{IB19}A. Iksanov and B. Mallein. A result on power moments of L\'evy-type perpetuities and its application to the $L_p$-convergence of Biggins' martingales in branching L\'evy processes. \emph{ALEA Lat. Am. J. Probab. Math. Stat.}, {\bf 16}(1) (2019), 315--331.
		
		
		\bibitem{King63} J. F. C. Kingman. Ergodic properties of continuous-time Markov processes and their discrete skeletons. \emph{Proc.
			Lond. Math. Soc.} {\bf 13}(1963),  593--604.
	
		
		\bibitem{Kyp99} A. E. Kyprianou. A note on branching L\'evy processes. \emph{Stochastic
		Process. Appl.}, {\bf 82}(1)(1999): 1--14.
		
		\bibitem{Liu2000} Q. Liu. On generalized multiplicative cascades. \emph{Stoch. Proc. Appl.},  {\bf 86} (2000), 263--286.
		
		\bibitem{Madaule} T. Madaule. Convergence in law for the Branching random walk seen from its tip. \emph{J. Theor. Probab.} {\bf 30}(2017), 27-63.
		
		\bibitem{MSh23} B. Mallein and  Q. Shi. A necessary and sufficient condition for the convergence of the derivative martingale in a branching L\'evy process.
		\emph{Bernoulli} {\bf 29}  (2023),  597-624.
		
		\bibitem{Sato} K.-I. Sato. \emph{L\'evy processes and infinitely divisible distributions}, Cambridge University Press,  2013.
		
		
		\bibitem{Yang} T. Yang.  Fluctuations of the additive martingales related to super-Brownian motion.
		\emph{Acta Math. Appl. Sin. Engl. Ser.}, https://doi.org/10.1007/s10255-025-0031-8.
		
		\bibitem{YR11} T. Yang and Y.-X. Ren. Limit theorem for derivative martingale at criticality w.r.t. branching Brownian motion. \emph{Statist. Prob. Lett.} {\bf 81} (2011), 195-200.
		
		\end{thebibliography}
\end{document}